\documentclass[11pt,a4paper,twoside]{article}

\usepackage{titling}
\usepackage{titletoc}
\usepackage{url}

\usepackage[utf8]{inputenc}
\usepackage[T1]{fontenc}
\usepackage{amsmath}
\usepackage{amsfonts}
\usepackage{amssymb}
\usepackage{amsthm}

\usepackage{mathrsfs}
\usepackage{mathtools}
\usepackage{stmaryrd}
\usepackage[all]{xy}
\usepackage{makeidx}
\usepackage{natbib}
\bibpunct{(}{)}{,}{a}{}{,}
\defcitealias{stacks-project}{Stacks Project}
\defcitealias{fstacks-project}{Stacks Project}
\usepackage[french,english]{babel}\frenchsetup{GlobalLayoutFrench=false}
\usepackage{xspace}
\usepackage{float}

\usepackage{etoolbox}

\usepackage{stixgras}
\RequirePackage{amsbsy}

\makeatletter\def\th@plain{\slshape}\makeatother
\makeatletter\patchcmd{\th@remark}{\itshape}{\slshape}{}{}\makeatother

\newcounter{bidon}
\newcommand{\rdb}{\refstepcounter{bidon}}

\newcommand\Section[1]{
\rdb\addcontentsline{toc}{section}{#1} \section*{#1}}

\newcommand\Subsubsection[1]{
\rdb\addcontentsline{toc}{subsubsection}{#1} 
\subsubsection*{#1}}

\renewcommand\paragraph[1]{
\rdb
\addcontentsline{toc}{paragraph}{#1} 
\medskip \noindent $\bullet$ \textbf{#1}}

\newcommand \sibrouillon[1]{}

\newcommand \hum[1] {\sibrouillon{\noindent {\sf hum: #1}}}

\usepackage[bookmarksopen=false,breaklinks=true,%
      backref=page,pagebackref=true,plainpages=false,%
      hyperindex=true,pdfstartview=FitH,colorlinks=true,%
      pdfpagelabels=true,linkcolor=blue,%
      citecolor=red,urlcolor=red,
      ]%
   {hyperref}

\begin{document}

\thispagestyle{empty}
~ 
\vspace{3cm}

\noindent In this file you find the English version starting on the page  numbered \pageref{beginenglish}.

\medskip \noindent  {\Large \bf Multivariate Hensel Lemma for ultrametric fields}

\bigskip \noindent  
Then the French version begins on the page numbered \pageref{beginfrench}.

\medskip\noindent   {\Large \bf Lemme de Hensel multivarié pour les corps ultramétriques} 

\smallskip \noindent Le lecteur ou la lectrice sera sans doute surprise de l'alternance des sexes ainsi que de l'orthographe du mot 'corolaire', avec d'autres innovations auxquelles elle n'est pas habituée. En fait, nous avons essayé de suivre au plus près les préconisations de l'orthographe nouvelle recommandée, telle qu'elle est enseignée aujourd'hui dans les écoles en France.  

\bigskip\noindent   {\large \bf Authors}  

\noindent M.-E. Alonso, {Universidad Complutense, Madrid, Espa\~na..\\
email: {\tt mariemi@mat.ucm.es}}

\smallskip \noindent Henri Lombardi
Université de Franche-Comté, Laboratoire de mathématiques de Besançon, UMR CNRS 6623, 16 route de Gray, 25000
Besançon, France. \\
email: {\tt henri.lombardi@univ-fcomte.fr}

\smallskip \noindent Stefan Neuwirth
Université de Franche-Comté, Laboratoire de mathématiques de Besançon, UMR CNRS 6623, 16 route de Gray, 25000
Besançon, France.. \\
email: {\tt stefan.neuwirth@univ-fcomte.fr}

\newpage
\thispagestyle{empty}

~

\pagestyle{headings}
\patchcmd{\sectionmark}{\MakeUppercase}{}{}{}
\setcounter{page}{0}\renewcommand\thepage{E\arabic{page}}

\begingroup
\startcontents[english]

\def\proofname{\textsl{Proof}}

\newcommand{\di}{\,\vert\,}
\newcommand{\ndi}{\nmid}

\newcommand \nzr {\neq 0}
\newcommand{\Vr}{\mathrm{Vr}}
\newcommand{\Rn}{\mathrm{Rn}}
\newcommand{\Nvr}{\mathrm{Nvr}}
\newcommand{\Nrn}{\mathrm{Nrn}}
\newcommand{\U}{\mathrm{U}}

\newcommand{\Gati}{\Gat\cup\so\infty }

\newcommand \Rzero {R_{=0}}
\newcommand \Rnz {R_{\nzr}}
\newcommand \Rvr {R_{\Vr}}
\newcommand \Ru {R_{\U}}
\newcommand \Rrn {R_{\Rn}}
\newcommand \Rnrn {R_{\Nrn}}
\newcommand \Rnvr {R_{\Nvr}}

\newcommand \zg {{\ZZ[G]}}
\newcommand \Izero {{\cI}_{=0}}
\newcommand \Mnz {{\cM}_{\nzr}}
\newcommand \Mpos {{\cM}_{\Pos}}
\newcommand \Cnng {{\cC}_{\Nng}}
\newcommand \Irn {{\cI}_{\Rn}}
\newcommand \Vvr {{\cV}_{\Vr}}
\newcommand \Mu {{\cM}_{\U}}
\newcommand \abg {{\rm Ab}(G)}
\newcommand \Hzero {{\cH}_{=0}}


\newcommand\Exists{\boldsymbol{\stixexists}}
\newcommand\Forall{\boldsymbol{\stixforall}}
\newcommand\VDash{\boldsymbol{\stixvdash}}
\newcommand\Land{\boldsymbol{\stixwedge}}
\newcommand\Lor{\boldsymbol{\stixvee}}
\newcommand\Top{\boldsymbol{\stixtop}}
\newcommand\Bot{\boldsymbol{\stixbot}}
\newcommand\bigLand{\boldsymbol{\stixbigwedge}}
\newcommand\bigLor{\boldsymbol{\stixbigvee}}

\newcommand\dotminus{\stixdotminus}

\newcommand \vdg{\VDash}
\newcommand \Vd {\,\vdg}
\newcommand \vd {\,\,\vdg}
\newcommand \vii{\Land}
\newcommand \vuu{\Lor}
\newcommand \vdi[1] {\mathrel{\,\,\vdg_{#1}}}
\newcommand \Vii{\bigLand}
\newcommand \Vuu{\bigLor}

\setcounter{tocdepth}{4}



\theoremstyle{plain}
\newtheorem{theoreme}{Theorem}[section]
\newtheorem{theorem}{Theorem}[subsection]
\newtheorem{thdef}[theorem]{Theorem and definition}
\newtheorem{lemma}[theorem]{Lemma}
\newtheorem{corollary}[theorem]{Corollary}
\newtheorem{proposition}[theorem]{Proposition}
\newtheorem{propdef}[theorem]{Proposition and definition}
\newtheorem{plcc}[theorem]{Concrete local-global principle}
\newtheorem{fact}[theorem]{Fact}
\newtheorem{vastf}[theorem]{Formal \vst}
\newtheorem{vast}[theorem]{\vst}
\newtheorem{valsatz}[theorem]{\vst}
\newtheorem{precision}[theorem]{Precision}
\newtheorem{convention}[theorem]{Convention}

\newtheorem{theoremc}[theorem]{Theorem\etoz}
\newtheorem{lemmac}{Lemma\etoz}
\newtheorem{corollaryc}{Corollairy\etoz}
\newtheorem{propositionc}{Proposition\etoz}
\newtheorem{factc}{Fait\etoz}
\newtheorem{propdefc}[theorem]{Proposition and definition\etoz}

\newtheorem*{Principleofcoveringbyquotients}{Principle of covering by quotients}

\theoremstyle{definition}
\newtheorem{note}[theorem]{Note}
\newtheorem{context}[theorem]{Context}
\newtheorem{conjecture}[theorem]{Conjecture}
\newtheorem{definition}[theorem]{Definition}
\newtheorem{definitions}[theorem]{Definitions}
\newtheorem{descri}[theorem]{Description}
\newtheorem{notation}[theorem]{Notation}
\newtheorem{definota}[theorem]{Definition and notation} 
\newtheorem{problem}[theorem]{Problem}
\newtheorem{question}[theorem]{Question}
\newtheorem{ter}[theorem]{Terminology}

\theoremstyle{remark}
\newtheorem{notes}[theorem]{Notes}
\newtheorem{remark}[theorem]{Remark}
\newtheorem{remarks}[theorem]{Remarks}
\newtheorem{comment}[theorem]{Comment}
\newtheorem{comments}[theorem]{Comments}
\newtheorem{example}[theorem]{Example}
\newtheorem{examples}[theorem]{Examples}

\newcommand \Glio {\MA{\mathsf{Liog}}}

\newcommand{\vou}{\MA{\tsbf{ or }}}
\newcommand{\Vou}{\MA{\tsbf{OR}}}
\newcommand \EXists[1] {\tsbf{Introduce }{#1}\tsbf{ such that }\,}
\newcommand \vet {{\tsbf{,}}\,}
\newcommand \Atcl {\mathrm{Clat}}
\newcommand \Atclv {\mathrm{VClat}}
\newcommand \tcl {\mathrm{Clt}}
\newcommand \Tcl {\tcl}

\newcommand\comm{\rdb
\noi{\sl Comment. }}

\newcommand\COM[1]{\rdb
\noi{\sl Comment #1. }}

\newcommand\comms{\rdb
\noi{\sl Comments. }}

\newcommand\Pb{\rdb
\noi{\bf Problem. }}

\newcommand \rem{\rdb
\noi{\sl Remark. }}

\newcommand \REM[1]{\rdb
\noi{\sl Remark #1. }}

\newcommand \rems{\rdb
\noi{\sl Remarks. }}

\newcommand \exl{\rdb
\noi{\bf Example. }}

\newcommand \EXL[1]{\rdb
\noi{\bf Example #1. }}

\newcommand \exls{\rdb
\noi{\bf Examples. }}

\newcommand\gui[1]{``{#1}''}

\newcommand \thref[1] {Theorem~\ref{#1}}
\newcommand \paref[1] {page~\pageref{#1}}
\newcommand \pstfref[1] {Positivstellensatz formel~\ref{#1}}
\newcommand \pstref[1] {Positivstellensatz~\ref{#1}}

\newcommand \num {{n$^{\mathrm{ o}}$}}

\newcommand\subsubsec[1] {\subsubsection*{#1}}

\newcommand \recu {induction\xspace}
\newcommand \hdr {induction hypothesis\xspace}
\newcommand \ssi {if and only if\xspace}
\newcommand \cnes {necessary and sufficient condition\xspace}
\newcommand \spdg {without loss of generality\xspace}
\newcommand \Propeq {T.f.a.e.\xspace}
\newcommand \propeq {t.f.a.e.\xspace}
\newcommand \disept {17$^{th}$ Hilbert's problem\xspace}

\newcommand \cad {i.e.\ }
\newcommand \Cad {I.e.\ }


\newcommand \Amo {$\gA$-module\xspace}
\newcommand \Amos {$\gA$-modules\xspace}

\newcommand \Bmo {$\gB$-module\xspace}
\newcommand \Bmos {$\gB$-modules\xspace}

\newcommand \kmo {$\gk$-module\xspace}
\newcommand \kmos {$\gk$-modules\xspace}

\newcommand \Kmo {$\gK$-module\xspace}
\newcommand \Kmos {$\gK$-modules\xspace}

\newcommand \Lmo {$\gL$-module\xspace}
\newcommand \Lmos {$\gL$-modules\xspace}

\newcommand \Vmo {$\gV$-module\xspace}
\newcommand \Vmos {$\gV$-modules\xspace}

\newcommand \Zmo {$\gZ$-module\xspace}
\newcommand \Zmos {$\gZ$-modules\xspace}

\newcommand \ZZmo {$\ZZ$-module\xspace}
\newcommand \ZZmos {$\ZZ$-modules\xspace}

\newcommand \Ali {$\gA$-\ali}
\newcommand \Alis {$\gA$-\alis}

\newcommand \Bli {$\gB$-\ali}
\newcommand \Blis {$\gB$-\alis}

\newcommand \Cli {$\gC$-\ali}
\newcommand \Clis {$\gC$-\alis}

\newcommand \Alg {$\gA$-\alg}
\newcommand \Algs {$\gA$-\algs}

\newcommand \Blg {$\gB$-\alg}
\newcommand \Blgs {$\gB$-\algs}

\newcommand \Clg {$\gC$-\alg}
\newcommand \Clgs {$\gC$-\algs}

\newcommand \klg {$\gk$-\alg}
\newcommand \klgs {$\gk$-\algs}

\newcommand \Klg {$\gK$-\alg}
\newcommand \Klgs {$\gK$-\algs}

\newcommand \Llg {$\gL$-\alg}
\newcommand \Llgs {$\gL$-\algs}

\newcommand \Vlg {$\gV$-\alg}
\newcommand \Vlgs {$\gV$-\algs}

\newcommand \kev {$\gk$-vector space\xspace}
\newcommand \kevs {$\gk$-vector spaces\xspace}

\newcommand \Kev {$\gK$-vector space\xspace}
\newcommand \Kevs {$\gK$-vector spaces\xspace}

\newcommand \Kli {$\gK$-\ali}
\newcommand \Klis {$\gK$-\alis}

\newcommand \ac {algebraically closed\xspace}
\newcommand \alc {\agq closure\xspace}

\newcommand \adv {valuation domain\xspace}
\newcommand \advs {valuation domains\xspace}

\newcommand \advl {domain with divisors\xspace} 
\newcommand \advls {domains with divisors\xspace} 

\newcommand \advd {discrete valuation ring\xspace}
\newcommand \advds {discrete valuation rings\xspace}

\newcommand \agq {algebraic\xspace}

\newcommand \alg {algebra\xspace}
\newcommand \algs {algebras\xspace}

\newcommand \alge {étale \alg}
\newcommand \alges {étale \algs}

\newcommand \agB {Boolean \alg}

\newcommand \algo{algorithm\xspace}
\newcommand \algos{algorithms\xspace}

\newcommand \algq{algorithmic\xspace}

\newcommand \alo {local ring\xspace}
\newcommand \alos {local rings\xspace}

\newcommand \aloh {henselian \alo}
\newcommand \alohs {henselian \alos}

\newcommand \ali {\lin map\xspace}
\newcommand \alis {\lin maps\xspace}

\newcommand \alrd {\dcd \alo}
\newcommand \alrds {\dcd \alos}

\newcommand \anar {\ari \ri}
\newcommand \anars {\ari \ris}
\newcommand \Anars {\Ari \ris}

\newcommand \ari{arithmetic\xspace}

\newcommand \arv {valuation ring\xspace}
\newcommand \arvs {valuation rings\xspace}

\newcommand \auto {automorphism\xspace}
\newcommand \autos {automorphisms\xspace}


\newcommand \cac {algebraically closed field\xspace}
\newcommand \cacs {algebraically closed fields\xspace}

\newcommand \cara{characteristic\xspace}  
\newcommand \caras{characteristics\xspace}  

\newcommand \carn{characterisation\xspace}  
\newcommand \carns{characterisations\xspace}

\newcommand \cdf{field of fractions\xspace}

\newcommand \cdH{Hensel code\xspace}
\newcommand \cdHs{Hensel codes\xspace}

\newcommand \cdi{discrete field\xspace}  
\newcommand \cdis{discrete fields\xspace}  

\newcommand \cdr{field of roots\xspace}

\newcommand \cdv{variables change\xspace}  
\newcommand \cdvs{variables changes\xspace}  
\newcommand \cdvu{homogeneous \cdv}  
\newcommand \cdvus{homogeneous \cdvs}  

\newcommand \cli{integral closure\xspace}  
\newcommand \clis{integral closures\xspace}  

\newcommand \codi {discrete ordered field\xspace}
\newcommand \codis {discrete ordered fields\xspace}

\newcommand \coe {coefficient\xspace}
\newcommand \coes {coefficients\xspace}

\newcommand \cof {\cov}

\newcommand \coh {coherent\xspace}

\newcommand \coli {linear combination\xspace}
\newcommand \colis {linear combinations\xspace}

\newcommand \colo {local \cou}
\newcommand \colos {local \cous}

\newcommand \com {comaximal\xspace}

\newcommand \coo {coordinate\xspace}
\newcommand \coos {coordinates\xspace}

\newcommand \cop {complementary\xspace}

\newcommand \cosv {conservative\xspace}

\newcommand \cou {pair\xspace}
\newcommand \couh {henselian \cou}
\newcommand \cous {pairs\xspace}
\newcommand \couhs {henselian \cous}

\newcommand \crdl {residual field\xspace}

\newcommand \crl {corollary\xspace}
\newcommand \crls {corollaries\xspace}

\newcommand \cval{valued field\xspace}
\newcommand \cvals{valued fields\xspace}

\newcommand \cvar{\fva}
\newcommand \cvars{\fvas}

\newcommand \cvard{valuated discrete field\xspace}
\newcommand \cvards{valuated discrete fields\xspace}

\newcommand \cvd {valued discrete field\xspace}
\newcommand \cvds {valued discrete fields\xspace}

\newcommand \cvdh {henselian \cvd}
\newcommand \cvdhs {henselian \cvds}

\newcommand \cvdsc {separably closed valued discrete field\xspace}
\newcommand \cvdscs {separably closed valued discrete fields\xspace}

\newcommand \cvdac {algebraicalle closed valued discrete field\xspace}
\newcommand \cvdacs {algebraicalle closed valued discrete fields\xspace}

\newcommand \cvdu{\ultm \cvard}
\newcommand \cvdus{\ultm \cvards}

\newcommand \cvu{\ultm field\xspace}
\newcommand \cvus{\ultm fields\xspace}


\newcommand \dcd {residually discrete\xspace}

\newcommand \ddp {Pr\"ufer domain\xspace}
\newcommand \ddps {Pr\"ufer domains\xspace}

\newcommand \ddk {Krull dimension\xspace}

\newcommand \demo {proof\xspace}
\newcommand \dems {proofs\xspace}
\newcommand \demos {\dems}

\newcommand \deter {determinant\xspace}
\newcommand \deters {determinants\xspace}

\newcommand \dfn{definition\xspace}  
\newcommand \Dfn{Definition\xspace}  
\newcommand \Dfns{Definitions\xspace}  
\newcommand \dfns{definitions\xspace}  

\newcommand \dil{differential\xspace}  
\newcommand \diles{differentials\xspace}  

\newcommand \dij{disjunctive\xspace}

\newcommand \dimm {immediate description\xspace}
\newcommand \dimms {immediate descriptions\xspace}

\newcommand \discri{discriminant\xspace}
\newcommand \discris{discriminants\xspace}

\newcommand \dok {Dedekind domain\xspace}
\newcommand \doks {Dedekind domains\xspace}

\newcommand \dve {divisibility\xspace}

\newcommand \dvz {zerodivisor\xspace}
\newcommand \dvzs {zerodivisors\xspace}

\newcommand \eco{\com \elts}  

\newcommand \egmt{also\xspace} 

\newcommand \egt{equality\xspace} 
\newcommand \egts{equalities\xspace} 

\newcommand \elr{elementary\xspace}  

\newcommand \elt{element\xspace}  
\newcommand \elts{elements\xspace}  

\def \endo {endomorphism\xspace}
\def \endos {endomorphisms\xspace}

\newcommand \entrel {entailment relation\xspace}
\newcommand \entrels {entailment relations\xspace}

\newcommand \eqn  {equation\xspace}
\newcommand \eqns  {equations\xspace}

\newcommand \eqv  {equivalent\xspace}

\newcommand \eqvc  {equivalence\xspace}
\newcommand \eqvcs  {equivalences\xspace}

\newcommand \eseq{essentially equivalent\xspace} 
\newcommand \Eseq{Essentially equivalent\xspace} 

\newcommand \esid{essentially identical\xspace} 
\newcommand \Esid{Essentially identical\xspace} 

\newcommand \evc{vector space\xspace} 
\newcommand \evcs{vector spaces\xspace} 


\newcommand \fab {bounded \fcn}
\newcommand \fabs {bounded \fcns}

\newcommand \facile{\begin{proof}
Left to the reader.
\end{proof}}

\newcommand \fap {partial \fcn}
\newcommand \faps {partial \fcns}

\newcommand \fat {total \fcn}

\newcommand \fcn {factorisation\xspace}
\newcommand \fcns {factorisations\xspace}

\newcommand \fdi{strongly discrete\xspace} 

\newcommand \fit {faithfully\xspace}
\newcommand \fpt {\fit flat\xspace}

\newcommand \fmt {formellement\xspace}

\newcommand \flw {following\xspace}

\newcommand \fva {field with an absolute value\xspace}
\newcommand \fvas {field with an absolute value\xspace}


\newcommand\gmq{geometric\xspace}

\newcommand\gne{generalised\xspace}

\newcommand\gnl{general\xspace}

\newcommand\gnlt{generally\xspace}

\newcommand\gnn{generalisation\xspace}
\newcommand\gnns{generalisations\xspace}

\newcommand\gnq{generic\xspace}

\newcommand\gnr{generalise\xspace}  

\newcommand\grl{$\ell$-group\xspace}
\newcommand\grls{$\ell$-groups\xspace}

\newcommand \gtr{generator\xspace}  
\newcommand \gtrs{generators\xspace}  


\newcommand \homo {homomorphism\xspace}
\newcommand \homos {homomorphisms\xspace}

\newcommand \icl {integrally closed\xspace}

\newcommand \id {ideal\xspace}
\newcommand \ids {ideals\xspace}

\newcommand \idd {determinantal \id}
\newcommand \idds {determinantal \ids}

\newcommand \idema {maximal \id}
\newcommand \idemas {maximal \ids}

\newcommand \idep {prime \id}
\newcommand \ideps {prime \ids}

\newcommand \idemi {minimal prime\xspace}
\newcommand \idemis {minimal primes\xspace}

\newcommand \idf {Fitting \id}
\newcommand \idfs {Fitting \ids}

\newcommand \idm {idempotent\xspace}
\newcommand \idms {idempotents\xspace}

\newcommand \idp {principal \id}
\newcommand \idps {principal \ids}

\newcommand \idtr {indeterminate\xspace}
\newcommand \idtrs {indeterminates\xspace}

\newcommand \ifr {fractional \id}
\newcommand \ifrs {fractional \ids}

\newcommand \inteq {intuitively \eqv}

\newcommand \ird {irreducible\xspace}

\newcommand \itf {\tf \id}
\newcommand \itfs {\tf \ids}

\newcommand \iso {isomorphism\xspace}
\newcommand \isos {isomorphisms\xspace}

\newcommand \iv {invertible\xspace}

\newcommand \lec {reader\xspace}
\newcommand \llec {the reader\xspace}

\newcommand \lgb {local global\xspace}

\newcommand \LHe {Hensel Lemma\xspace}
\newcommand \lHe {Hensel lemma\xspace}
\newcommand \LHm {Multivariate \LHe}

\newcommand \lin {linear\xspace}

\newcommand \lint {linearly\xspace}

\newcommand \lon {localisation\xspace}
\newcommand \lons {localisations\xspace}

\newcommand \lop {\lot principal\xspace}

\newcommand \losd {\lsdz}
\newcommand \lsdz {\lot \sdz\xspace}

\def \lot {locally\xspace}

\newcommand \mlp {principal \lon matrix\xspace}
\newcommand \mlps {principal \lon matrices\xspace}

\newcommand \mnp {manipulation\xspace}
\newcommand \mnps {manipulations\xspace}
\newcommand \mnr {\elr \mnp}
\newcommand \mnrs {\elr \mnps}

\newcommand \mo {monoid\xspace}
\newcommand \mos {monoids\xspace}
\newcommand \moco {\com \mos}

\newcommand \mpf {\pf module\xspace}
\newcommand \mpfs {\pf modules\xspace}

\newcommand \mpn {\pn matrix\xspace}
\newcommand \mpns {\pn matrices\xspace}

\newcommand \mpr {\pro module\xspace}
\newcommand \mprs {\pro modules\xspace}

\newcommand \mprn {\prn matrix\xspace}
\newcommand \mprns {\prn matrices\xspace}

\newcommand \mptf {\ptf module\xspace}
\newcommand \mptfs {\ptf modules\xspace}

\newcommand \mrc {projective module of constant rank\xspace}
\newcommand \mrcs {projective modules of constant rank\xspace}


\newcommand \ncr{necessary\xspace}

\newcommand \ncrt{necessarily\xspace}

\newcommand \ndz {regular\xspace}

\newcommand \noe {noetherian\xspace}
\newcommand \noco {\noe\coh}

\newcommand \nst {Nullstellensatz\xspace}
\newcommand \nsts {Nullstellens\"atze\xspace}

\newcommand \odz {Zariski open set\xspace}

\newcommand \oqc {\qc open set\xspace}
\newcommand \oqcs {\qc open sets\xspace}

\newcommand \ort {orthogonal\xspace}


\newcommand \pa {saturated pair\xspace}
\newcommand \pas {saturated pairs\xspace}

\newcommand \pb{problem\xspace}  
\newcommand \pbs{problems\xspace}

\newcommand \peq {purely equational\xspace}

\newcommand \pf {finitely presented\xspace}

\newcommand \pgn {Newton polygon\xspace}
\newcommand \pgns {Newton polygons\xspace}

\newcommand \plg {\lgb principle\xspace}
\newcommand \plgs {\lgb principles\xspace}

\newcommand \pn {presentation\xspace}
\newcommand \pns {presentations\xspace}

\newcommand \pol {polynomial\xspace}
\newcommand \pols {polynomials\xspace}

\newcommand \polcar {characteristic \pol}
\newcommand \polcars {characteristic \pols}

\newcommand \polH {Hensel \pol}
\newcommand \polHs {Hensel \pols}

\newcommand \polmin {minimal \pol}
\newcommand \polmins {minimal \pols}

\newcommand \poll {polynomial\xspace}
\newcommand \Poll {Polynomial\xspace}

\newcommand \polu {monic \pol}
\newcommand \polus {monic \pols}

\newcommand \prc {rank constant \pro}

\newcommand \prmt {precisely\xspace}
\newcommand \Prmt {Precisely\xspace}

\newcommand \prn {projection\xspace}
\newcommand \prns {projections\xspace}

\newcommand \pro {projective\xspace}

\newcommand \proi {potential prime\xspace}
\newcommand \prois {potential primes\xspace}

\newcommand \proc {potential chain\xspace}
\newcommand \procs {potential chains\xspace}

\newcommand \proel {elementary \proc}
\newcommand \proels {elementary \procs}
\newcommand \proelo {\proel of length }
\newcommand \proelos {\proels of length }

\newcommand \prolo {\proc of length }
\newcommand \prolos {\procs of length }

\newcommand \prt {property\xspace}
\newcommand \prts {properties\xspace}

\newcommand \pst {Positivstellensatz\xspace}
\newcommand \psts {Positivstellens\"atze\xspace}

\newcommand \ptf {\tf \pro}


\newcommand \qc {quasi-compact\xspace}

\newcommand \qi {quasi integral\xspace}
\newcommand \qiri {pp-ring\xspace}
\newcommand \qiris {pp-rings\xspace}

\newcommand \ralg {Horn rule\xspace}
\newcommand \ralgs {Horn rules\xspace}

\newcommand \rcf {real closed field\xspace}
\newcommand \rcfs {real closed fields\xspace}

\newcommand \rdl {linear dependance relation\xspace}
\newcommand \rdls {linear dependance relations\xspace}

\newcommand \rdi {integral dependance relation\xspace}
\newcommand \rdis {integral dependance relations\xspace}

\newcommand \rdij {\dij rule\xspace}
\newcommand \rdijs {\dij rules\xspace}

\newcommand \rdv {valuative divisibility relation\xspace}
\newcommand \rdvs {valuative divisibility relations\xspace}

\newcommand \rdy {dynamical rule\xspace}
\newcommand \rdys {dynamical rules\xspace}

\newcommand \red {direct rule\xspace}
\newcommand \reds {direct rules\xspace}

\newcommand \rex {existential rule\xspace}
\newcommand \rexs {existential rules\xspace}

\newcommand \ri {ring\xspace}
\newcommand \ris {rings\xspace}


\newcommand \sad {dynamical algebraic structure\xspace}
\newcommand \sads {dynamical algebraic structures\xspace}
\newcommand \SAD {Dynamical algebraic structure\xspace}
\newcommand \SADs {Dynamical algebraic structures\xspace}

\newcommand \salg {algebraic structure\xspace}
\newcommand \salgs {algebraic structures\xspace}

\newcommand \sdz {without \dvz}

\newcommand \sfio {fundamental system of orthogonal idempotents\xspace}

\newcommand \sgr {\gtr set\xspace}
\newcommand \sgrs {\gtr sets\xspace}

\newcommand \sli {\lin \sys}
\newcommand \slis {\lin \syss}

\newcommand \spb {separable\xspace}
\newcommand \spl {separable\xspace}

\newcommand \sps {spectral space\xspace}
\newcommand \spss {spectral spaces\xspace}

\newcommand \ste {strictly étale\xspace}

\newcommand \stf {strictly finite\xspace}

\newcommand \sul {direct factor\xspace}

\newcommand \Syp {\Poll \sys}
\newcommand \sype {étale \syp}
\newcommand \sypes {étale \syps}
\newcommand \syp {\poll \sys}
\newcommand \syps {\poll \syss}

\newcommand \sys {system\xspace}
\newcommand \syss {systems\xspace}

\newcommand \sysN {Newton \syp}
\newcommand \SysN {Newton \syp}
\newcommand \sysNs {Newton \syps}

\newcommand \sysNe {étale Newton \sys}
\newcommand \sysNes {étale Newton \syss}

\newcommand \talg {Horn theory\xspace}
\newcommand \talgs {Horn theories\xspace}

\newcommand \tco {coherent theory\xspace}
\newcommand \tcos {coherent theories\xspace}

\newcommand \tdij {\dij theory\xspace}
\newcommand \tdijs {\dij theories\xspace}

\newcommand \Tdy {dynamical theory\xspace}
\newcommand \tdy {dynamical theory\xspace}
\newcommand \tdys {dynamical theories\xspace}

\newcommand \tel {regular theory\xspace}
\newcommand \tels {regular theories\xspace}

\newcommand \telri {cartesian theory\xspace}
\newcommand \telris {cartesian theories\xspace}

\newcommand \tf {finitely generated\xspace}

\newcommand \tfo {formal theory\xspace}
\newcommand \tfos {theory formelles\xspace}

\newcommand \tgm {\gmq theory\xspace}
\newcommand \tgms {\gmq theories\xspace}

\newcommand \Tho {Theorem\xspace}
\newcommand \Thos {Theorems\xspace}
\newcommand \tho {theorem\xspace}
\newcommand \thos {theorems\xspace}

\newcommand \tpe {purely equational theory\xspace}
\newcommand \tpes {purely equational theories\xspace}

\newcommand \trdi {distributive lattice\xspace}
\newcommand \trdis {distributive lattices\xspace}


\newcommand \ultm {ultrametric\xspace}
\newcommand \ultmf {ultrametric field\xspace}

\newcommand \unt {monic\xspace}

\newcommand \uvl {universal\xspace}

\newcommand \vala {absolute value\xspace}
\newcommand \valas {absolute values\xspace}

\newcommand \valn {valuation\xspace}
\newcommand \valns {valuations\xspace}

\newcommand \valu {\ultm \vala}
\newcommand \valus {\ultm \valas}

\newcommand \vfn {verification\xspace}
\newcommand \vfns {verifications\xspace}

\newcommand \vmd {vector \umd}
\newcommand \vmds {vectors \umds}

\newcommand \vst {Valuativstellensatz\xspace}
\newcommand \vsts {Valuativstellensätze\xspace}

\newcommand \vstf {formal \vst}
\newcommand \vstfs {formal \vsts}


\newcommand \zeH {Hensel zero\xspace}
\newcommand \zeHs {Hensel zeroes\xspace}

\newcommand \zed {zero-dimensional\xspace}
\newcommand \zedr {zero-dimensional reduced\xspace}

\newcommand \zmt {\tho de Zariski-Grothen\-dieck\xspace}


\newcommand \cov {constructive\xspace}

\newcommand \coma {\cov \maths}
\newcommand \clama {classical \maths}

\renewcommand \cot {constructively\xspace}

\newcommand \mathe {mathematical\xspace}
\newcommand \maths {mathematics\xspace}

\newcommand \matn {mathematician\xspace}

\newcommand \pte {excluded middle principle\xspace}

\newcommand \prco {\cov proof\xspace}
\newcommand \prcos {constructive proofs\xspace}

\newcommand \tcg {compactness theorem\xspace}
\newcommand \Tcgi {The \tcg implies the following result. }
%

\newcommand {\junk}[1]{}

\newcommand{\Cadre}[2]{%

\medskip%
\newskip\oldleftskip
\newskip\oldrightskip
\oldleftskip=\leftskip%
\oldrightskip=\rightskip%
\leftskip=-\tabcolsep%
\rightskip=-\tabcolsep%
\begin{center}\fbox{%
\begin{tabular}%
{p{#1\textwidth}}
\setlength{\parindent}{5mm}%
\vspace{-1.5mm}#2\vspace{1mm}%
\end{tabular}}\end{center}\par\medskip%
\leftskip=\oldleftskip%
\rightskip=\oldrightskip%
\setlength{\parindent}{6mm}}

\newcommand\boite[2]{\begin{minipage}[c]{#1cm}
     \centering {#2} \end{minipage}}
\newcommand\Boite[3]{\parbox[t][#1cm][c]{#2cm}{\boite{#2}{#3}}}

\newcommand{\Encadre}[1]{\Cadre{.8}{#1}}

\newcommand{\Cencadre}[1]{\Encadre{\vspace{-3mm}
\begin{center}
#1 \end{center}\vspace{-8mm}}}

\newcommand{\cen}{\centerline}
\newcommand \Grandcadre[1]{%
\begin{center}
\begin{tabular}{|c|}
\hline
~\\[-3mm]
#1\\[-3mm]
~\\
\hline
\end{tabular}
\end{center}}

\newcommand\dsp{\displaystyle}
\newcommand\ndsp{\textstyle}

\newcommand{\eop}{\hfill \mbox{$\Box$}}

\newcommand \noi {\noindent}
\renewcommand \ss {\smallskip}
\newcommand \sni {\ss\noi}
\newcommand \snii {\noi}
\newcommand \ms {\medskip}
\newcommand \mni {\ms\noi}
\newcommand \bs {\bigskip}
\newcommand \bni {\bs\noi}
\newcommand \hs {\qquad}
\newcommand \alb {\allowbreak}
\newcommand \ce {\centerline}


\renewcommand \le{\leqslant}
\renewcommand \leq{\leqslant}
\renewcommand \preceq{\preccurlyeq}
\renewcommand \ge{\geqslant}
\renewcommand \geq{\geqslant}
\renewcommand \succeq{\succcurlyeq}
\newcommand   \nneq {\mathrel{\#}}
\newcommand   \ineq {$_{\,\mathrel{\#}}$}

\newcommand\eti{^\times}
\newcommand \epr{^\perp}
\newcommand \etl{^*}
\newcommand \sta{^\star}
\newcommand \bu {{$\bullet$}}
\newcommand \bl {^\bullet}
\newcommand{\bul}{^{\bullet}}
\newcommand \eci {^\circ}
\newcommand \uci{\mathring}
\newcommand \ep[1]{^{(#1)}}
\newcommand \esh{^\sharp}
\newcommand \efl{^\flat}
\newcommand \eto{$^*$ }
\newcommand \etoz{$^*$}
\newcommand \ist{_\star}

\newcommand \Ast {\gA^{\!\star}}
\newcommand \Bst {\gB^{\star}}
\newcommand \iBA {_{\gB/\!\gA}}
\newcommand \iWV {_{\gW\!/\gV}}
\newcommand \Bo{\BB\mathrm{o}}
\newcommand \Ati {\gA^{\!\times}}
\newcommand \Bti {\gB^{\times}}
\newcommand \Vti {\gV^{\times}}
\newcommand \Atl {\gA^{\!*}}
\newcommand \Btl {\gB^{*}}
\newcommand \Ktl {\gK^{*}}
\newcommand \Vtl {\gV^{*}}
\newcommand{\KAt}{\gK\etl\!\sur{\Ati}}
\newcommand{\AAt}{\Atl\!\sur{\Ati}}

\newcommand \divi {\mid}

\newcommand\equidef{\buildrel{{\rm def}}\over{\;\Longleftrightarrow\;}}
\newcommand\eqdef{\buildrel{\rm def}\over {\;=\;}}
\newcommand\eqdefi{\buildrel{\rm def}\over {\;=\;}}


\newcommand \fraC[2] {{{#1}\over {#2}}}
\newcommand \formule[1]{{\left\{ {\arraycolsep2pt\begin{array}{lll} #1 \end{array}}\right.}}
\newcommand \formul[1]{{\left\{ {\arraycolsep2pt\begin{array}{rcl} #1 \end{array}}\right.}}
\newcommand \formu[1]{\arraycolsep2pt\begin{array}{rcl} #1 \end{array}}

\newcommand\mapright[1]{\smash{\mathop{\longrightarrow}\limits^{#1}}}
\newcommand\maprightto[1]{\smash{\mathop{\longmapsto}\limits^{#1}}}
\newcommand\mapdown[1]{\downarrow\rlap{$\vcenter{\hbox{$\scriptstyle
#1$}}$}}
\newcommand{\pref}[1]{\textup{\hbox{\normalfont(\ref{#1})}}}

\newcommand \abs[1] {\left|{#1}\right|}
\newcommand \abS[1] {\big|{#1}\big|}
\newcommand \aqo[2] {#1\sur{\gen{#2}}\!}
\newcommand \aQo[2] {#1/{\gEn{#2}}\!}
\newcommand \Aqo[2] {#1\sur{\big\langle{#2}\big\rangle}\!}
\newcommand \Al[1] {\Vi^{\!#1}}
\newcommand \ci[1] {{{#1}^\circ}}
\newcommand \crac[2] {\cro {\frac{#1}{#2}}}
\newcommand \cro[1] {\left[#1\right]}
\newcommand \eqdf[1] {\buildrel{#1}\over =}
\newcommand \equivdf[1] {\buildrel{#1}\over \longleftrightarrow}
\newcommand \frt[1] {\!\left|_{#1}\right.\!}
\newcommand \impdef[1] {\buildrel{#1}\over \Longrightarrow}
\newcommand \norme[1]{\left\lVert #1 \right\rVert}
\newcommand \Norme[1]{\big\lVert #1 \big\rVert}
\newcommand \tra[1] {{\,^{\rm t}\!#1}}
\newcommand \gen[1] {\left\langle{#1}\right\rangle}
\newcommand \gEn[1] {\langle{#1}\rangle}
\newcommand \geN[1] {\big\langle{#1}\big\rangle}
\newcommand \sing[1] {\left\{{#1}\right\}}
\newcommand \so[1] {\left\{\,{#1}\, \right\}}
\newcommand \soo[1] {\{\,{#1}\,\}}
\newcommand \sO[1]{\big\{{#1}\big\}}
\newcommand \sotq[2]{\so{#1\mathrel{;}#2}}
\newcommand \sootq[2]{\soo{#1\mathrel{;}#2}}
\newcommand \sotQ[2]{\sO{#1\mathrel{;}#2}}
\newcommand \sur[1] {\!\left/#1\right.}
\newcommand \und[1] {\underline{#1}}

\newcommand \Sqr {\mathrm{Sqr}}

\newcommand \idg[1] {|\,#1\,|}
\newcommand \idG[1] {\big|\,#1\,\big|}

\newcommand \norm[1] {\Vert\,#1\,\Vert}

\newcommand \dex[1] {[\,#1\,]}
\newcommand \deX[1] {\big[\,#1\,\big]}

\newcommand \lst[1] {[\,#1\,]}
\newcommand \lsT[1] {\big[\,#1\,\big]}

\newcommand{\mt}{\mapsto}

\newcommand{\llongrightarrow}{\relbar\joinrel\mkern-1mu\longrightarrow}
\newcommand{\lllongrightarrow}{\relbar\joinrel\mkern-1mu\llongrightarrow}
\newcommand{\llllongrightarrow}{\relbar\joinrel\mkern-1mu\lllongrightarrow}
\newcommand\simarrow{\vers{_\sim}}
\newcommand\isosim{\buildrel{_\sim}\over \longleftrightarrow }
\newcommand\vers[1]{\buildrel{#1}\over \longrightarrow }
\newcommand\vvers[1]{\buildrel{#1}\over \llongrightarrow }
\newcommand\vvvers[1]{\buildrel{#1}\over \lllongrightarrow }
\newcommand \lora {\longrightarrow}
\newcommand \llra {\llongrightarrow}
\newcommand \lllra {\lllongrightarrow}

\renewcommand \leq{\leqslant}
\renewcommand \geq{\geqslant}

\newcommand \som {\sum\nolimits}
\newcommand \Ex {{\exists}}
\newcommand \Tt {{\forall}}
\newcommand \te {\otimes}
\newcommand \vep{{\varepsilon}}


\newcommand\lra[1]{\langle{#1}\rangle}
\newcommand\lrb[1] {\llbracket #1 \rrbracket}
\newcommand\lrbd {\lrb{1..d}}
\newcommand\lrbn {\lrb{1..n}}
\newcommand\lrbl {\lrb{1..\ell}}
\newcommand\lrbm {\lrb{1..m}}
\newcommand\lrbk {\lrb{1..k}}
\newcommand\lrbp {\lrb{1..p}}
\newcommand\lrbq {\lrb{1..q}}
\newcommand\lrbr {\lrb{1..r}}
\newcommand\lrbs {\lrb{1..s}}


\newcommand \vda {\,\vdash\,}

\newcommand \vdw {\,\vdash_w\,}
\newcommand \dar[1] {\MA{\downarrow \!#1}}
\newcommand \uar[1] {\MA{\uparrow \!#1}}
\newcommand \clps[1] {{\downarrow #1\,\downarrow}}

\newcommand\tsbf[1]{\textsf{\textbf{\textup{#1}}}}
\newcommand\lab[1]{\item[\tsbf{#1}]}
\newcommand\Lab[1]{\rdb\item[\tsbf{#1}]\label{Ax#1}}
\newcommand\Tsbf[1]{\hyperref[Ax#1]{\tsbf{#1}}}

\newcommand\SA[1]{\rdb\sa{#1}\label{theorie#1}}
\newcommand\Sa[1]{\hyperref[theorie#1]{\sa{#1}}}
\newcommand\sa[1]{\hbox{\usefont{T1}{pzc}{m}{it}#1}\,}
\newcommand\sab[1]{\hbox{\usefont{T1}{pzc}{m}{it}#1}\,\,}
\newcommand\sA[1]{\hbox{\small\usefont{T1}{pzc}{m}{it}#1}\,}
\newcommand\sAb[1]{\hbox{\small\usefont{T1}{pzc}{m}{it}#1}\,\,}

\newcommand \snic[1] {\sni\centerline{$#1$}

\ss}

\newcommand \snac[1]{\sni
{\small\centering$#1$\par}

\ss}

\newcommand \eoe {\hbox{}\nobreak\hfill
\vrule width .5em height .5em depth 0mm \par \smallskip}

\newcommand \bal[1] {^\rK_{#1}}
\newcommand \ul[1] {_\rK^{#1}}

\newcommand \ov[1] {\overline{#1}}

\newcommand \wh[1] {\widehat{#1} }
\newcommand \wi[1] {\widetilde{#1} }

\newcommand\dessus[2]{{\textstyle {#1} \atop \textstyle {#2}}}

\newcommand\carray[2]{{\left[\begin{array}{#1} #2 \end{array}\right]}}
\newcommand\cmatrix[1]{\left[\matrix{#1}\right]}
\newcommand\clmatrix[1]{{\left[\begin{array}{lllllll} #1 \end{array}\right]}}
\newcommand\dmatrix[1]{\abs{\matrix{#1}}}
\newcommand\Cmatrix[2]{\setlength{\arraycolsep}{#1}\left[\matrix{#2}\right]}
\newcommand\Dmatrix[2]{\setlength{\arraycolsep}{#1}\left|\matrix{#2}\right|}

\newcommand \bloc[4]{\left[
\begin{array}{cc}
#1 & #2   \\
#3 & #4
\end{array}
\right]}

\newcommand \cm{em}

\newcommand{\blocs}[8]{%
{\setlength{\unitlength}{.0833\textwidth}
\tabcolsep0pt\renewcommand{\arraystretch}{0}%
\begin{tabular}{|c|c|}
\hline
\parbox[t][#3\cm][c]{#1\cm}{\begin{minipage}[c]{#1\cm}
\centering#5
\end{minipage}}&
\parbox[t][#3\cm][c]{#2\cm}{\begin{minipage}[c]{#2\cm}
\centering#6
\end{minipage}}\\
\hline
\parbox[t][#4\cm][c]{#1\cm}{\begin{minipage}[c]{#1\cm}
\centering#7
\end{minipage}}&
\parbox[t][#4\cm][c]{#2\cm}{\begin{minipage}[c]{#2\cm}
\centering#8
\end{minipage}}\\
\hline
\end{tabular}
}}

\newcommand \UneCol[1]{%
\sni\mbox{\hspace{.02\textwidth}%
\parbox[t]{.98\textwidth}{#1}%
}}

\newcommand \Unecol[1]{%
\sni\mbox{\hspace{.1\textwidth}%
\parbox[t]{.9\textwidth}{#1}%
}}

\newcommand \Deuxcol[4]{%
\sni\mbox{\parbox[t]{#1\textwidth}{#3}%
\hspace{.05\textwidth}%
\parbox[t]{#2\textwidth}{#4}}}

\newcommand \DeuxCol[2]{%
\Deuxcol{.475}{.475}{#1}{#2}}

\newcommand \DeuxCols[2]{%
\sni\mbox{\hspace{.02\textwidth}%
\parbox[t]{.475\textwidth}{#1}%
\hspace{.03\textwidth}%
\parbox[t]{.475\textwidth}{#2}}}

\newcommand \DeuxRegles[2]{%
\vspace{-1em}\DeuxCols
{\begin{enumerate}\itemsep=0pt  #1
\end{enumerate}
}
{\begin{enumerate}\itemsep=0pt  #2
\end{enumerate}
}
\vspace{-.3em}
}

\newcommand \UneRegle[2]{%
\vspace{-1em}\UneCol{
\begin{enumerate}
\lab{#1}{#2}
\end{enumerate}
}
\vspace{-.3em}
}

\newcommand \Regles[1]{%
\vspace{-1em}\UneCol{
\begin{enumerate}\itemsep=0pt
{#1}
\end{enumerate}
}
\vspace{-.3em}
}

\newcommand \regles[1]{%
\vspace{-1em}\Unecol{
\begin{enumerate}\itemsep=0pt
{#1}
\end{enumerate}
}
\vspace{-.3em}
}

\newcommand \itbu {\item[$\bullet$]}
\newcommand \labu {\lab{$\bullet$}}

\newcommand \Deuxbu[2]{%
\vspace{-.6em}\DeuxCols
{\begin{itemize}  #1
\end{itemize}}
{\begin{itemize}  #2
\end{itemize}}}

\makeatletter
\def\revddots{\mathinner{\mkern1mu\raise\p@
\vbox{\kern7\p@\hbox{.}}\mkern2mu
\raise4\p@\hbox{.}\mkern2mu\raise7\p@\hbox{.}\mkern1mu}}
\makeatother

\newcommand \BB{\mathbb {B}}
\newcommand \CC{\mathbb {C}}
\newcommand \FF{\mathbb {F}}
\newcommand \II{\mathbb {I}}
\newcommand \KK{\mathbb {K}}
\newcommand \MM{\mathbb {M}}
\newcommand \NN{\mathbb {N}}
\newcommand \ZZ{\mathbb {Z}}
\newcommand \OO{\mathbb {O}}
\newcommand \PP{\mathbb {P}}
\newcommand \QQ{\mathbb {Q}}
\newcommand \RR{\mathbb {R}}

\newcommand \gk {\mathbf{k}}
\newcommand \gkb {\ov\gk}
\newcommand \gl {\mathbf{l}}
\newcommand \gA {\mathbf{A}}
\newcommand \gB {\mathbf{B}}
\newcommand \gC {\mathbf{C}}
\newcommand \gD {\mathbf{D}}
\newcommand \gd {\mathbf{d}}
\newcommand \gDd {\mathbf{Dd}}
\newcommand \gE {\mathbf{E}}
\newcommand \gF {\mathbf{F}}
\newcommand \gG {\mathbf{G}}
\newcommand \gI {\mathbf{I}}
\newcommand \gIo {\mathbf{Io}}
\newcommand \gK {\mathbf{K}}
\newcommand \gKb {\ov\gK}
\newcommand \gAh {\widehat\gA}
\newcommand \gKh {\widehat\gK}
\newcommand \gKt {\widetilde\gK}
\newcommand \rhe {^{\mathrm{h}}}
\newcommand \rHe {^{\mathrm{H}}}
\newcommand \rNe {^{\mathrm{N}}}
\newcommand \rhs {^{\mathrm{hs}}}
\newcommand \Khe {\gK\rhe}
\newcommand \KHe {\gK\rHe}
\newcommand \KNe {\gK\rNe}
\newcommand \Kxi {\gK[\xi]}
\newcommand \gKp {{\gK'}}
\newcommand \gKw {\widetilde\gK}
\newcommand \gL {\mathbf{L}}
\newcommand \gLb {\ov\gL}
\newcommand \gLh {\widehat\gL}
\newcommand \gLw {\widetilde\gL}
\newcommand \Lhe {\gL\rhe}
\newcommand \gM {\mathbf{M}}
\newcommand \gP {\mathbf{P}}
\newcommand \gR {\mathbf{R}}
\newcommand \gRa {\mathbf{R}_\mathrm{a}}
\newcommand \gS {\mathbf{S}}
\newcommand \gT {\mathbf{T}}
\newcommand \gV {\mathbf{V}}
\newcommand \Ahe {\gA\!\rhe}
\newcommand \KV {(\gK,\gV)}
\newcommand \LW {(\gL,\gW)}

\newcommand \fmhe {\fm\rhe}
\newcommand \fmHe {\fm\rHe}
\newcommand \mhe {\fmhe}
\newcommand \mHe {\fmHe}
\newcommand \Ahs {\gA\rhs}
\newcommand \fmhs {\fm\rhs}
\newcommand \fmh {\widehat\fm}

\newcommand \fmA {{\fm_\gA}}
\newcommand \fmB {{\fm_\gB}}
\newcommand \fmC {{\fm_\gC}}
\newcommand \fmV {{\fm_\gV}}
\newcommand \fmVhe {\fmV\rhe}
\newcommand \fmVHe {\fm\gV\rHe}
\newcommand \fmVh {\widehat{\fmV}}
\newcommand \fmAh {\widehat{\fmA}}
\newcommand \fmti {\widetilde\fm}
\newcommand \fmVti {\widetilde\fmV}

\newcommand \Vhe {\gV\rhe}
\newcommand \VHe {\gV\rHe}
\newcommand \VNe {\gV\rNe}
\newcommand \gVh {\widehat\gV}
\newcommand \gAt {\widetilde\gA}
\newcommand \gVt {\widetilde\gV}
\newcommand \gVp {{\gV'}}
\newcommand \gW {\mathbf{W}}
\newcommand \gX {\mathbf{X}}
\newcommand \gZ {\mathbf{Z}}

\newcommand \Gat {\wh{\Gamma}}
\newcommand \Ksep {{\gK^\mathrm{sep}}}
\newcommand \ksep {{\gk^\mathrm{sep}}}
\newcommand \Kac {{\gK^\mathrm{ac}}}
\newcommand \Vsep {{\gV^\mathrm{sep}}}
\newcommand \fmsep {{\fm^\mathrm{sep}}}
\newcommand \vsep {{v_\mathrm{sep}}}
\newcommand \Vac {{\gV^\mathrm{ac}}}
\newcommand \fmac {{\fm^\mathrm{ac}}}

\newdimen\xyrowsp
\xyrowsp=3pt
\newcommand{\SCO}[6]{
\xymatrix @R = \xyrowsp {
                                  &1 \ar@{-}[dl] \ar@{-}[dr] \\
#3 \ar@{-}[ddr]                   &   & #6 \ar@{-}[ddl] \\
                                  &\bullet\ar@{-}[d] \\
                                  &\bullet   \\
#2 \ar@{-}[ddr] \ar@{-}[uur]      &   & #5 \ar@{-}[ddl] \ar@{-}[uul] \\
                                  &\bullet \ar@{-}[d] \\
                                  &\bullet  \\
#1 \ar@{-}[uur]                   &   & #4 \ar@{-}[uul] \\
                                  & 0 \ar@{-}[ul] \ar@{-}[ur] \\
}
}

\newcommand \Adj {\MA{\mathrm{Adj}}}
\newcommand \adj {\MA{\mathrm{adj}}}
\newcommand \Adu {\MA{\mathrm{Adu}}}
\newcommand \Ann {\mathrm{Ann}}
\newcommand \Atom {\mathrm{Atom}}
\newcommand \Aut {\MA{\mathrm{Aut}}}
\newcommand \BZ {\MA{\mathrm{BZ}}}
\newcommand \car {\MA{\mathrm{car}}}
\newcommand \Cl {\MA{\mathrm{Cl}}}
\newcommand \ClW {\MA{\mathrm{Cl}_{\mathrm{W}}}}
\newcommand \Coker {\MA{\mathrm{Coker}}}
\newcommand \Cont{\mathrm{Co}}
\newcommand \Dc {\MA{\mathrm{Dc}}}
\newcommand \DDiv {\MA{\mathrm{Dv}}}
\renewcommand \det {\MA{\mathrm{det}}}
\renewcommand \deg {\MA{\mathrm{deg}}}
\newcommand \Diag {\MA{\mathrm{Diag}}}
\newcommand \disc {\MA{\mathrm{disc}}}
\newcommand \Disc {\MA{\mathrm{Disc}}}
\newcommand \Div {\MA{\mathrm{Div}}}
\newcommand \DivA {\Div\gA }
\newcommand \DivAp {(\Div\gA)^{+} }
\newcommand \DivB {\Div\gB }
\newcommand \DivBp {(\Div\gB)^{+} }
\newcommand \DkM {\MA{\mathrm{DkM}}}
\newcommand \dv {\MA{\mathrm{div}} }
\newcommand \dvA {\dv_\gA }
\newcommand \dvB {\dv_\gB }
\newcommand \ev {{\mathrm{ev}}}
\newcommand \End {\MA{\mathrm{End}}}
\newcommand \fsac {\MA{\mathrm{fsa}}}
\newcommand \Fix {\MA{\mathrm{Fix}}}
\newcommand \Frac {\MA{\mathrm{Frac}}}
\newcommand \Gal {\MA{\mathrm{Gal}}}
\newcommand \Gfr {\MA{\mathrm{Gfr}}}
\newcommand \Gr {\MA{\mathrm{Gr}}}
\newcommand \gr {\MA{\mathrm{gr}}}
\newcommand \Gram {\MA{\mathrm{Gram}}}
\newcommand \gram {\MA{\mathrm{gram}}}
\newcommand \Grl {\MA{\mathrm{Grl}}}
\newcommand \hauteur {\mathrm{hauteur}}
\newcommand \Hom {\MA{\mathrm{Hom}}}
\newcommand \Id {\MA{\mathrm{Id}}}
\newcommand \Iv {\MA{\mathrm{Iv}}}
\newcommand \I {\mathrm{I}}
\newcommand \Idif {\MA{\mathrm{Idif}}}
\newcommand \Idv {\MA{\mathrm{Idv}}}
\newcommand \Ifr {\MA{\mathrm{Ifr}}}
\newcommand \Icl {\MA{\mathrm{Icl}}}
\renewcommand \Im {\MA{\mathrm{Im}}}
\newcommand \Inf {\MA{\mathrm{Inf}}}
\newcommand \Itf {\MA{\mathrm{Itf}}}
\newcommand \Ker {\MA{\mathrm{Ker}}}
\newcommand \Lst {\MA{\mathrm{Lst}}}
\newcommand \LIN {\mathrm{Lin}}
\newcommand \Lsf {\MA{\mathrm{Lsf}}}
\newcommand \Mat {\MA{\mathrm{Mat}}}
\newcommand \Mip {\mathrm{Min}}
\newcommand \md {\mathrm{md}}
\newcommand \Mgcd {\MA{\mathrm{Mgcd}}}
\renewcommand \mod {\;\mathrm{mod}\;}
\newcommand \Mor {\MA{\mathrm{Mor}}}
\newcommand \NDc {\MA{\mathrm{NDc}}}
\newcommand \poids {\mathrm{poids}}
\newcommand \poles {\hbox {\rm p\^oles}}
\newcommand \pgcd {\MA{\mathrm{pgcd}}}
\newcommand \ppcm {\MA{\mathrm{ppcm}}}
\newcommand \Rad {\MA{\mathrm{Rad}}}
\newcommand \Reg {\MA{\mathrm{Reg}}}
\newcommand \rg{\MA{\mathrm{rg}}}
\newcommand \rgst {\mathrm{rgst}}
\newcommand \Res {\mathrm{Res}}
\newcommand \Rs {\MA{\mathrm{Rs}}}
\newcommand \rPr{\MA{\mathrm{Pr}}}
\newcommand \Rv {\mathrm{Rv}}
\newcommand \Sli {\MA{\mathrm{Sli}}}
\newcommand \Som {\MA{\mathrm{Som}}}
\newcommand \Sup {\MA{\mathrm{Sup}}}
\newcommand \Sace {\MA{\mathrm{Sace}}}
\newcommand \Smtf {\MA{\mathrm{Smtf}}}
\newcommand \Stp {\MA{\mathrm{Stp}}}
\newcommand \St {\mathrm{St}}
\newcommand \Tri {\MA{\mathrm{Tri}}}
\newcommand \Tor {\MA{\mathrm{Tor}}}
\newcommand \tr {\MA{\mathrm{tr}}}
\newcommand \Tr {\MA{\mathrm{Tr}}}
\newcommand \Tsc {\MA{\mathrm{Tsch}}}
\newcommand \Um {\MA{\mathrm{Um}}}
\newcommand \val {\MA{\mathrm{val}}}

\newcommand \SIPD {\MA{\mathrm{SIPD}}}
\newcommand \ARC {\MA{\mathrm{ARC}}}
\newcommand \AFR {\MA{\mathrm{AFR}}}
\newcommand \AFRRV {\MA{\mathrm{AFRRV}}}
\newcommand \AFRNZ {\MA{\mathrm{AFRNZ}}}
\newcommand \PPM {\MA{\mathrm{PPM}}}
\newcommand \PB {\MA{\mathrm{PB}}}

\newcommand \Suslin{{\rm Suslin}}
\newcommand{\DBxk}{{\Der \gk\gB\xi}}%
\newcommand{\DAxk}{{\Der \gk\gA\xi}}%
\newcommand{\DkXxk}{{\Der \gk\kuX\xi}}%
\newcommand \DAbul {\rD_{\!\Abul}}

\newcommand \sfP {\mathsf{P}}

\newcommand\MA[1]{\mathop{#1}\nolimits}

\newcommand \Cdim {\MA{\mathsf{Cdim}}}
\newcommand \Divdim {\MA{\mathsf{Divdim}}}
\newcommand \Glo {\MA{\mathsf{Glo}}}
\newcommand \GK {\MA{\mathsf{GK}}}
\newcommand \GKO {\MA{\mathsf{GK}_0}}
\newcommand \HO {\MA{\mathsf{H}_0}}
\newcommand \HOp {\MA{\mathsf{H}_0^+}}
\newcommand \Hdim {\MA{\mathsf{Hdim}}}
\newcommand \HeA {{\Heit\gA}}
\newcommand \Heit {\MA{\mathsf{Heit}}}
\newcommand \Hspec {\MA{\mathsf{Hspec}}}
\newcommand \Jdim {\MA{\mathsf{Jdim}}}
\newcommand \jdim {\MA{\mathsf{jdim}}}
\newcommand \Jspec {\MA{\mathsf{Jspec}}}
\newcommand \jspec {\MA{\mathsf{jspec}}}
\newcommand \KO {\MA{\mathsf{K}_0}}
\newcommand \KOp {\MA{\mathsf{K}_0^+}}
\newcommand \KTO {\wi{\mathsf{K}}_0}
\newcommand \Kdim {\MA{\mathsf{Kdim}}}
\newcommand \Lin {\mathsf{L}}
\newcommand \Max {\MA{\mathsf{Max}}}
\newcommand \Min {\MA{\mathsf{Min}}}
\newcommand \OQC {\MA{\mathsf{Oqc}}}
\newcommand \Pic {\MA{\mathsf{Pic}}}
\newcommand \Reel {\MA{\mathsf{Reel}}}
\newcommand \Spec {\MA{\mathsf{Spec}}}
\newcommand \Speclin {\MA{\mathsf{Speclin}}}
\newcommand \Spv {\MA{\mathsf{Spv}}}
\newcommand \Spev {\MA{\mathsf{Spev}}}
\newcommand \Sper {\MA{\mathsf{Sper}}}
\newcommand \Val {\MA{\mathsf{Val}}}
\newcommand \Valp {\MA{\mathsf{Val'}}}
\newcommand \SpecA {\Spec\gA}
\newcommand \SpevA {\Spev\gA}
\newcommand \SpvA {\Spv\gA}
\newcommand \SperA {\Sper\gA}
\newcommand \SpecT {\Spec\gT}
\newcommand \Zar {\MA{\mathsf{Zar}}}
\newcommand \ZF {\MA{\mathsf{ZF}}}
\newcommand \ValA {{\Val\gA}}
\newcommand \ValpA {{\Valp\gA}}
\newcommand \ZarA {{\Zar\gA}}

\newcommand \cA {{\cal A}}
\newcommand \cB {{\cal B}}
\newcommand \cC {{\cal C}}
\newcommand \cD {{\cal D}}
\newcommand \cI {{\cal I}}
\newcommand \cJ {{\cal J}}
\newcommand \cF {{\cal F}}
\newcommand \cH {{\cal H}}
\newcommand \cK {{\cal K}}
\newcommand \cL {{\cal L}}
\newcommand \cM {{\cal M}}
\newcommand \cN {{\cal N}}
\newcommand \cP {{\cal P}}
\newcommand \cQ {{\cal Q}}
\newcommand \cR {{\cal R}}
\newcommand \cS {{\cal S}}
\newcommand \cT {{\cal T}}
\newcommand \cV {{\cal V}}

\newcommand \ccd{\mathcal{CD}}
\newcommand \cco{\mathcal{CO}}

\newcommand \Cin{C^{\infty}}

\newcommand \SK {\cS^\rK}
\newcommand \IK {\cI^\rK}
\newcommand \JK {\cJ^\rK}
\newcommand \IH {\cI^\rH}
\newcommand \JH {\cJ^\rH}

\newcommand \JAC {J}
\newcommand \Jac {\mathrm{Jac}}
\newcommand \rc {\mathrm{c}}
\newcommand \Df {\MA{\mathrm{Df}}}
\newcommand \Dfr {\MA{\mathrm{Df}}^\mathrm{R}}
\newcommand \Dfmc {\MA{\mathrm{Dfmc}}}
\newcommand \rd {\mathrm{d}}
\newcommand \rv {\mathrm{v}}
\newcommand \rI {\mathrm{I}}
\newcommand \Ic {\mathrm{Ic}}
\newcommand \rC {\mathrm{C}}
\newcommand \rD {\mathrm{D}}
\newcommand \rF {\mathrm{F}}
\newcommand \rG {\mathrm{G}}
\newcommand \rH {\mathrm{H}}
\newcommand \rJ {\mathrm{J}}
\newcommand \Li {\MA{\mathrm{Li}}}
\newcommand \rK {\mathrm{K}}
\newcommand \rL {\mathrm{L}}
\newcommand \Mc {\mathrm{Mc}}
\newcommand \rN {\mathrm{N}}
\newcommand \rP {\mathrm{P}}
\newcommand \rR {\mathrm{R}}
\newcommand \rmSa {\MA{\mathrm{Sa}}}
\newcommand \rmSamc {\MA{\mathrm{Samc}}}
\newcommand \rS {\mathrm{S}}
\newcommand \rU {\mathrm{U}}
\newcommand \rV {\mathrm{V}}
\newcommand \DA {\rD_{\!\gA}}
\newcommand \JA {\rJ_\gA}
\newcommand \JT {\rJ_\gT}

\newcommand\fa{\mathfrak{a}}
\newcommand\fb{\mathfrak{b}}
\newcommand\fc{\mathfrak{c}}
\newcommand\fA{\mathfrak{A}}
\newcommand\fB{\mathfrak{B}}
\newcommand\fD{\mathfrak{D}}
\newcommand\fI{\mathfrak{i}}
\newcommand\fII{\mathfrak{I}}
\newcommand\fj{\mathfrak{j}}
\newcommand\fJ{\mathfrak{J}}
\newcommand\fF{\mathfrak{F}}
\newcommand\ff{\mathfrak{f}}
\newcommand\ffg{\mathfrak{g}}
\newcommand\fG{\mathfrak{G}}
\newcommand\fh{\mathfrak{h}}
\newcommand\fl{\mathfrak{l}}
\newcommand\fm{\mathfrak{m}}
\newcommand\mV{{\fm_\gV}}
\newcommand\mW{{\fm_\gW}}
\newcommand\fM{\mathfrak{M}}
\newcommand\fp{\mathfrak{p}}
\newcommand\fP{\mathfrak{P}}
\newcommand\fq{\mathfrak{q}}
\newcommand\fU{\mathfrak{U}}
\newcommand\fV{\mathfrak{V}}
\newcommand\fx{\mathfrak{x}}
\newcommand\fy{\mathfrak{y}}

\newcommand \scC{\mathscr{C}}
\newcommand \scR{\mathscr{R}}

\newcommand{\bma}{\bm{a}}
\newcommand{\bmb}{\bm{b}}
\newcommand{\bmc}{\bm{c}}
\newcommand{\bmd}{\bm{d}}
\newcommand{\bme}{\bm{e}}
\newcommand{\bmf}{\bm{f}}
\newcommand{\bmu}{\bm{u}}
\newcommand{\bmv}{\bm{v}}
\newcommand{\bmw}{\bm{w}}
\newcommand{\bmy}{\bm{y}}
\newcommand{\bmx}{\bm{x}}
\newcommand{\bmz}{\bm{z}}

\newcommand \LLPO {\tsbf{LLPO}}
\newcommand \LPO  {\tsbf{LPO}}

\newcommand \Zg {{\Z[G]}}

\newcommand \vu {\vee} 
\newcommand \vi {\wedge} 
\newcommand \Vu {\bigvee}
\newcommand \Vi {\bigwedge}
\newcommand \im {\rightarrow} 
\newcommand \da {\,\downarrow\!}

\newcommand \vdu[1] {\vdash^{#1}}
\newcommand \vdb[1] {\vdash_{#1}}

\newcommand \Vrai {\mathsf{True}}
\newcommand \Faux {\mathsf{False}}
\newcommand \Un {\mathbf{1}}
\newcommand \Deux {\mathbf{2}}
\newcommand \Trois {\mathbf{3}}
\newcommand \Quatre {\mathbf{4}}
\newcommand \Cinq {\mathbf{5}}

\newcommand \una {{\underline{a}}}
\newcommand \ual {{\underline{\alpha}}}
\newcommand \ua  {{\underline{a}}}
\newcommand \ub  {{\underline{b}}}
\newcommand \ube {{\underline{\beta}}}
\newcommand \uc{{\underline{c}}}
\newcommand \ud  {{\underline{d}}}
\newcommand \udel{{\underline{\delta}}}
\newcommand \ue  {{\underline{e}}}
\newcommand \uf  {{\underline{f}}}
\newcommand \ug  {{\underline{g}}}
\newcommand \uh  {{\underline{h}}}
\newcommand \uga {{\underline{\gamma}}}
\newcommand \uP  {{\underline{P}}}
\newcommand \ur{{\underline{r}}}
\newcommand \ut{{\underline{t}}}
\newcommand \uu{{\underline{u}}}
\newcommand \ux {{\underline{x}}}
\newcommand \uxi {{\underline{\xi}}}
\newcommand \uX {\underline{X}}
\newcommand \uy{{\underline{y}}}
\newcommand \uY  {{\underline{Y}}}
\newcommand \uz{{\underline{z}}}
\newcommand \uze {{\underline{0}}}

\newcommand \ak {a_1,\ldots,a_k}
\newcommand \am {a_1,\ldots,a_m}
\newcommand \an {a_1,\ldots,a_n}
\newcommand \aq {a_1,\ldots,a_q}
\newcommand \aln {\alpha_1,\ldots,\alpha_n}
\newcommand \bn {b_1,\ldots,b_n}
\newcommand \bzn {b_0,\ldots,b_n}
\newcommand \bbm {b_1,\ldots,b_m}
\newcommand \ck {c_1,\ldots,c_k}
\newcommand \cq {c_1,\ldots,c_q}
\newcommand \rcr {c_1,\ldots,c_r}
\newcommand \gan {\gamma_1,\ldots,\gamma_n}
\newcommand \un {u_1,\ldots,u_n}
\newcommand \xk {x_1,\ldots,x_k}
\newcommand \Xk {X_1,\ldots,X_k}
\newcommand \xm {x_1,\ldots,x_m}
\newcommand \Xm {X_1,\ldots,X_m}
\newcommand \fn {f_1,\ldots,f_n}
\newcommand \lfm {f_1,\ldots,f_m}
\newcommand \Fn {F_1,\ldots,F_n}
\newcommand \lFm {F_1,\ldots,F_m}
\newcommand \Fp {\FF_p}
\newcommand \Zp {\ZZ_p}
\newcommand \Qp {\QQ_p}
\newcommand \Cp {\CC_p}
\newcommand \sn {s_1,\ldots,s_n}
\newcommand \xn {x_1,\ldots,x_n}
\newcommand \xzn {x_0,\ldots,x_n}
\newcommand \xhn {x_0:\ldots:x_n}
\newcommand \Xn {X_1,\ldots,X_n}
\newcommand \xr {x_1,\ldots,x_r}
\newcommand \Xr {X_1,\ldots,X_r}
\newcommand \xin {\xi_1,\ldots,\xi_n}
\newcommand \xizn {\xi_0,\ldots,\xi_n}
\newcommand \xihn {\xi_0:\ldots:\xi_n}
\newcommand \ym {y_1,\ldots,y_m}
\newcommand \yr {y_1,\ldots,y_r}
\newcommand \Yr {Y_1,\ldots,Y_r}
\newcommand \Yn {Y_1,\ldots,Y_n}
\newcommand \Ym {Y_1,\ldots,Y_m}
\newcommand \yn {y_1,\ldots,y_n}

\newcommand \AT {\gA[T]}
\newcommand \AX {\gA[X]}
\newcommand \Ax {\gA[x]}
\newcommand \Ared {\gA\red}
\newcommand \AuX {\gA[\uX]}
\newcommand \Aux {\gA[\ux]}
\newcommand \ArX {\gA\lra X}
\newcommand \Axn {\gA[\xn]}
\newcommand \AXn {\gA[\Xn]}

\newcommand \AXm {\gA[\Xm]}
\newcommand \KKXm {{\KK[\Xm]}}
\newcommand \KKuX {{\KK[\uX]}}

\newcommand \AY {\gA[Y]}
\newcommand \Ayn {\gA[\yn]}

\newcommand \BuX {\gB[\uX]}
\newcommand \BuY {\gB[\uY]}
\newcommand \BX {{\gB[X]}}
\newcommand \BT {{\gB[T]}}
\newcommand \BY {{\gB[Y]}}
\newcommand \Bxn {\gB[\xn]}
\newcommand \BXn {{\gB[\Xn]}}
\newcommand \BYm {\gB[\Ym]}

\newcommand \CT {{\gC[T]}}
\newcommand \CX {{\gC[X]}}
\newcommand \CXn {{\gC[\Xn]}}

\newcommand \kX {{\gk[X]}}
\newcommand \KX {\gK[X]}
\newcommand \Kx {\gK[x]}
\newcommand \VX {\gV[X]}
\newcommand \Vx {\gV[x]}
\newcommand \KT {\gK[T]}
\newcommand \KuX {\gK[\uX]}
\newcommand \kuX {\gk[\uX]}
\newcommand \VuX {\gV[\uX]}
\newcommand \Vuxi {\gV[\uxi]}
\newcommand \Vux {\gV[\ux]}
\newcommand \Kux {\gK[\ux]}
\newcommand \kux {\gk[\ux]}
\newcommand \KXk {\gK[\Xk]}
\newcommand \KXm {\gK[\Xm]}
\newcommand \KXn {\gK[\Xn]}
\newcommand \kxm {\gk[\xm]}
\newcommand \kxn {\gk[\xn]}
\newcommand \Ky {\gK[y]}
\newcommand \KY {\gK[Y]}
\newcommand \Kz {\gK[z]}
\newcommand \KZ {\gK[Z]}
\newcommand \KYn {\gK[\Yn]}
\newcommand \KYm {\gK[\Ym]}
\newcommand \kXr {\gk[\Xr]}
\newcommand \KXr {\gK[\Xr]}
\newcommand \Kxr {\gK[\xr]}

\newcommand \Vxn {\gV[\xn]}
\newcommand \VXn {\gV[\Xn]}

\newcommand \KuY {\gK[\uY]}
\newcommand \Kuy {\gK[\uy]}
\newcommand \Kyn {\gK[\yn]}
\newcommand \Kyr {\gK[\yr]}
\newcommand \kYr {\gk[\Yr]}
\newcommand \KYr {\gK[\Yr]}
\newcommand \Kxn {\gK[\xn]}

\newcommand \LuX {\gL[\uX]}
\newcommand \lXn {\gl[\Xn]}
\newcommand \lxn {\gl[\xn]}
\newcommand \LXn {\gL[\Xn]}
\newcommand \lXr {\gl[\Xr]}
\newcommand \LXr {\gL[\Xr]}
\newcommand \lYr {\gl[\Yr]}
\newcommand \LYr {\gL[\Yr]}

\newcommand \QQXn {\QQ[\Xn]}

\newcommand \Rx {\gR[x]}
\newcommand \Rux {\gR[\ux]}
\newcommand \RuX {\gR[\uX]}
\newcommand \RXk {{\gR[\Xk]}}
\newcommand \RXm {{\gR[\Xm]}}
\newcommand \Rxm {{\gR[\xm]}}
\newcommand \RXn {{\gR[\Xn]}}
\newcommand \Rxn {{\gR[\xn]}}
\newcommand \RXzn {{\gR[\Xzn]}}
\newcommand \Rxzn {{\gR[\xzn]}}
\newcommand \RXr {\gR[\Xr]}

\newcommand \Ruy {\gR[\uy]}
\newcommand \Ryn {{\gR[\yn]}}
\newcommand \RRX {\RR[X]}
\newcommand \RRXn {\RR[\Xn]}
\newcommand \RRxn {\RR[\xn]}
\newcommand \RYr {\gR[\Yr]}
\newcommand \RRuX {\RR[\uX]}
\newcommand \RRux {\RR[\ux]}

\newcommand \CCX {\CC[X]}

\newcommand \ZZXn {\ZZ[\Xn]}
\newcommand \ZG {\ZZ[G]}

\newcommand \RRXk {{\RR[\Xk]}}
\newcommand \RRxk {{\RR[\xk]}}
\newcommand \RRXm {{\RR[\Xm]}}
\newcommand \RRxm {{\RR[\xm]}}

\newcommand \lfs {f_1,\ldots,f_s}
\newcommand \lfn {f_1,\ldots,f_n}

\newcommand \Gn  {\gG_n}
\newcommand \Gnk {\gG^n_{k}}
\newcommand \Gnr {\gG^n_{r}}
\newcommand \cGn {\cG_n}
\newcommand \cGnk{\cG_{n,k}}
\newcommand \GGn {\GG^n}
\newcommand \GGnk{\GGn_{k}} 
\newcommand \GGnr{\GGn_{r}}
\newcommand \GA  {\mathbb{GA}}
\newcommand \GAn {\GA^n}  
\newcommand \GAq {\GA^q}
\newcommand \GAnk{\GAn_{k}}
\newcommand \GAnr{\GAn_{r}}
\newcommand \GL {\mathbb{GL}}
\newcommand \GLn {{\GL_n}}
\newcommand \SL {\mathbb{SL}}
\newcommand \SLn {{\SL_n}}
\newcommand \EE {\mathbb{E}}
\newcommand \En {\EE_n}
\newcommand \Pn {\PP^n}
\newcommand \An {\AA^n}
\newcommand \Sl {\mathbf{SL}}
\newcommand \Sln {{\Sl_n}}

\newcommand \Mm {\MM_{m}}
\newcommand \Mn {\MM_{n}}
\newcommand \Mk {\MM_{k}}
\newcommand \Mq {\MM_{q}}
\newcommand \Mr {\MM_{r}}
\newcommand \MMn {\MM_{n}}

\newcommand \Pf {{{\cal P}_{\mathrm{f}}}}
\newcommand \Pfe {{\rm P}_{{\rm fe}}}

\newcommand\hsz{\\ }
\newcommand\hsu{\\ \hspace*{4mm}}
\newcommand\hsd{\\ \hspace*{8mm}}
\newcommand\hst{\\ \hspace*{1,2cm}}
\newcommand\hsq{\\ \hspace*{1,6cm}}
\newcommand\hsc{\\ \hspace*{2cm}}
\newcommand\hsix{\\ \hspace*{2,4cm}}
\newcommand\hsept{\\ \hspace*{2,8cm}}

\title{Multivariate Hensel Lemma for ultrametric fields}

\author{M.-E. Alonso\thanks{Universidad Complutense, Madrid, Espa\~na. {\tt mariemi@mat.ucm.es}}
 \and  Henri Lombardi
\thanks{Université de Franche-Comté, Laboratoire de mathématiques de Besançon, UMR CNRS 6623, 16 route de Gray, 25000
Besançon, France. {\tt henri.lombardi@univ-fcomte.fr}}
\and
Stefan Neuwirth
\thanks{Université de Franche-Comté, Laboratoire de mathématiques de Besançon, UMR CNRS 6623, 16 route de Gray, 25000
Besançon, France. {\tt stefan.neuwirth@univ-fcomte.fr}}}

\date{\today}
\maketitle

\vspace{-1em}

\begin{abstract} 
The Multivariate Hensel Lemma for local rings is usually proved as a consequence of the Grothendieck version of Zariski's Main Theorem. This version deals with a more general situation that is a priori much more difficult.
In this paper, we give a direct  proof of the Multivariate Hensel Lemma for ultrametric fields, in the framework of constructive mathematics and without using~ZMT\@. In the framework of classical mathematics, our result entails the Lemma for rank-one valued fields. 
\end{abstract}

\rdb
\label{beginenglish}

\medskip \noindent {\bf Keywords:} Multivariate Hensel Lemma,
Valued field,
Ultrametric field,
Henselisation of a local ring,
Henselisation of a valued field,
Constructive mathematics. 

\smallskip  \noindent {\bf MSC2020:} 13B40, 13J15, 03F65.

\setcounter{tocdepth}{4}
\markboth{Contents}{Contents}

\small
\printcontents[english]{}{1}{}
\normalsize

\newpage
\markboth{Introduction}{Introduction}
\Section{Introduction}
This paper is written in the style of \coma à la Bishop (\citealt*{Bi67,BB85,BR1987,MRR,CACM,Yen2015,ACMC}).

It is a natural sequel to the articles \citealt*{CLR01,KL00,KLP03,ALP08} and, to a lesser extent, \citealt*{ALN2021,CL2016,CL2016b,LM2022}.

\medskip  
Hensel’s lemma was introduced in mathematics in the context of  Hensel’s work on $p$-adic fields $\Qp$. In his work, valued fields are in fact ultrametric fields, i.e.\ fields $\gK$ given with a nonarchimedean absolute value $x\mapsto\abs x$, $\gK\to\RR$ (see Section 2). 

In 1930, \cite{Kru1930} and \cite{Deu1931}  introduced the notion of general valued fields, where ultrametric fields correspond to rank-one valued fields. The general notion of valued field was necessary to obtain Krull's fundamental theorem stating that the integral closure of a domain is the intersection of its overrings which are valuation domains.    

The abstract notion of an henselian ultrametric field (an ultrametric field satisfying Hensel's Lemma) along with the henselisation of an ultrametric field were introduced by Ostrowski in the seminal paper \citealt*{Ost1934}. Ostrowski introduced the henselisation of $({\gK},\abs. )$  as the separable closure  $\gKt$ of ${\gK}$ in its completion $\hat{{\gK}}$.
  See \citealt*{Roq02}  for more details on this topic.
   From a modern viewpoint, Ostrowski dealt with \cvals in the case of  rank-one valuations.

At a first glance one could consider the \LHm to be obvious for  ultrametric fields, insofar as a zero of a Newton polynomial system
can be calculated by Newton's method in a completion of the ultrametric field. However the fact that this zero of the Newton polynomial system  belongs to the valuation ring of the henselisation turns out to be difficult to prove, even for the most  simple discrete valuation rings.
Indeed,  reference  books in  \clama for the theory of \advs,  as \citealt*{Nagata62}, \citealt*{bourbaki72}, and \citealt*{EP2005}, do not pay any attention to the \LHm, even in the exercises.

%
%
%
%

 In the fifties the more general notion of \alohs was introduced by \cite{Azu1951} and \cite{Nag1953}, becoming afterwards a very important tool in Algebraic Geometry.
In this framework,  the \LHm states that on a \alo, a \sysN always has a zero with \coos in the henselisation of the ring.
The \LHm for local rings is usually proved as a consequence of the Grothendieck version of Zariski's Main Theorem (ZMT). This ZMT version  deals with a more general situation that is a priori much more difficult
(for a \cof treatment, see \citealt*{ACL2014}).

\smallskip The aim of this paper is to  give a direct \prco of the \LHm for ultrametric discrete fields, without using~ZMT\@. In classical mathematics, this provides an ``ad-hoc'' proof of the  Lemma for henselian rank-one valuations.

The plan is the following.

\smallskip In the first section we introduce the general framework of the \LHm for \alos. 
We recall the notions of Newton polynomial system, étale \syp and étale algebra.
Then we present the constructive versions of  structure \thos for \stf étale algebras over \cdis   
and the structure \tho for unramified \pf algebras over \cdis, as given in \citealt*[Theorems VI-1.7, \hbox{VI-1.9} and Corollary VI-6.15]{CACM}.

\smallskip 
In Section 2 we deal with ultrametric fields.

\setcounter{section}{2}
\setcounter{subsection}{2}

Let  $(\gK,\abs\cdot)$ be an ultrametric field with a nontrivial \vala. We  denote by~$\gKt$ the separable closure of~$\gK$ in its completion $\gKh$,
by $\gVt=\sootq{x\in\gKt}{\abs x \leq 1}$ its \gui{valuation ring},   and
let $\fmti=\sootq{x\in\gKt}{\abs x < 1}$. 

Our crucial result is the following \tho, where we compare  the henselisation $(\KHe,\VHe)$ constructed  in \citealt*{KL00} with~$(\gKt,\gVt)$. 

\setcounter{theorem}{4}

\begin{theorem}[two equivalent versions of the henselisation of an ultrametric discrete field] \label{thMRRKLintro}~\\
Let $(\gK,\abs\cdot)$ be an ultrametric discrete field.
The henselisation  $(\KHe,\VHe)$ is isomorphic to~$(\gKt,\gVt)$. More \prmt, there exists a unique $\gK$-\homo  $\KHe\to\gKt$ sending $\VHe$ into~$\gVt$, and this \homo is an  \iso.
\end{theorem}

We have got that an arbitrary \elt of $\gKt$ is \prmt the image of an \elt $\gamma$
in a field $\gK[\xi]\subseteq \KHe$, where 
$\xi$ is the special zero of a special polynomial.

Finally we prove our \LHm.

\setcounter{theorem}{5}

\begin{theorem}[\LHm  for an ultrametric discrete field] \label{LHM2intro}~\\ 
Let $(\gK,\abs\cdot)$ be an ultrametric discrete field and 
$(f_1,\dots,f_n)$ a \sysN at $(\uze)$ over $(\gV,\fm)$. This system admits a unique zero with  \coos in $\fmti$.  It admits also a unique zero with  \coos in $\fm\VHe$.
\end{theorem}

\setcounter{section}{0}
\setcounter{subsection}{0}
\setcounter{theorem}{0}

\section{Newton polynomial systems, étale \algs}\label{secEtales}

\Subsubsection{Constructive terminology}

In constructive mathematics, a \textsl{local ring} is defined as a ring where for all $x$, $x$ or $1-x$ is invertible (with an explicit \gui{or}). 

The \textsl{Jacobson radical of a ring}~$\gA$ is the ideal $\Rad(\gA)=\sotq{x\in\gA}{1+x\gA\subseteq \Ati}$. 

For a local ring, the Jacobson radical is its unique maximal ideal,\footnote{In \clama.} generally denoted by $\fm_\gA$ or $\fm$.
We shall simply say  \textsl{the local ring $(\gA,\fm)$}. 
 
An \textsl{Heyting field} is a nontrivial local ring whose Jacobson radical is $0$. The \textsl{residual field} of~a nontrivial local ring $(\gA,\fm)$ is the Heyting field $\gA/\fm$ also denoted by~\(\ov\gA\).

A \textsl{discrete field} is a nontrivial ring in which any element is zero or invertible. It is the same thing as an Heyting field with a zero test. 

The local ring is said to be \textsl{residually discrete} if its residual field is discrete. This amounts to saying we have explicitly the disjunction $x\in\Ati$ or $x\in\fm_\gA$ for all $x$ in $\gA$. 

For local rings $(\gA,\fmA)$ and $(\gB,\fmB)$ a ring morphism $\varphi\colon\gA\to\gB$ is said to be \textsl{local} when $\varphi^{-1}(\gB^\times)\subseteq \Ati$. In this case we say that $(\gB,\fmB)$ is an $(\gA,\fmA)$-\alg.

\Subsubsection{Hensel codes over local rings}

An \textsl{\cdH} over a local ring  $(\gA,\fm)$ is a pair $(f,a)\in \AX\times \gA$ where $f$ is a \polu, $ f(a)\in \fm $ and $f'(a)\in \Ati$ (in other words $\ov a$ is a simple zero of $\ov f$). In this case we say that  $f$ is an \textsl{\polH} or a \textsl{Nagata polynomial}. 

A \textsl{special \pol} is a \pol 
$h(X)=X^n-X^{n-1}+\sum_{k=0}^{n-2} a_kX^k$ with the $a_k$'s $\in\fm$. In this case $(h,1)$ is an \cdH.

An \textsl{\zeH for the \cdH $(f,a)$} in an $(\gA,\fmA )$-\alg $(\gB,\fmB)$, is an \elt $\xi\in\gB$ such that $\varphi_\star(f)(\xi)=0$ and $\xi-\varphi(a)\in\fmB$. 

A local  ring  $(\gA, \fm)$ is said to be \textsl{henselian} if any Hensel code has an \textsl{Hensel zero} $\alpha$ in $(\gA, \fm)$.
In this case we say that \textsl{the zero $\alpha$ of~$f$ is the lifting in~$\gA$ of  the simple zero $a$} of $\ov f$ in $\ov\gA$. The \zeH for the \cdH $(h,1)$ of a special \pol $h$ is called the \textsl{special zero} of the \pol.

An \zeH for an \cdH $(f,a)$ is \ncrt unique, regardless of the hypothesis that $f$ be monic: one writes 
\[0=f(\alpha')-f(\alpha)=f'(\alpha)\mu+b\mu^2=\mu\cdot(f'(\alpha)+b\mu) \quad \hbox{with }b\in\gA,\eqno(+)\]
where $\mu=\alpha'-\alpha\in\fm$, and as $f'(\alpha)+b\mu\in \Ati$ we get $\mu=0$.

The \textsl{completion} of  $\gA$ for the $\fm$-adic topology, \cad the projective limit of  $(\gA/\fm^k)_{k\in\NN}$, is denoted by $\gAh$.
We have a natural morphism $\varphi\colon\gA\to\gAh$ and we let $\fmh:=\varphi(\fm)\gAh$.   
The morphism~$\varphi$ is injective \ssi $\bigcap_{k\in\NN}\fm^k=0$.
If $(\gA,\fm)$ is a \dcd local ring,  $(\gAh,\fmh)$ is a local ring and it is henselian. 
Indeed, Newton's process (as in \thref{thNewtonQuad}) allows us to compute an \zeH for any \cdH.

\Subsubsection{Henselisation of a local ring}

The notion of henselisation of a local ring  corresponds to the solution of the \uvl \pb  for the full subcategory of henselian  local rings  in the category of local rings and local homomorphisms.

For a \dcd \alo the paper \citealt*{ALP08} constructs this henselisation by iteratively adding formal \zeHs for \polHs. These building blocks are \textsl{Nagata extensions} of $\gA$ \citep*[Definition~6.1]{ALP08}. 
We denote by $\Ahe$ the \textsl{henselisation} of $(\gA,\fmA)$.

Our first result explains how to lift a  smooth residual zero of a not \ncrt monic (viz.\ not Hensel) \pol on $\gA$ to an henselian extension. It is stated for \alrds in \citealt*[Lemma~5.3 and Proposition~5.4]{ALP08}.
Here we state it in the context of a local ring  $(\gA,\fmA)$, a \pol $f\in\AX$ 
and a local morphism from $(\gA,\fmA)$ to an henselian local ring $(\gB,\fmB)$.
\begin{lemma}[Hervé's trick] \label{lemTrick1} Let $(\gA,\fmA)$ be a local ring, $\varphi\colon\gA\to\gB$ a local  morphism with $(\gB,\fmB)$ henselian, 
and $f(X)=\sum_{k=0}^na_kX^k\in\AX$ with $a_0\in\fmA$ and $a_1\in\Ati$.
The \pol~$f$ has a zero  $\gamma$ in $a_0\gB\subseteq \fmB$, with $f'(\gamma)$ invertible in $\gB$. It is the unique zero of~$f$
in~$\fmB$.
In particular, if  $(\gA,\fmA)$ is henselian,~$f$ admits a zero in $a_0\gA\subseteq \fmA$, and  it is the unique zero of~$f$ in $\fmA$.
\end{lemma}
%
\begin{proof} We define the special \pol 
\[
\formu{  
    g(X)&=&X^n-X^{n-1}+a_{0}\cdot
\big(\som_{j=2}^{n}(-1)^j   a_{j}a_{0}^{j-2} a_{1}^{-j}  X^{n-j}\big) \\[.4em]
&=& X^n-X^{n-1}+a_{0}\ell(X)  \quad \text{with } \ell(X)\in\AX.
  } 
\]
The \flw \egt is correct in $\gA[X,1/X]$
\[  
  a_0g(X) = X^nf\left(\frac{-a_0a_1^{-1}}X\right)\text.  \eqno (*)
  \]  
Let $\delta=1+\alpha$ (where  $\alpha\in\fmB$) be the special zero of the special \pol $g$. Then $\delta\in
{\gB}^\times$.  
Let $\gamma={\frac{-a_0  a_1^{-1}}\delta}=-a_0  (a_1\delta)^{-1}\in\fmB$. Applying~$(*)$ we see that
$-a_0  g(\delta)=\delta^n  f(\gamma)$, so $f(\gamma)=0$. Moreover $f'(\gamma)\in{\gB}^\times$ because $f'(0)\in\Ati$ and $\gamma\in\fmB$.

\noindent Uniqueness is already proved, see $(+)$.
\end{proof}
%

\begin{remark}
In the case of a \alrd $(\gA,\fmA)$ with 
henselisation $\Ahe$, the \elt $\delta$ in the previous proof appears in the construction of the henselisation $\Ahe$ as an \elt of 
\[\gA_g:=S^{-1}\Ax \text{  where }\Ax=\aqo \AX g\text{  and } S=\sotq{s(x)\in\Ax}{s(1)\in\Ati}. 
\]
The \ri $\gA_g$ is a \alrd, \fpt over $\gA$,
and $\Rad(\gA_g)=\fm\gA_g$. So~$\delta$ is the image of $x\in\Ax$ in $\gA_g$ via the localisation morphism.
As the canonical morphism $\gA\to\gA_g$ is injective, we can identify $\gA$ with a subring of $\gA_g$. 
These rings $\gA_g$ are the \elr building blocks in the  construction of the henselisation of the \alrd $(\gA,\fmA)$ given in \citealt*{ALP08}. 
\eoe
\end{remark}
\subsection{Newton \poll \sys in a local ring}

Given  $(\gA,\fm)$ a local ring, a \textsl{\sysN} (or \textsl{Hensel \syp)  at the point $(\ua)=(\an)\in\gA^n$}  is given by a \syp $(\uf)=(f_1,\dots,f_n)$ in $\AXn^n$  when  $(\ua)$ is   an
\textsl{approximate simple zero} in the following sense:
\begin{itemize}
\item  
the $f_j(\ua)$ are in~$\fm$;
\item  the Jacobian matrix of the \sys is invertible at $(\ua)$ modulo $\fm$.

%
\end{itemize}

The second condition amounts to saying that the Jacobian \deter $\Jac(\ua)$ is invertible in $\gA$ modulo $\fm$.\footnote{Note that because \((\gA,\fm)\) is a local ring, invertibility and invertibility modulo~\(\fm\) coincide.}

The \textsl{\LHm} says  that a \sysN  at $(\ua)$ on a local ring  $(\gA,\fm)$ has  a zero $(\uxi)$ with \coos in the henselisation~$\Ahe$ with $\xi_j-a_j\in\fm\Ahe$ for each~$j$.


By \dfn, an \cdH is a one-variable \sysN. 

\smallskip We remark here  that  \cite{Lafon63} and the \citetalias[Section 15.11]{stacks-project},  
give the \LHm in a slightly hidden form for  henselian \cous (see e.g.\ implication $(5)\Rightarrow(2)$ in Lemma 15.11.6 in the \citetalias{stacks-project}) but they use  ZMT for their \demo.

\Subsubsection{Uniqueness of the \zeH}

\begin{lemma} \label{lem-Unicite-zeH}
Let $(\gA,\fm)$ be a local ring and $(\uf)=(f_1,\dots,f_n)$ a \sysN at $(\ua)\in\gA^n$.  If $(\ual)$ and $(\uga)$ are \zeHs at $(\ua)$ for this \sys, then $(\ual)=(\uga)$.      
\end{lemma}
\begin{proof} Let $J(\uX)$ be the Jacobian  matrix of the \sys.  Write $(\uga)=(\ual)+(\udel)$ with the $\delta_i$'s in~$\fm$. 
We see $(\uf(\ual))=(f_i(\ual))_{i\in\lrbn}$, $(\uf(\uga))$, and $(\udel)$ as column vectors.
Taylor  formulas in several variables for \pols yield an \egt
\[(\uf(\uga))=(\uf(\ual))+J(\ual)(\udel)+ M(\udel),
\]
where  $M$ is a square matrix with \coes in the \id $\gen{{\udel}}$ and 
$(\uf(\uga))=(\uf(\ual))=(\uze)$. So $(\rI_n+J(\ual)^{-1}M)(\udel)=(\uze)$, and
$(\udel)=(\uze)$.
\end{proof}

The image of a \sysN by a local morphism  is a \sysN.
The uniqueness of a zero (if it exists) in an $(\gA,\fm)$-algebra is proved in the same way.

\subsection{Newton's method}

The following \tho describes the so-called quadratic Newton method in
a purely \agq context.

This result holds true also for  a general ideal $\fa$ of any ring $\gA$ (instead of $\fm_\gA$ for a \alo), as explained in  \citet*[Theorem III-10.3]{CACM}.  In many cases it is  relevant to consider the ideal $\fa$  generated by the $f_i(\ua)$'s. 
\begin{theorem}[quadratic Newton method] \label{thNewtonQuad}~ \\
Let $(\gA,\fm)$ be a local ring, and 
 $(\uf)=(f_1,\dots,f_n)$ a \sysN in  $\AXn$ at $(\ua)=(a_1,\ldots ,a_n)\in\gA^n$. We denote by $\JAC(\uxi)$ the Jacobian matrix of $(\uf)$ at $(\uxi)$. Let $U$ be an inverse of $\JAC(\ua)$ modulo $\fm$.  
We define sequences $(\ua^{(m)})_{m\geq
0}$ in $\gA^n$  and $(U^{(m)})_{m\geq 0}$ in $\Mn(\gA)$
by the following iteration:
$$\begin{array}{lcl}
\ua^{(0)}=\ua,&  \quad \quad  & \ua^{(m+1)}=\ua^{(m)}- U^{(m)} \cdot
\uf(\ua^{(m)}),    \\[1mm]
U^{(0)}=U,&   & U^{(m+1)}=U^{(m)}\,\left(2\I_n-\JAC(\ua^{(m+1)})U^{(m)}\right).
\end{array}$$
Then for each integer $m$ we get the following congruences:
$$
\begin{array}{lcll}
 \ua^{(m+1)}\equiv\ua^{(m)} & \,\,\text{and}\,\, &
  U^{(m+1)}\equiv U^{(m)} &\,\mod \,\fm^{2^m}  ,  \\[1mm]
  \uf(\ua^{(m)})\equiv 0 &  \,\,\text{and}\,\,  &
   U^{(m)}\,\JAC(\ua^{(m)})\equiv \rI_n  &\,\mod \,\fm^{2^m}.
\end{array}
$$
\end{theorem}

This \tho says that if we have a \sysN   $(\uf)=(f_1,\dots,f_n)$ in $\AXn$ with $(\ua)$ a 1-approximate simple zero in the local ring $(\gA,\fm)$, we can find  ~$(\ua^{(m)})$ which is  a $2^m$-approximate simple zero of $(\uf)$,  with $\ua^{(m)}\equiv\ua \mod \,\fm$ ($m> 1$).


\smallskip We shall use the following terminology in \coma: a local ring  $(\gA,\fm)$  
is said to be \textsl{quasi-\noe} if $\gA$ is a \coh \fdi \ri, $\fm$ is a \itf, and 
$\bigcap_{n\in\NN}\fm^n=0$.\footnote{In some \prcos, it will be important that, for any $x\neq 0$, the integer $k$ such that $x\in\fm^k\setminus\fm^{k+1}$ be known.} In this context, each $\fm^n$ is a \coh \fdi \Amo and the  natural morphism $\gA\to\gAh$ is injective.
This happens in the following case:  $\gA=\gB_{1+\fm_\gB}$ where $\gB$ is a  \pf \alg over~$\ZZ$  or over a \cdi, $\fm_\gB$ is a \itf, and~$\fm_\gA=\fm_\gB\gA$.\footnote{A quasi-\noe \colo is said to be \noe if the ring $\gA$ is \noe. We shall not use this notion.}

\smallskip Applying Newton's method in two contexts, namely, 
to a quasi-noetherian local ring, and to the  valuation ring of an ultrametric discrete field, we can get what we call ``two weak forms of the \LHm'', where the \coos of the zeroes of the \pol \sys are not asserted to be in the henselisation of the ring:

\begin{corollary}[\LHm, first weak form] \label{LHM1} ~\\
Let $(\gK,\abs \cdot )$ be an ultrametric discrete field\ \footnote{See Section~\ref{seccvdu}.} and $(f_1,\dots,f_n)$ a \syp as in   \thref{thNewtonQuad}, with 
$\gA=\gV=\sotq{x\in\gK}{\abs x \leq 1}$ and $\fm=\sotq{x\in\gK}{\abs x < 1}$. Then the system has a unique zero $(\uxi)=(\xi_1,\dots,\xi_n)$ with \coos in $\gKh$ satisfying $\xi_i-a_i\in \fmh$ for $i=1,\dots,n$.  
\end{corollary}

\begin{corollary}[\LHm, second weak form] \label{LHMquasinoet} ~\\
Let $(\gA,\fm )$ be a quasi-noetherian local ring and $(f_1,\dots,f_n)$ a \syp as in \thref{thNewtonQuad}. Then the \sys has a unique zero $(\uxi)=(\xi_1,\dots,\xi_n)$ with \coos in $\gAh$ satisfying 
$\xi_i-a_i\in \fm\gAh$ for $i=1,\dots,n$.
\end{corollary}

\subsection{Étale \algs
 }

The context of Newton's method can be generalised and formalised under the name of basic étale \alg.

\begin{definition}[étale polynomial system, étale \Alg] \label{defiBasicEtale} 
~\\
Let~$\gA$ be a commutative \ri.
\begin{enumerate}
\item Let
 $(\uf)= (f_1,\dots,f_n)$ be a \sys of $n$ \pols in  $\AXn$ and 
 \hbox{$\gB=\aQo\AuX \uf$}. 
 The \Alg $\gB=\Axn$ is said  to be \textsl{basic étale} if the Jacobian matrix $\JAC({\ux})$ of the \sys $(\uf)$  is
\iv in~$\gB$. In this case we say that the \syp $(f_1,\dots,f_n)$ is an \textsl{\sype}.
\item A  \pf \Alg $\gC=\aqo \AXm {g_1,\dots,g_s}$ is said to be \textsl{étale} if we know \eco $(u_i)_{i\in I}$ in $\gC$ such that each \Alg $\gC[1/u_i]$ is basic étale for a convenient \pn \((\underline{f_i})=(f_{i,1},\dots,f_{i,n_i})\).
%
%
\end{enumerate}
\end{definition}

The notion of \alge is a fundamental concept of  Commutative Algebra, introduced by Grothendieck. 
One  result of that theory asserts that any \alge is 
basic étale for a convenient \pn. 

Let us remark that  a  finite product of étale \algs is étale 
and a localisation $\gC[1/s]$ of an étale \alg $\gC$ is étale. In particular the trivial \alg is étale.


\begin{lemma} \label{note0}
It is always possible to replace a \sysN $((f_1,\dots,f_n),(\uze))$ over a local ring $(\gA,\fm)$ 
with a basic étale \sysN \(((f_1,\dots,f_{n+1}),(\uze,0))\), i.e.\ with a \sysN over $(\gA,\fm)$ which is a basic étale \syp over $\gA$. Moreover for the new \syp
the notion of  \zeH  is unchanged.
\end{lemma}

\begin{proof}
We denote by $\Jac(\uX)$ the Jacobian \deter, we add an \idtr $X_{n+1}$ and the \pol $f_{n+1}:=(1+X_{n+1})\Jac(\uze)^{-1}\Jac(\uX)-1$.
We have $f_{n+1}(\uze,0)=0$. The new \syp has its Jacobian \deter~$\Jac_1(\ux,x_{n+1})$ \iv in the new quotient \Alg, and an \zeH~$(\uxi)$ of the first \sys with \coos in an \Alg gives the \zeH $(\uxi,\eta)$  for the new \sysN with $(1+\eta)\* \Jac(\uze)^{-1}\Jac(\uxi)=1$, \cad 
$\eta=\Jac(\uze)/\Jac(\uxi)-1$.
\end{proof}

\smallskip This \elr manipulation gives some \prts (of the \coos) of the zero. For example, as a consequence of \thref{thNewtonCorpsDiscret}, if $(\gA,\fm)$ is an integral \alo with~\hbox{$\gK=\Frac\gA$}, the \coos of the \zeH of a \sysN over $(\gA,\fm)$ in a \Klg are always \spl over $\gK$.


\Subsubsection{Structure of étale \algs over a \cdi} 

Let  $\gK$ be a \cdi and $\gKp$ a field extension which is a \tf \Kev. 
We say in this case that $\gKp$ is a \textsl{finite algebraic extension} of $\gK$. 

More \gnlt a \Klg $\gC$ which is a \tf \Kev is said to be \textsl{finite over $\gK$}. The \elts of $\gC$ are \agq over $\gK$, but perhaps we don't know the dimension of~$\gC$ as \Kev.
If $\gC$ is generated by~$n$ \elts  as \Kev, we write $[\gC:\gK]\leq n$.
 If $\gC$ contains $m$ $\gK$-\lint independent \elts, we write $[\gC:\gK]\geq m$.
Finally  $\gC$ is said to be 
 \textsl{\stf} over $\gK$ when the dimension of $\gC$ as \Kev is known, and we write $[\gC:\gK]= n$. In this case,
\begin{itemize}
\item we know how to compute the \polmin over $\gK$ of each \elt of $\gC$;
\item  any intermediate finite \Klg $\gD$ is \stf over $\gK$;
\item  if moreover $\gD$ is a field, $\gC$ is \stf over $\gD$ and we get the usual formula $[\gC:\gK]= [\gC:\gD]\* [\gD:\gK]$.
\end{itemize}

\begin{definition}[\ste \algs over a \cdi] \label{defiaste}~\\ 
Let $\gK$ be a \cdi. A \stf \Klg $\gB$ is  \textsl{\ste} if the trace  form $\phi(x,y)=\Tr_{\gB/\gK}(xy)\colon\gB\times \gB\to\gK$ is 
\textsl{nondegenerate}, \cad letting $\varphi (x):=\phi(x,\bullet)$, the \Kli $\varphi$ defines an \iso from the \Kev $\gB$ onto its dual.\footnote{These \dfns may be generalised to \algs over an arbitrary commutative ring \cite[Theorem VI-5.5]{CACM}.}  
\end{definition}

We give a first structure \tho for \ste  \Klgs.
A \prco is given in \citealt*[Theorems VI-1.7 and VI-1.9]{CACM}. 
It shows in particular that \ste \algs over a \cdi are \stf étale \algs.
Moreover \thref{thNewtonCorpsDiscret} will show that étale \algs over a \cdi are \ste \algs. 

Note that this \cov \thref{thstrucste} is more precise than the classical ones, and the \prco is rather subtle. In fact hypotheses are given in a form allowing to obtain an \algo for the conclusions.  

\begin{theorem}[primitive element theorem] \label{thstrucste} ~\\
Let $\gK$  be a \cdi and $\gB$ a  \stf \Klg. 
\begin{enumerate}
\item 
The following are equivalent.
\begin{enumerate}
\item  $\gB$ is \ste.
\item $\gB$ is generated by \elts which are \spl
over $\gK$.
\item  All \elts of $\gB$ are \spl
over $\gK$.
\item $\gB$ is isomorphic to a finite product of \Klgs $\aqo{\KX}{h_i}$ with  \spl \polus
 $h_i$.
\end{enumerate}
In particular \ste \Klgs  are étale.
\item  When $\gB$ is a \cdi or $\gK$ is infinite, properties of Item 1 are also equivalent to: 
\begin{enumerate}\setcounter{enumii}{4}
\item $\gB$ is isomorphic to a \Klg $\gK[\zeta]=\aqo \KZ g$, where 
$g$ is a \spl \polu in~$\KZ$.
\end{enumerate}
\end{enumerate}
\end{theorem}

In the last case, if $g$ is factorised as $g=g_1\cdots g_r$, we have a canonical \iso  $\gB\simeq \prod_{j=1}^r\aqo\KZ{g_j}$.

\smallskip
An \Alg $\gB$  is said to be \textsl{unramified} (or \textsl{neat}) if it is \pf and its \textsl{module of (K\"ahler) \diles} reduces to zero. The \Bmo of \diles is isomorphic to the cokernel of the transpose of the Jacobian matrix (seen in $\gB$). In other words, the module of \diles reduces to zero \ssi the transpose of the Jacobian matrix is surjective. 
It is clear that an \alge over an arbitrary commutative ring is unramified.

The following important \tho gives a strong converse to  
\thref{thstrucste}. Concerning the notion of simple isolated zero in this \tho, see the \cof approach in \citealt*[Section IX-4]{CACM}.

\begin{theorem}[unramified \alg over a \cdi] \label{thNewtonCorpsDiscret}~\\
Over a \cdi $\gK$ any unramified \Klg  is \stf, étale, \ste.  
In particular, for an étale \syp, all zeroes of the corresponding variety in an \ac overfield are isolated, simple, with \coos  \spl
over $\gK$.
\end{theorem}

A \prco is in \citealt*[Corollary~\hbox{VI-6.15}]{CACM}. In the second French edition \citealt*{ACMC}, a more elementary \prco is given at the end of Section VI-6.

\smallskip We now give a more precise description of the situation. 

\newcommand{\Kual}{\gK[\ual]}
\begin{descri}[étale \syp over a \cdi, precisions] \label{descriEtaleCdi}~\\ 
Let
 $(\uf)= (f_1,\dots,f_n)$ be an étale \syp  in  $\KXn$ over an infinite \cdi  $\gK$\footnote{When $\gK$ is finite, or more \gnlt when we don't know whether it is infinite, slight modifications are to be introduced in this description, using Item \textsl{1.(d)} of \thref{thstrucste}. A valued field with a nontrivial valuation is always infinite.}  
 and let  
 \[\gD=\Aqo\KuX \uf=\Kxn\] 
be the quotient \Klg. 
\begin{itemize}
\item We can construct a primitive \elt $\zeta$ of $\gD$. Its \polmin $g$ over $\gK$ is \spl, so $\gD=\gK[z]\simeq \aqo{\KZ}{g(Z)}$. 
So we have \pols $q_i\in\KZ$ such that $x_i=q_i(z)$ in~$\gD$.
\item  A zero $(\ual)=(\alpha_1,\dots,\alpha_n)$ of the \syp in a \Klg  $\gC$ gives a $\gK$-morphism $\varphi\colon\gD\to\Kual \subseteq \gC$ satisfying $\varphi(\ux)=(\ual)$.
The \alg $\Kual$ is then isomorphic to
a quotient of~$\gD$.
If $\Kual$ is connected and  nontrivial, it is a \cdi, for it is \zedr \citep[Fact IV-8.8]{CACM}. 

\item  In the following we assume  $\Kual$ to be connected and nontrivial.
\begin{itemize}
\item [$\bullet$] If we know a \fcn of $g$ as a product of $r$ \ird \pols~$g_j$ over~$\gK$, then $\gD\simeq \gL_1\times \dots\times \gL_r$ with $\gL_j\simeq \aqo \gD{g_j(z)}$, and we get a corresponding \sfio $(e_1(z),\dots,e_r(z))$ in $\gD$.\footnote{The \id $\gen{g_j(z)}$ is generated by the \idm $1-e_j(z)$. The  case $r=0$ remains possible. An \sype may be impossible.} And $\gK[\ual]$ is isomorphic to one of the \cdis $\gL_j$ via $\varphi$, with $\varphi(e_j(z))=1$.
\item [$\bullet$] Otherwise if $\gL$ is a \spb extension of $\gK$
and $g$ is completely factorised over~$\gL$,
the quotient \alg seen over~$\gL$, \cad $\gL\otimes_\gK\gD\simeq\aqo{\gL[Z]}{g(Z)} $, is isomorphic to~$\gL^{d}$, \hbox{where $d=\deg(g)$}. So the \syp has exactly $d$ zeroes with \coos in~$\gL$.
If we embed $\gL\subseteq\Ksep\subseteq\Kac$ we get in this way all the zeroes with \coos in~$\Kac$.\footnote{Here \(\Ksep\) and \(\Kac\) are respectively a separable and an algebraic closure of~\(\gK\).}
\eoe 
\end{itemize}
\end{itemize}
\end{descri}

What happens in the more general situation where we don't know a factorisation of~$g$ over~$\gK$? This case is not very different from the preceding one. Indeed, if at a certain step we have found an \idm $e\neq 0,1$ in $\gD$,\footnote{This happens each time an \elt $\neq 0$ of $\gD$ is not invertible, \cad when its \polmin has degree $>1$ and its constant \coe is zero.} it is equal to $1$ or $0$ in $\Kual$ and we can replace~$\gD$ by $\gD[\frac 1 e]$ or~$\gD[\frac 1{1-e}]$ (this amounts to replacing $g$ with a strict divisor). The new \alg remains \ste. All previous computations remain valid, and the new version of~$g$ is better. The possible improvements of this type are bounded in number by $\deg(g)$. In consequence most of the concrete results obtained when we assume that we know a factorisation of~$g$ remain valid without this hypothesis. We are here at the heart of the dynamical method in algebra.    

\section{The case of ultrametric discrete fields}\label{seccvdu}

\subsection{Definitions}

\Subsubsection{Valued discrete fields}

First we recall the  \cov \dfn of a \textsl{\cvd} $\KV$ as in \citealt{KL00} or \citealt*{CLR01}: 
\begin{itemize}
\item $\gK$  is a \cdi;
\item $\gV$ is a subring of $\gK$;
\item for all $x\in\gK\eti$ we have  $x$ or $1/x\in\gV$; 
\item  \dve in $\gV$ is explicit.%
\footnote{This means that for $x,y\in\gV$ we have a test for the existence of a $z\in\gV$ such that $yz=x$. This amounts to saying that $\gV$ is a detachable subring of $\gK$. This is aways true in \clama by the Law of Excluded Middle.}
\end{itemize}

In this case~$\gV$ is an  \icl \alrd. 

This definition is equivalent in \clama to the usual one. We have added decidability hypotheses in order to facilitate computations.

\begin{definition}[henselisation of a \cvd] ~ \\
Let $\KV$ be a \cvd.
\begin{itemize}

\item  An \textsl{extension} of $\KV$ is a \cvd $(\gL,\gW)$ 
together with an homomorphism $\phi\colon\gK \rightarrow \gL$ such that $\gV = \gK \cap  \phi^{-1}(\gW)$

\item  An \textsl{henselisation of $\KV $} is an extension which is an henselian valued discrete field $(\KHe,\VHe)$ 
and 
such that the corresponding homomorphism
$\phi^H\colon\gK \rightarrow \KHe$ factorises in a unique way through every
extension of $(\gK, \gV)$ which is an henselian valued discrete field.
\end{itemize}
\end{definition}

Being the solution of a universal \pb, an henselisation of a \cvd is unique up to unique \iso. 

\cite{KL00}  construct the henselisation of a \cvd $\KV$%
, denoted by $(\KHe,\VHe)$. It is obtained by adding successively \zeHs of \polHs. 

Indeed, given $ f = X^n + a_{n-1}X^{n-1} + \cdots + a_0 $ a Nagata polynomial in $\gV[X]$,
the authors describe explicitly an extension $(\gK[\alpha], \gV_{\alpha})$ of $\KV$ for which the image of $f$ in $\gK [\alpha][X]$
has an henselian zero $\alpha$, and such that the extension map $\gK \rightarrow  \gK[\alpha]$ factorises
in a unique way every extension of $\KV$ to $(\gL,\gW)$ if the image of $f$
in $\gL[X]$ has an henselian zero. Furthermore, the residue field and the value group
of $(\gK[\alpha], \gV_{\alpha}) $ are canonically isomorphic to the residue field and the value group
of $\KV$ respectively. This explicit construction is based on the study of the \textsl{Newton polygon} of $f$.
Note that we do not assume to know if the base field contains a special zero of $f$.
So, the finite extension $\gK[\alpha]$ which is constructed is a \cdi but it is not \ncrt a \stf extension of $\gK$. The construction of the henselisation of a \cvd is   
 very similar to the construction of the real closure of \codi, which works even when we are not able to decide if an arbitrary \pol has a zero in the base field.

\begin{remark}
  We point out that, even though the construction of the henselisation of a valued field  given in \cite{KL00} 
is  very similar  to the construction of the henselisation of a local ring, the tools used in the two cases are  different. Indeed, the first one is  entirely in the framework of \cvds, and it is a priori less \gnl than the construction in the framework of \alrds given in \citealt*{ALP08}. 
So, a priori we should use two distinct notations for these ``two henselisations'', namely $\Vhe$ (henselisation as local ring), and $\VHe$ for the valuation ring of the henselisation of the valued field. Although the nonobvious fact that they coincide seems to be accepted, we have not found a proof in the literature. 
\eoe
\end{remark}

\Subsubsection{Ultrametric fields}

In the book \citealt*{MRR}, the theory of \valas is treated \cot using  the following \dfn which is the usual one in \clama for fields with an absolute value.\footnote{In fact, they start with an Heyting field $\gK$, but the two definitions are clearly equivalent. In \clama one starts also usually with $\gK$ a \cdi.} 

\begin{definition}[\fva, \cvu] \label{defivala}~
\begin{enumerate}
\item 
An \textsl{\vala} over a ring $\gK$ is a function
$\gK\to \RR^{\geq 0}, \,x \mt |x|$ satisfying the following  \prts.
\begin{itemize}\itemsep=.1em
\item $|x| = 0$ \ssi $x=0$;
\item $|x| > 0$ \ssi $x$ is invertible;
\item $|xy| = |x||y|$;
\item $|x + y| \leq |x| + |y|$.
\end{itemize}
\vspace{.1em}
\noindent One says that $(\gK,\abs\cdot)$ is a 
\textsl{field with an absolute value}.\footnote{The book \citealt*{MRR} uses  \gui{\cval}, as very often in the english literature, but this is in conflict with our  terminology for general \cvds, which follows Krull and Bourbaki.}
%
\item 
The \vala is said to be \textsl{ultrametric} if $|x + y| \leq \sup(|x|,|y|)$ for all $x,y$.\footnote{Note that from a \cov viewpoint, $z=\sup(x,y)$ is well defined for real numbers, but it cannot be proved that $z=x$ or $z=y$ with an explicit \gui{or}.}
In this case we speak of an \textsl{\ultm field}. Note that if $\abs x<\abs y$ then $\abs{x+y}=\abs y$.
\item 
An \valu defines a \gui{\adv} $\gV$: 
\[
\formul{
\gV\hphantom{\eti}&:=&\sotq{x\in\gK}{|x|\leq 1} ,\text{ with} 
\\ \gV\eti&:=&\sotq{x\in\gK}{|x|= 1}\text{ and }
\\ \fm\hphantom{\eti}&:=&\Rad(\gV)=\sotq{x\in\gK}{|x|< 1}.
}
\]
\item 
Two  nontrivial \valus on $\gK$ defining the same \adv are said to be \textsl{\eqv}, and each is a positive power of the other \citep*[see][Theorem XII-1.2]{MRR}.
\item 
The distance $d(x,y)=|x-y|$ makes $\gK$ a metric space, whose completion is denoted by~$\gKh$. The \vala extends uniquely to $\gKh$, and $(\gKh,\abs{\cdot})$ is also an \ultm field\footnote{In the archimedean case, $(\gKh,\abs{\cdot})$ is also a field with an absolute value.}. The completions of $\gV$ and $\fm$ are denoted by $\gVh$ and $\fmh$.
The image of $\gK\eti$ in $(\RR^{>0},\times )$ is the \textsl{value group} of $(\gK,\abs\cdot)$.

\end{enumerate}
\end{definition}

A \cvar $\gK$ is a nontrivial \alo whose Jacobson radical is reduced to $0$, \cad an \textsl{Heyting field}.
\footnote{The ring $\gK$ is not \ncrt a \cdi. So we have given the \dfn for a ring $\gK$. This avoids to recall first the \dfn of a Heyting field.}
 
The \dfn does not require $\gK$ to be a \cdi. In \gnl the completion~$\gKh$ is not discrete, even when $\gK$ is discrete.

\begin{remark} \label{remdefivala} 
In Item 4  we have put \gui{\adv} in quotes because $\gV$ needs not, from the \cof viewpoint, be a \alo, nor discrete, nor \dcd.
\end{remark}

\smallskip Let us consider an ultrametric field $(\gK,\abs\cdot)$.

The Heyting field  $\gK$ is discrete \ssi for all $x\in\gK$ we have the disjunction \hbox{\gui{$\abs x = 0 \vuu \abs x >0$}}. This amounts to saying that  $\gV$ is an integral domain. 

The \cou $\KV$ is a \cvd in the \cof meaning if moreover the disjunction \gui{$\abs x = 1 \vuu \abs x <1 \vuu \abs x >1$} is valid for all $x\in\gK$.

This amounts to saying that $\gV$ is an integral \alrd.
In this case we say that $(\gK,\abs\cdot)$ is an \textsl{ultrametric discrete field}.
  
\Subsubsection{Translation in terms of valuations}

Let us consider the map \(\ell\colon(\RR^{\geq 0},\times)\to(\RR\cup\sing{+\infty},+)\) defined by \(\ell(r)=-\log(r)\) for $r\neq 0$ and by $\ell(0)=+\infty$; endow $\RR\cup\sing{+\infty}$ with the topology that
makes  $\ell$ an homeomorphism. 
For an ultrametric field we define the \valn $v\colon\gK\to (\RR\cup\sing{+\infty},+)$ by $v(x)=\ell(\abs x)$. We simply translate the \prts of $({x\mapsto \abs x},\gK\to \RR^{\geq 0})$  into \prts of $v$, reversing the order relation, replacing multiplication with addition and  $\sup$ with $\inf$.
This gives the following \prts:
\begin{itemize}\itemsep=.1em
\item $ v(x) = \infty$ \ssi $x=0$;
\item $v(x) \neq  \infty$ \ssi $x$ is invertible;
\item $v(xy) = v(x)+v(y)$;
\item $v(x + y) \geq \inf(v(x),v(y))$, with \egt if $v(x)\neq v(y)$;
\item $ v(x) \geq 0$ \ssi $x\in \gV$;
\item $ v(x) > 0$ \ssi $x\in \fm$.
%
%
\end{itemize}

When replacing the \vala with an \eqv \vala, the valuation $v$ is simply multiplied by a constant $r>0$.

\smallskip For an ultrametric discrete field, $\gK$ and the residual ring $\gV/\fm$ are \cdis, the subgroup $\sotq{\abs x}{x\in\gK\eti}$ is a discrete multiplicative subgroup  $\Delta$ of  $\RR^{\geq 0}$, $\Gamma=\sotq{v (x)}{x\in\gK\eti}$ is a discrete additive subgroup of $(\RR,+)$, and the union
$\Gamma\cup\sing{+\infty}$ is a disjoint union: the topology is discrete. 

In this case, we have the following useful result, also valid for \gnl \cvds:
\begin{itemize}
\item if $\sum_{i=1}^nx_i=0$, with  the $x_i$'s not all zero, the infimum of the $v(x_i)$'s is attained at least for two  distinct $i$.
\end{itemize}

\Subsubsection{Three basic examples}

In \coma we define a  \textsl{\advd} (a DVR in short)  as an integral domain  $\gV$ (with \cdf $\gK$) in which we give an  \elt  $\pi\neq 0$ (called a  \textsl{regular parameter}) such that each \elt of $\Vtl$ is written as $a=u\pi^k$ with $u\in\gV\eti$ and~$k\in\NN$. This makes $\KV$ a \cvd with \valn $v(a)=k$.
Letting $\abS {u\pi^k}= e^{-k}$ for a fixed real number $e>0$, $(\gK,\abs\cdot)$ is an ultrametric discrete field.   

Three basic examples are now given. In Examples 2 and 3, the \vala is not in $\RR^{\geq 0}$ but in a submonoid of 
$(\gK,\times)$ isomorphic~to the closure of $\sotq{1/2^n}{n\in\NN}$ 
in $\RR^{\geq 0}$. 

\begin{enumerate}
\item Here $\gK=\QQ$, the standard $p$-adic \vala is $|r|_p=p^{-k}$ if~$r=\frac m n\,p^k$ with $m$ and $n\in\ZZ$ coprime with $p$. The corresponding \advd is $\gV=\ZZ_{1+p\ZZ}$ (the ring~$\ZZ$ localised  at the prime $\gen{p}$) with $\Rad\gV=p\gV$, the regular parameter is~$p$, and the \crdl is~$\Fp$.
\item Here $\gK=\QQ(t)$, the standard $t$-adic \vala is $|r|_t=t^{-k}$ if~$r=\frac m n\,t^k$ \hbox{with $m,n\in\QQ[t]$}, $m(0)$ and $n(0)\neq 0$ in $\QQ$.
The corresponding \advd is~\hbox{$\gV=(\QQ[t])_{1+t\QQ[t]}$}  (the  ring $\QQ[t]$ localised  at the prime $\gen{t}$) with $\Rad\gV=t\gV$,  the regular parameter is $t$, and the \crdl is~$\QQ$. 
\item Here $\gK=\Fp(t)$, the standard $t$-adic \vala is $|r|_t=t^{-k}$ if~$r=\frac m n\,t^k$ \hbox{with $m,n\in\Fp[t]$}, $m(0)$ and $n(0)\neq 0$ in $\Fp$.
The corresponding \advd is \hbox{$\gV=(\Fp[t])_{1+t\Fp[t]}$}  (the ring  $\Fp[t]$ localised at the prime $\gen{t}$) with $\Rad\gV=t\gV$,  the regular parameter is $t$, and the \crdl is~$\Fp$.

\end{enumerate}

In the second example, the field $\QQ$ may be replaced with an arbitrary \cdi (as in the third example). 

\subsection{The \LHm for ultrametric discrete fields}

\Subsubsection{A  crucial result in \texorpdfstring{\citealt*{MRR}}{Mines, Richman, and Ruitenburg 1988}, and a simpler proof}

\begin{notation} \label{notaKtilde}
Let  $(\gK,\abs\cdot)$ be an ultrametric field with a nontrivial \vala.\footnote{\Cad there exists an $x$ such that $\abs x\neq 0,1$.} We  let 
\begin{itemize}\itemsep=.1em
\item $\gKt$ be the separable closure of~$\gK$ in its completion $\gKh$,
\item $\gVt=\sootq{x\in\gKt}{\abs x \leq 1}$ be its \gui{valuation ring},   and
\item $\fmti=\sootq{x\in\gKt}{\abs x < 1}$.
\end{itemize}
\end{notation}

\citet*{MRR} prove that for an ultrametric discrete field $(\gK,\abs \cdot)$,  $(\gKt,\gVt)$ is an \cvdh with the usual meaning (any \polH has an \zeH).

The \demo of this result is rather complicated because \citet*{MRR} prove   \gnl results concerning the nondiscrete case.
Therefore we consider appropriate to include  here a simpler \demo of the result, which is done in the following two lemmas. 



\begin{lemma} \label{lemComplete}
Let $(\gK,\abs\cdot)$ be an ultrametric discrete field. 
\begin{enumerate}
\item Then  $(\gKh,\abs\cdot)$ is an ultrametric field, $\gVh$ is a \alo with Jacobson radical $\fmh$, and the residual ring $\gVh/\fmh$  is isomorphic to $\gV/\fm$ (it is a \cdi).
\item The \alo  $(\gVh,\fmh)$ is henselian. More \gnlt any \sysN $(\lfn)\in\gVh[\uX]^n$ at $(\an)\in {\gVh}^n$ has a zero $(\xin)$
with $\xi_k\in a_k+\fmh$ ($k\in\lrbn$).
%
%
\end{enumerate}
\end{lemma}
%
\begin{proof} The first item is easy. So $\gKh$ is an Heyting field, a fortiori a \alo. The second item is Corollary~\ref{LHM1}, a consequence of Newton's method explained in \thref{thNewtonQuad}. 
\end{proof}
%

\begin{lemma} \label{lemKtilde} 

\noindent 
The elements of $\gKt$ form a discrete subring of $\gKh$ and $(\gKt,\abs\cdot)$ is an henselian ultrametric discrete field. 
\end{lemma}
%
\begin{proof} The fact that $\gKt$ is a subring is classical.
Let us now consider an \elt $\xi\in\gKh$ annihilating a \spl \pol  $f\in\KX$. 
We let $\Kx=\aqo\KX f$. We have a \Klg morphism
$\varphi\colon\Kx \to \gK[\xi]$ satisfying $\varphi(x)=\xi$. 
As~$\gK[\xi]$ is connected (it is a subring of the \alo $\gKh$), 
it is a \cdi as quotient of a \ste \Klg:  a connected \zedr ring is a \cdi \citep[Fact IV-8.8]{CACM}. So  $\gKt$ is a \cdi. 
\noindent \hum{Il faut en plus justifier le fait que le corps résiduel is lui-aussi discret} 
\end{proof}
%

\Subsubsection{The \iso between two variations on henselisation}

First we recall a variant in the Hensel-Newton style of  Hensel's Lemma  for univariate \pols \cite[Proposition XII-7.6]{Lang}.
This works for all \cvdhs.
\begin{lemma}[Hensel-Newton Lemma for \cvds] \label{lemNewtonHensel} ~\\
Let $\KV$ be a \cvd. We let $v\colon\gK\to\Gamma\cup\so{+\infty}$ be the associated \valn. Let $F(x)=\sum_{k=0}^na_kx^k\in\Vx$ with $a_1\neq 0$ and  $v(a_0)>2v(a_1)$. 
\begin{enumerate}
\item The \pol $F$ has a zero in $\frac{a_0}{a_1}\cdot\VHe\subseteq a_{1}\cdot\fm \VHe$, and it is the unique zero of $F$ in $a_{1}\cdot\fm \VHe$.
\item In particular if $\KV$ is henselian, $F$ has a zero $\xi\in\frac{a_0}{a_1}\cdot\gV\subseteq a_{1}\cdot\fm$, and it is the unique zero of~$F$ in $a_{1}\cdot\fm$.
\end{enumerate}
\end{lemma}
%
\begin{proof}
Let us consider the \pol 
\[
f(x)=\frac{1}{a_1^2}\,F(a_1x)=\frac{a_0}{a_1^2}+X+\som_{k=2}^na_ka_1^{k-2}X^k.
\]
The hypotheses of Hervé's trick (Lemma~\ref{lemTrick1}) are satisfied, so the \pol $f$ has a zero~$\zeta\in\frac{a_0}{a_1^2}\cdot\VHe\subseteq \fm\VHe$, and it is the unique zero of $f$ in $\fm\VHe$. This yields the zero $\xi=a_1\zeta$ of $F$ in the ideal 
$\frac{a_0}{a_1}\cdot\VHe\subseteq a_1\cdot\fm\VHe$, 
and it is the unique zero of $F$ in $a_1\cdot\fm \VHe$.    
\end{proof}

Rereading the \demo of Lemma~\ref{lemTrick1}, we see that the computation shows that  $\zeta$ is in the image of an initial stage of the construction of the henselisation $\Vhe$, obtained by adding the special zero of a special \pol, but $\zeta$ itself is not \ncrt an \zeH of an \cdH. Moreover $v(\zeta)>0$.     

\begin{theorem} [two equivalent versions of the henselisation of an ultrametric discrete field] \label{thMRRKL}~\\
Let $(\gK,\abs\cdot)$ be an ultrametric discrete field.
\begin{enumerate}
\item The henselisation  $(\KHe,\VHe)$ of $\KV$ constructed as in \citealt*{KL00} is isomorphic to~$(\gKt,\gVt)$.
\end{enumerate}
More \prmt, \setcounter{enumi}{1}
\begin{enumerate}
\item there exists a unique $\gK$-\homo  $\KHe\to\gKt$ sending $\VHe$ into~$\gVt$, and this \homo is an  \iso;
\item the natural morphism $\Vhe\to\gVt$ is onto and the morphism $\VHe\to\gVt$ is an \iso.
\end{enumerate}
\end{theorem}
\begin{proof}
Since $(\gKt,\gVt)$ is an \cvdh extension of~$\KV$, we have a unique  $\KV$-morphism $\varphi\colon(\KHe,\VHe)\to(\gKt,\gVt)$. It is injective because~$\KHe$ is a \cdi and~$\gKt$ is not trivial.
We have to show that $\varphi$ is surjective.

It is sufficient to prove Item \textsl{2}. 

We call $v$ the valuation of $\gKh$. 
Let us consider an \elt $\xi\in\gVt$: $\xi\in\gKh$, $v(\xi)\geq 0$, and $f(\xi)=0$ for a \pol $f\in\VX$ which is \spl in $\KX$. In particular $f'(\xi)\neq 0$, i.e. $v(f'(\xi))<+\infty$.

Since $\xi\in\gKh$ we know arbitrarily precise approximations of $\xi$ in $\gK$, \cad written as 
\hbox{$a=\xi+\zeta\in\gK$}  with $v(\zeta)$ arbitrarily large. Since $v(\xi)\geq 0$, we have $v(a)\geq 0$.

For such an $a$ we consider the \pol 
$F_a(X)=f(-X+a) \in \VX$.
The \coe $c_0=F_a(0)=f(a)$ is an arbitrarily precise approximation of $f(\xi)=0$, \cad $v(c_0)$ is arbitrarily large. 
The \coe $c_1=F'_a(0)=f'(a)$ is an arbitrarily precise approximation of~$f'(\xi)$, so for $v(\zeta)$ sufficiently large $v(c_1)=v(f'(a))=v(f'(\xi))<+\infty$. For $v(\zeta)$ sufficiently large we get $v(c_0)>2v(c_1)$. Lemma~\ref{lemNewtonHensel} says that $F_a$ has a zero  $\alpha\in\frac{c_0}{c_1}\cdot\VHe\subseteq c_1\cdot\fm \VHe$. 

The image of $\alpha$ in $\gKt$ is equal to $\zeta$, for these are two \zeHs in $\gKt$ for the same \cdH $(F_a,0)$.
So the image  of $a-\alpha \in \Vhe$ is $a-\zeta=\xi$.
Moreover, as explained just after Lemma~\ref{lemNewtonHensel}, $\xi$ is
also in the range of the natural morphism  $\Vhe\to\gVt$.    
\end{proof}

So we have got that any \elt of $\gKt$ is \prmt the image of an   \elt $\gamma$
belonging to a field $\gK[\xi]\subseteq \KHe$ where 
$\xi$ is the special zero of a special \pol. 
Since $\KHe\simeq\gKt$, this implies that any \elt of $\KHe$
can be obtained as an \elt at the first stage of some construction of~$\KHe$.
This result was not clear a priori in \citealt{KL00} (but there the \cvd is arbitrary).

\Subsubsection{\LHm}

Now we get the desired result.

\begin{theorem}[\LHm  for an ultrametric discrete field] \label{LHM2} ~\\
Let $(\gK,\abs\cdot)$ be an ultrametric discrete field and 
$(f_1,\dots,f_n)$ a \sysN at $(\uze)$ over $\KV$. This system admits a unique zero with  \coos in $\fmti$.  It admits also a unique zero with  \coos in $\fm\VHe$.
\end{theorem}
%
\begin{proof}
We assume that the \sysN is étale (Lemma~\ref{note0}). 
Newton's method constructs a zero $(\ual)$ with \coos in $\fmh\subseteq \gKh$ (Corollary~\ref{LHM1}).
We let $\gD=\Kxn=\aqo\KXn{\lfn}$ be the quotient \Klg of the \syp.
It is \stf, \ste (see \Thos~\ref{thNewtonCorpsDiscret} and~\ref{thstrucste}).
We have the natural morphism of \Klgs $\varphi\colon\gD\to\gK[\ual]\subseteq \gKh$, where $(\ual)$ is the \zeH of the \syp  ($\varphi(x_i)=\alpha_i$ for each~$i$).  
So the~$\alpha_i$'s are \spl over $\gK$.  
Thus, for each $i$,  $\alpha_i\in\gKt$ and $v(\alpha_i)>0$: the \coos of $(\ual)$  are in $\fmti$. 
\\
On the other hand, as we have proved (\thref{thMRRKL}), $\gVt$ is canonically isomorphic to~$\VHe$.  The elements  in $\VHe$ that correspond to the $\alpha_i$'s give the \zeH of the \syp with \coos in~$\fm\VHe$.
\end{proof}

Note that according to Description~\ref{descriEtaleCdi},
since the \Klg $\gK[\ual]$ is connected nontrivial, $\gK[\ual]$ is a \cdi isomorphic to a quotient of $\gD$, but it seems that there is no general algorithm for determining the dimension of  $\gK[\ual]$ as \Kev. 

\medskip 
\centerline{--------------------------}

\medskip \noindent \textsl{Final remark.}
Note that many classical texts using variants of \LHm as \textsl{e.g.}  \cite{Fis1997} and \cite{Sma1998} give the solution under the form of zeroes with coordinates in some completion of the local ring, and not in the henselisation of the local ring. 
\\
Papers \cite{Kuh2011} and \cite{PCR2000} give a non \algq solution by using the notion of spherically complete field. Their \demos are very difficult to interpret \cot.

\addcontentsline{toc}{section}{References}
\rdb
\small
\bibliographystyle{plainnat}

\normalsize
\endgroup
\stopcontents[english]

\clearpage
\newpage
\thispagestyle{empty}

\clearpage
\newpage
~

\setcounter{page}{1}\renewcommand\thepage{F\arabic{page}}\renewcommand\theHsection{F\arabic{section}}

\begingroup

\startcontents[french]

\setcounter{page}{0}

\newcommand{\di}{\,\vert\,}
\newcommand{\ndi}{\nmid}

\newcommand \nzr {\neq 0}
\newcommand{\Vr}{\mathrm{Vr}}
\newcommand{\Rn}{\mathrm{Rn}}
\newcommand{\Nvr}{\mathrm{Nvr}}
\newcommand{\Nrn}{\mathrm{Nrn}}
\newcommand{\U}{\mathrm{Un}}

\newcommand{\Gati}{\Gat\cup\so\infty }

\newcommand \Rzero {R_{=0}}
\newcommand \Rnz {R_{\nzr}}
\newcommand \Rvr {R_{\Vr}}
\newcommand \Ru {R_{\U}}
\newcommand \Rrn {R_{\Rn}}
\newcommand \Rnrn {R_{\Nrn}}
\newcommand \Rnvr {R_{\Nvr}}

\newcommand \zg {{\ZZ[G]}}
\newcommand \Izero {{\cI}_{=0}}
\newcommand \Mnz {{\cM}_{\nzr}}
\newcommand \Mpos {{\cM}_{\Pos}}
\newcommand \Cnng {{\cC}_{\Nng}}
\newcommand \Irn {{\cI}_{\Rn}}
\newcommand \Vvr {{\cV}_{\Vr}}
\newcommand \Mu {{\cM}_{\U}}
\newcommand \abg {{\rm Ab}(G)}
\newcommand \Hzero {{\cH}_{=0}}


\newcommand\Exists{\boldsymbol{\stixexists}}
\newcommand\Forall{\boldsymbol{\stixforall}}
\newcommand\VDash{\boldsymbol{\stixvdash}}
\newcommand\Land{\boldsymbol{\stixwedge}}
\newcommand\Lor{\boldsymbol{\stixvee}}
\newcommand\Top{\boldsymbol{\stixtop}}
\newcommand\Bot{\boldsymbol{\stixbot}}
\newcommand\bigLand{\boldsymbol{\stixbigwedge}}
\newcommand\bigLor{\boldsymbol{\stixbigvee}}

\newcommand\dotminus{\stixdotminus}

\newcommand \vdg{\VDash}
\newcommand \Vd {\,\vdg}
\newcommand \vd {\,\,\vdg}
\newcommand \vii{\Land}
\newcommand \vuu{\Lor}
\newcommand \vdi[1] {\mathrel{\,\,\vdg_{#1}}}
\newcommand \Vii{\bigLand}
\newcommand \Vuu{\bigLor}


\newtheorem{ftheoreme}{Théorème}[section]
\newtheorem{ftheorem}{Théorème}[subsection]
\newtheorem{fthdef}[ftheorem]{Théorème et définition}
\newtheorem{fvastf}[ftheorem]{\vst formel}
\newtheorem{fvast}[ftheorem]{\vst}
\newtheorem{pstf}[ftheorem]{Positivstellensatz formel}
\newtheorem{fpst}[ftheorem]{Positivstellensatz}
\newtheorem{flemma}[ftheorem]{Lemme}
\newtheorem{fcorollary}[ftheorem]{Corolaire}
\newtheorem{fconjecture}[ftheorem]{Conjecture}
\newtheorem{fproposition}[ftheorem]{Proposition}
\newtheorem{fprpta}[ftheorem]{Propriétés attendues}
\newtheorem{fpropdef}[ftheorem]{Proposition et définition}
\newtheorem{ffact}[ftheorem]{Fait}
\newtheorem{fconvention}[ftheorem]{Convention}
\newtheorem{fplcc}[ftheorem]{Principe local-global concret}
\newtheorem{fprecision}[ftheorem]{Précision}

\newtheorem{ftheoremc}[ftheorem]{Th\'{e}or\`{e}me\etoz}
\newtheorem{flemmac}[ftheorem]{Lemme\etoz}
\newtheorem{fcorollaryc}[ftheorem]{Corolaire\etoz}
\newtheorem{fproprietec}[ftheorem]{Propri\'{e}t\'{e}\etoz}
\newtheorem{fpropositionc}[ftheorem]{Proposition\etoz}
\newtheorem{ffactc}[ftheorem]{Fait\etoz}
\newtheorem{fpropdefc}[ftheorem]{Proposition et définition\etoz}

\theoremstyle{definition}
\newtheorem{fdefinition}[ftheorem]{Définition}
\newtheorem{fdfni}[ftheorem]{Définition informelle}
\newtheorem{fdefinitions}[ftheorem]{Définitions}
\newtheorem{fexample}[ftheorem]{Exemple}
\newtheorem{fexamples}[ftheorem]{Exemples}
\newtheorem{fnotation}[ftheorem]{Notation}
\newtheorem{fter}[ftheorem]{Terminologie}
\newtheorem{fproblem}[ftheorem]{Problème}
\newtheorem{fquestion}[ftheorem]{Question}
\newtheorem{fquestions}[ftheorem]{Questions}
\newtheorem{fcontext}[ftheorem]{Contexte}
\newtheorem{fdescri}[ftheorem]{Description}

\newtheorem{fdefinitionc}[ftheorem]{Définition\etoz}
\newtheorem{fdefinota}[ftheorem]{Définition et notation}

\theoremstyle{remark}
\newtheorem{fnote}[ftheorem]{Note}
\newtheorem{fnotes}[ftheorem]{Notes}
\newtheorem{fremark}[ftheorem]{Remarque}
\newtheorem{fremarks}[ftheorem]{Remarques}
\newtheorem{fcomment}[ftheorem]{Commentaire}

\newcommand{\vou}{\MA{\tsbf{ ou }}}
\newcommand{\Vou}{\MA{\tsbf{OU}}}
\newcommand \EXists[1] {\tsbf{Introduire }{#1}\tsbf{ tel que }\,}
\newcommand \vet {\tsbf{,}\;}
\newcommand \Atcl {\mathrm{Atcl}}



\newcounter{MF}
\newcommand\stMF{\stepcounter{MF}}

\newcommand{\lec}{\stMF\ifodd\value{MF}lecteur\xspace\else 
lectrice\xspace\fi}

\newcommand{\lecs}{\stMF\ifodd\value{MF}lecteurs\xspace\else 
lectrices\xspace\fi}

\newcommand{\alec}{\stMF\ifodd\value{MF}au lecteur\xspace\else%
à la lectrice\xspace\fi}

\newcommand{\dlec}{\stMF\ifodd\value{MF}du lecteur\xspace\else%
de la lectrice\xspace\fi}

\newcommand{\llec}{\stMF\ifodd\value{MF}le lecteur\xspace\else la lectrice\xspace\fi}

\newcommand{\Llec}{\stMF\ifodd\value{MF}Le lecteur\xspace\else La lectrice\xspace\fi}

\newcommand{\lui}{\ifodd\value{MF}lui\xspace\else
elle\xspace\fi}

\newcommand{\celui}{\ifodd\value{MF}celui\xspace\else
celle\xspace\fi}

\newcommand{\ceux}{\ifodd\value{MF}ceux\xspace\else
celles\xspace\fi}

\newcommand{\er}{\ifodd\value{MF}er\xspace\else
ère\xspace\fi}

\newcommand{\eux}{\ifodd\value{MF}eux\xspace\else
elles\xspace\fi}

\newcommand{\eUx}{\ifodd\value{MF}eux\xspace\else
euse\xspace\fi}

\newcommand{\eUX}{\ifodd\value{MF}eux\xspace\else
euses\xspace\fi}

\newcommand{\leux}{\ifodd\value{MF}leux\xspace\else
leuse\xspace\fi}

\newcommand{\il}{\ifodd\value{MF}il\xspace\else
elle\xspace\fi}

\newcommand{\ien}{\ifodd\value{MF}ien\xspace\else
ienne\xspace\fi}

\newcommand{\e}{\ifodd\value{MF}\xspace \else e\xspace\fi}

\newcommand{\n}{\ifodd\value{MF}n\xspace\else nne\xspace\fi}

\makeatletter
\newcommand{\la}{\@ifstar{\ifodd\value{MF}le\else
la\fi}{\stMF\ifodd\value{MF}le\else la\fi}}
\makeatother

\newcommand \rem{\rdb
\noi{\it Remarque. }}

\newcommand \REM[1]{\rdb
\noi{\it Remarque#1. }}

\newcommand \rems{\rdb
\noi{\it Remarques. }}

\newcommand \exl{\rdb
\noi{\bf Exemple. }}

\newcommand \EXL[1]{\rdb
\noi{\bf Exemple: #1. }}

\newcommand \exls{\rdb
\noi{\bf Exemples. }}

\newcommand \thref[1] {théorème~\ref{#1}}
\newcommand \paref[1] {page~\pageref{#1}}
\newcommand \pstfref[1] {Posi\-tiv\-stel\-lensatz formel~\ref{#1}}
\newcommand \pstref[1] {Posi\-tiv\-stel\-lensatz~\ref{#1}}

\newcommand\oge{\leavevmode\raise.3ex\hbox{$\scriptscriptstyle\langle\!\langle\,$}}
\newcommand\feg{\leavevmode\raise.3ex\hbox{$\scriptscriptstyle\,\rangle\!\rangle$}}

\newcommand\gui[1]{\oge{#1}\feg}

\newcommand \facile{\begin{proof}
La démonstration est laissée \alec.
\end{proof}
}

\def \num {\no} 

\newcommand\comm{\rdb
\noi{\sl Commentaire. }}

\newcommand\COM[1]{\rdb
\noi{\sl Commentaire #1. }}

\newcommand\comms{\rdb
\noi{\sl Commentaires. }}

\newcommand\Pb{\rdb
\noi{\bf Problème. }}

\newcommand\eoq{\hbox{}\nobreak
\vrule width 1.4mm height 1.4mm depth 0mm}

\newcommand \Cad {C'est-à-dire\xspace}
\newcommand \recu {récur\-rence\xspace}
\newcommand \hdr {hypo\-thèse de \recu}
\newcommand \cad {c'est-à-dire\xspace}
\newcommand \cade {c'est-à-dire en\-co\-re\xspace}
\newcommand \ssi {si, et seu\-lement si, }
\newcommand \ssiz {si, et seu\-lement si,~}
\newcommand \cnes {con\-di\-tion néces\-sai\-re et suf\-fi\-san\-te\xspace}
\newcommand \spdg {sans per\-te de géné\-ra\-lité\xspace}
\newcommand \Spdg {Sans per\-te de géné\-ra\-lité\xspace}

\newcommand \Propeq {Les pro\-pri\-é\-tés sui\-van\-tes sont 
équi\-va\-len\-tes.}
\newcommand \propeq {les pro\-pri\-é\-tés sui\-van\-tes sont 
équi\-va\-len\-tes.}

\newcommand \Kev {$\gK$-\evc}
\newcommand \Kbev {$\gKb$-\evc}
\newcommand \Kevs {$\gK$-\evcs}

\newcommand \Lev {$\gL$-\evc}
\newcommand \Levs {$\gL$-\evcs}

\newcommand \Qev {$\QQ$-\evc}
\newcommand \Qevs {$\QQ$-\evcs}

\newcommand \kev {$\gk$-\evc}
\newcommand \kevs {$\gk$-\evcs}

\newcommand \lev {$\gl$-\evc}
\newcommand \levs {$\gl$-\evcs}

\newcommand \Alg {$\gA$-\alg}
\newcommand \Algs {$\gA$-\algs}

\newcommand \Blg {$\gB$-\alg}
\newcommand \Blgs {$\gB$-\algs}

\newcommand \Clg {$\gC$-\alg}
\newcommand \Clgs {$\gC$-\algs}

\newcommand \klg {$\gk$-\alg}
\newcommand \klgs {$\gk$-\algs}

\newcommand \llg {$\gl$-\alg}
\newcommand \llgs {$\gl$-\algs}

\newcommand \Klg {$\gK$-\alg}
\newcommand \Klgs {$\gK$-\algs}

\newcommand \Llg {$\gL$-\alg}
\newcommand \Llgs {$\gL$-\algs}

\newcommand \QQlg {$\QQ$-\alg}
\newcommand \QQlgs {$\QQ$-\algs}

\newcommand \Rlg {$\gR$-\alg}
\newcommand \Rlgs {$\gR$-\algs}

\newcommand \RRlg {$\RR$-\alg}
\newcommand \RRlgs {$\RR$-\algs}

\newcommand \Vlg {$\gV$-\alg}
\newcommand \Vlgs {$\gV$-\algs}

\newcommand \ZZlg {$\ZZ$-\alg}
\newcommand \ZZlgs {$\ZZ$-\algs}

\newcommand \Amo {$\gA$-module\xspace}
\newcommand \Amos {$\gA$-modules\xspace}

\newcommand \Bmo {$\gB$-module\xspace}
\newcommand \Bmos {$\gB$-modules\xspace}

\newcommand \Cmo {$\gC$-module\xspace}
\newcommand \Cmos {$\gC$-modules\xspace}

\newcommand \kmo {$\gk$-module\xspace}
\newcommand \kmos {$\gk$-modules\xspace}

\newcommand \Kmo {$\gK$-module\xspace}
\newcommand \Kmos {$\gK$-modules\xspace}

\newcommand \Lmo {$\gL$-module\xspace}
\newcommand \Lmos {$\gL$-modules\xspace}

\newcommand \Rmo {$\gR$-module\xspace}
\newcommand \Rmos {$\gR$-modules\xspace}

\newcommand \Vmo {$\gV$-module\xspace}
\newcommand \Vmos {$\gV$-modules\xspace}

\newcommand \Ali {appli\-ca\-tion $\gA$-\lin}
\newcommand \Alis {appli\-ca\-tions $\gA$-\lins}

\newcommand \Kli {appli\-ca\-tion $\gK$-\lin}
\newcommand \Klis {appli\-ca\-tions $\gK$-\lins}

\newcommand \Bli {appli\-ca\-tion $\gB$-\lin}
\newcommand \Blis {appli\-ca\-tions $\gB$-\lins}

\newcommand \Cli {appli\-ca\-tion $\gC$-\lin}
\newcommand \Clis {appli\-ca\-tions $\gC$-\lins}

\newcommand \ac{algébriquement clos\xspace}  

\newcommand \acl {an\-neau \icl}
\newcommand \acls {an\-neaux \icl}

\newcommand \adp {an\-neau de Pr\"u\-fer\xspace}
\newcommand \adps {an\-neaux de Pr\"u\-fer\xspace}

\newcommand \adpc {\adp \coh}
\newcommand \adpcs {\adps \cohs}

\newcommand \adu {\alg de décom\-po\-sition univer\-selle\xspace}
\newcommand \adus {\algs de décom\-po\-sition univer\-selle\xspace}

\newcommand \adv {anneau de valuation\xspace}
\newcommand \advs {anneaux de valuation\xspace}

\newcommand \advd {anneau de valuation discrète\xspace}
\newcommand \advds {anneaux de valuation discrète\xspace}

\newcommand \advl {anneau \dvla} 
\newcommand \advls {anneaux \dvlas} 

\newcommand \Afr {Anneau \frl}
\newcommand \Afrs {Anneaux \frls}
\newcommand \afr {anneau \frl}
\newcommand \afrs {anneaux \frls}

\newcommand \aFr {\hyperref[theorieAfr]{anneau \frl}}
\newcommand \aFrs {\hyperref[theorieAfr]{anneau \frls}}

\newcommand \afrr {\afr réduit\xspace}
\newcommand \afrrs {\afrs réduits\xspace}
\newcommand \Afrrs {\Afrs réduits\xspace}

\newcommand \afrvr {\afr avec \ravs}
\newcommand \aFrvr {\hyperref[theorieAfrrv]{\afrvr}}
\newcommand \afrvrs {\afrs avec \ravs}

\newcommand \aftr {anneau réticulé \ftm réel\xspace}
\newcommand \aftrs {anneaux réticulés \ftm réels\xspace}

\newcommand \aG {\alg galoisienne\xspace}
\newcommand \aGs {\algs galoisiennes\xspace}

\newcommand \agB {\alg de Boole\xspace}
\newcommand \agBs {\algs de Boole\xspace}

\newcommand \agH {\alg de Heyting\xspace}
\newcommand \agHs {\algs de Heyting\xspace}

\newcommand \agq{algé\-bri\-que\xspace}
\newcommand \agqs{algé\-bri\-ques\xspace}

\newcommand \agqt{algé\-bri\-que\-ment\xspace}

\newcommand \aKr {anneau de Krull\xspace}
\newcommand \aKrs {anneaux de Krull\xspace}

\newcommand \ale {\alg étale\xspace}
\newcommand \ales {\algs étales\xspace}

\newcommand \alg {algèbre\xspace}
\newcommand \algs {algèbres\xspace}

\newcommand \algo{algo\-rithme\xspace}
\newcommand \algos{algo\-rithmes\xspace}

\newcommand \algq{al\-go\-rith\-mi\-que\xspace}
\newcommand \algqs{al\-go\-rith\-mi\-ques\xspace}

\newcommand \ali {appli\-ca\-tion \lin}
\newcommand \alis {appli\-ca\-tions \lins}

\newcommand \alo {anneau local\xspace}
\newcommand \alos {anneaux locaux\xspace}

\newcommand \algb {an\-neau \lgb}
\newcommand \algbs {an\-neaux \lgbs}

\newcommand \alrd {\alo \dcd}
\newcommand \alrds {\alos \dcds}

\newcommand \alrdh{\alrd hensélien\xspace}
\newcommand \alrdhs{\alrds henséliens\xspace}

\newcommand \anar {anneau \ari}
\newcommand \anars {anneaux \aris}

\newcommand \anor {an\-neau nor\-mal\xspace}
\newcommand \anors {an\-neaux nor\-maux\xspace}

\newcommand \apf {\alg \pf}
\newcommand \apfs {\algs \pf}

\newcommand \apG {\alg pré\-galoisienne\xspace}
\newcommand \apGs {\algs pré\-galoisiennes\xspace}

\newcommand \arch {archimédien\xspace}
\newcommand \arche {archimédienne\xspace}
\newcommand \archs {archimédiens\xspace}
\newcommand \arches {archimédiennes\xspace}

\newcommand \arc {anneau réel clos\xspace}
\newcommand \aRc {\hyperref[theorieArc]{\arc}}
\newcommand \arcs {anneaux réels clos\xspace}

\newcommand \ari{arith\-mé\-tique\xspace}  
\newcommand \aris{arith\-mé\-tiques\xspace}  

\newcommand \arv {\adv}
\newcommand \arvs {\advs}

\newcommand \Asr {Anneau \str}
\newcommand \Asrs {Anneaux \strs}
\newcommand \asr {anneau \str}
\newcommand \asrs {anneaux \strs}

\newcommand \asrvr {\asr avec \ravs}
\newcommand \asrvrs {\asrs avec \ravs}

\newcommand \atf {\alg \tf}
\newcommand \atfs {\algs \tf}

\newcommand \atfr {anneau total de fractions\xspace}
\newcommand \atfrs {anneaux totals de fractions\xspace}

\newcommand \auto {auto\-mor\-phisme\xspace}
\newcommand \autos {auto\-mor\-phismes\xspace}


\newcommand \bdg {base de Gr\"obner\xspace}
\newcommand \bdgs {bases de Gr\"obner\xspace}

\newcommand \bdp {base de \dcn partielle\xspace}
\newcommand \bdps {bases de \dcn partielle\xspace}

\newcommand \bdf {base de \fap}

\newcommand \Bif {Borne infé\-rieure\xspace} %
\newcommand \bif {borne infé\-rieure\xspace} %
\newcommand \bifs {bornes infé\-rieures\xspace} %

\newcommand \bsp {borne supé\-rieure\xspace} %
\newcommand \bsps {borne supé\-rieures\xspace} %


\newcommand \cac{corps \ac}  

\newcommand \calf{calcul formel\xspace}  

\newcommand \cara{carac\-té\-ris\-tique\xspace}  
\newcommand \caras{carac\-té\-ris\-tiques\xspace}  

\newcommand \carn{carac\-té\-ri\-sation\xspace}  
\newcommand \carns{carac\-té\-ri\-sations\xspace}  

\newcommand \carar{carac\-té\-riser\xspace}

\newcommand \carf{de carac\-tère fini\xspace}  

\newcommand \cdac{\cdi \ac}  
\newcommand \cdacs{\cdis \ac}  
\newcommand \cdi{corps discret\xspace}
\newcommand \cdis{corps discrets\xspace}

\newcommand \cdf{corps de fractions\xspace}

\newcommand \cdH{code de Hensel\xspace}
\newcommand \cdHs{codes de Hensel\xspace}

\newcommand \cdr{corps de racines\xspace}
  
\newcommand \cdv{changement de variables\xspace}  
\newcommand \cdvs{changements de variables\xspace}

\newcommand \cla {clôture \agq}
\newcommand \clas {clôtures \agqs}

\newcommand \cli {clôture intégrale\xspace}
\newcommand \clis {clôtures intégrales\xspace}

\newcommand \clr {clôture réelle\xspace}
\newcommand \clrs {clôtures réelles\xspace}

\newcommand \clsep {clôture \spl}
\newcommand \clseps {clôtures \spls}

\newcommand \codi {corps ordonné discret\xspace}
\newcommand \codis {corps ordonnés discrets\xspace}

\newcommand \coe {coef\-fi\-cient\xspace}
\newcommand \coes {coef\-fi\-cients\xspace}

\newcommand \coh {co\-hé\-rent\xspace}
\newcommand \cohs {co\-hé\-rents\xspace}

\newcommand \cohc {co\-hé\-rence\xspace}

\newcommand \colH {\colo hensélien\xspace}
\newcommand \colHs {\colos henséliens\xspace}

\newcommand \coli {combinaison \lin}
\newcommand \colis {combinaisons \lins}

\newcommand \com {co\-ma\-xi\-maux\xspace}
\newcommand \come {co\-ma\-xi\-males\xspace}

\newcommand \colo {\cou local\xspace}
\newcommand \colos {\cous locaux\xspace}

\newcommand \coo {coor\-donnée\xspace}
\newcommand \coos {coor\-données\xspace}

\newcommand \cop {complé\-men\-taire\xspace}
\newcommand \cops {complé\-men\-taires\xspace}

\newcommand \cosv {conser\-vative\xspace}
\newcommand \cosvs {conser\-vatives\xspace}

\newcommand \cOsv {\hyperref[defithconserv]{conser\-vative\xspace}}
\newcommand \cOsvs {\hyperref[defithconserv]{conser\-vatives\xspace}}

\newcommand \cou {couple\xspace}
\newcommand \cous {couples\xspace}

\newcommand \covr {corps ordonné avec \ravs}
\newcommand \covrs {corps ordonnés avec \ravs}

\newcommand \cpb {compa\-tible\xspace} 
\newcommand \cpbs {compa\-tibles\xspace} 

\newcommand \cpbt {compa\-tibi\-lité\xspace} 
\newcommand \cpbtz {compa\-tibi\-lité}

\newcommand \crdl {corps résiduel\xspace}
\newcommand \crdls {corps résiduels\xspace}

\newcommand \crc {corps réel clos\xspace}
\newcommand \crcs {corps réels clos\xspace}

\newcommand \crcd {corps réel clos discret\xspace}
\newcommand \crcds {corps réels clos discrets\xspace}

\newcommand \crl {corolaire\xspace}
\newcommand \crls {corolaires\xspace}

\newcommand \cval{corps valué\xspace}
\newcommand \cvals{corps valués\xspace}

\newcommand \cvar{corps valorisé\xspace}
\newcommand \cvars{corps valorisés\xspace}

\newcommand \cvd{\cval discret\xspace}
\newcommand \cvds{\cvals discrets\xspace}
\newcommand \Cvd{Corps valué discret\xspace}
\newcommand \Cvds{Corps valués discrets\xspace}

\newcommand \cvdac{\cvd\ac}
\newcommand \cvdacs{\cvds\ac}

\newcommand \cvdh{\cvd hensélien\xspace}
\newcommand \cvdhs{\cvds henséliens\xspace}

\newcommand \cvdsc{\cvd \splc}
\newcommand \cvdscs{\cvds \splc}

\newcommand \cvdu{\cvar discret ultramétrique\xspace}
\newcommand \cvdus{\cvars discrets ultramétriques\xspace}

\newcommand \cud{corps \ultm discret\xspace}
\newcommand \cuds{corps \ultms discrets\xspace}

\newcommand \cvu{\cvar \ultm}
\newcommand \cvus{\cvars \ultms}


\newcommand \dcd {résiduellement dis\-cret\xspace}
\newcommand \dcds {résiduellement dis\-crets\xspace}

\newcommand \dcn {décomposition\xspace}
\newcommand \dcns {décompositions\xspace}

\newcommand \dcnb {\dcn bornée\xspace}

\newcommand \dcnc {\dcn complète\xspace}

\newcommand \dcnp {\dcn partielle\xspace}

\newcommand \dcp {décom\-posa\-ble\xspace}
\newcommand \dcps {décom\-posa\-bles\xspace}

\newcommand \ddk {dimension de~Krull\xspace}
\newcommand \ddi {de dimension infé\-rieure ou égale à~}

\newcommand \dimm {description immédiate\xspace}
\newcommand \dimms {descriptions immédiates\xspace}

\newcommand \ddp {domaine de Pr\"u\-fer\xspace}
\newcommand \ddps {domaines de Pr\"u\-fer\xspace}

\newcommand \Demo{Démonstration\xspace}     
\newcommand \Demos{Démonstrations\xspace}     

\newcommand \demo{démon\-stra\-tion\xspace}     
\newcommand \demos{démon\-stra\-tions\xspace}     

\newcommand \dems{démons\-tra\-tions\xspace}

\newcommand \deno{déno\-mi\-nateur\xspace}     
\newcommand \denos{déno\-mi\-nateurs\xspace}   

\newcommand \deter {déter\-mi\-nant\xspace}  
\newcommand \deters {déter\-mi\-nants\xspace}  
  
\newcommand \Dfn{Défi\-nition\xspace}  
\newcommand \Dfns{Défi\-nitions\xspace}  
\newcommand \dfn{défi\-nition\xspace}  
\newcommand \dfns{défi\-nitions\xspace}  

\newcommand \dftr {droite réticulée \ftm réelle\xspace}
\newcommand \dftrs {droites réticulées \ftm réelles\xspace}
  
\newcommand \dil{diffé\-rentiel\xspace}  
\newcommand \dils{diffé\-rentiels\xspace}  
\newcommand \dile{diffé\-ren\-tielle\xspace}  
\newcommand \diles{diffé\-ren\-tielles\xspace}  

\newcommand \dip{diviseur principal\xspace}
\newcommand \dips{diviseurs principaux\xspace}

\newcommand \discri{discri\-minant\xspace}  
\newcommand \discris{discri\-minants\xspace}  

\newcommand \divle {dimension divisorielle\xspace}

\newcommand \dit{distri\-bu\-ti\-vité\xspace}

\newcommand \dlg{d'élar\-gis\-sement\xspace}  

\newcommand \dok {domaine de Dedekind\xspace}
\newcommand \doks {domaines de Dedekind\xspace}

\newcommand \dvla {à diviseurs\xspace}
\newcommand \dvlas {à diviseurs\xspace}

\newcommand \dvld {\dvlt décom\-posé\xspace} %
\newcommand \dvlds {\dvlt décom\-posés\xspace} %

\newcommand \dvlg {diviso\-riel\xspace} 
\newcommand \dvlgs {diviso\-riels\xspace} 

\newcommand \dvli {\dvlt inver\-sible\xspace} 
\newcommand \dvlis {\dvlt inver\-sibles\xspace} 

\newcommand \dvlt {diviso\-riel\-lement\xspace} %

\newcommand \dvz {di\-viseur de zéro\xspace}
\newcommand \dvzs {di\-viseurs de zéro\xspace}

\newcommand \dve {divi\-si\-bi\-lité\xspace}

\newcommand \dvee {à \dve explicite\xspace}

\newcommand \dvr {diviseur\xspace}
\newcommand \dvrs {diviseurs\xspace}


\newcommand \Eds {Exten\-sion des sca\-laires\xspace}
\newcommand \edss {exten\-sions des sca\-laires\xspace}
\newcommand \eds {exten\-sion des sca\-laires\xspace}

\newcommand \eco {\elts \com}

\newcommand \egmt {éga\-lement\xspace}

\newcommand \egt {éga\-li\-té\xspace}
\newcommand \egts {éga\-li\-tés\xspace}

\newcommand \eimm {extension immédiate\xspace}
\newcommand \eimms {extensions immédiates\xspace}

\newcommand \eli{élimi\-nation\xspace}  

\newcommand \elr{élé\-men\-taire\xspace}  
\newcommand \elrs{élé\-men\-taires\xspace}  

\newcommand \elrt{élé\-men\-tai\-rement\xspace}  

\newcommand \elt{élé\-ment\xspace}  
\newcommand \elts{élé\-ments\xspace}  

\def \endo {en\-do\-mor\-phisme\xspace}
\def \endos {en\-do\-mor\-phismes\xspace}

\newcommand \entrel {rela\-tion impli\-ca\-tive\xspace}
\newcommand \entrels {rela\-tions impli\-ca\-tives\xspace}

\newcommand\evc{es\-pa\-ce vec\-to\-riel\xspace} 
\newcommand\evcs{es\-pa\-ces vec\-to\-riels\xspace} 

\newcommand \eqn {équation\xspace}  
\newcommand \eqns {équations\xspace}  

\newcommand \eqv {équi\-valent\xspace}  
\newcommand \eqve {équi\-va\-lente\xspace}  
\newcommand \eqvs {équi\-valents\xspace}  
\newcommand \eqves {équi\-val\-entes\xspace}  

\newcommand \eqvc {équi\-va\-lence\xspace}  
\newcommand \eqvcs {équi\-va\-lences\xspace}  

\newcommand \esid {essen\-tiel\-lement iden\-tique\xspace}  
\newcommand \esids {essen\-tiel\-lement iden\-tiques\xspace}  

\newcommand \Esid {\hyperref[defitdyesidentiques]{\esid}}  
\newcommand \Esids {\hyperref[defitdyesidentiques]{\esids}}  

\newcommand \eseq {essen\-tiel\-lement \eqve}  
\newcommand \eseqs {essen\-tiel\-lement \eqves}  

\newcommand \Eseq {\hyperref[defitheseq]{\eseq}}  
\newcommand \Eseqs {\hyperref[defitheseq]{\eseqs}}

\newcommand \fab {\fcn bornée\xspace}
\newcommand \fabs {\fcns bornées\xspace}

\newcommand \fat {\fcn totale\xspace}
\newcommand \fats {\fcn totales\xspace}

\newcommand \fap {\fcn partielle\xspace}
\newcommand \faps {\fcns partielles\xspace}

\newcommand \fip {filtre pre\-mier\xspace}
\newcommand \fips {filtres pre\-miers\xspace}

\newcommand \fipma {\fip maximal\xspace}
\newcommand \fipmas {\fips maximaux\xspace}

\newcommand \fcn {factorisation\xspace}
\newcommand \fcns {factorisations\xspace}

\newcommand \fdi {for\-te\-ment dis\-cret\xspace}
\newcommand \fdis {for\-te\-ment dis\-crets\xspace}

\newcommand \fsa {fermé \sagq}
\newcommand \fsas {fermés \sagqs}

\newcommand \fsagc {fonction \sagc}
\newcommand \fsagcs {fonctions \sagcs}

\newcommand \fmt {formellement\xspace}

\newcommand \fit {fidèlement\xspace}
\newcommand \fpt {\fit plat\xspace}
\newcommand \fpte {\fit plate\xspace}
\newcommand \fpts {\fit plats\xspace}
\newcommand \fptes {\fit plates\xspace}

\newcommand \frl {fortement réticulé\xspace}
\newcommand \frle {fortement réticulée\xspace}
\newcommand \frls {fortement réticulés\xspace}

\newcommand \ftm {fortement\xspace}

\newcommand\gmt{géométrie\xspace}  
\newcommand\gmts{géométries\xspace}  

\newcommand\gaq{\gmt \agq}  

\newcommand\gmq{géomé\-trique\xspace}  
\newcommand\gmqs{géomé\-triques\xspace}  

\newcommand\gmqt{géomé\-tri\-quement\xspace}  

\newcommand\gne{géné\-ra\-lisé\xspace}  
\newcommand\gnee{géné\-ra\-lisée\xspace}  
\newcommand\gnes{géné\-ra\-lisés\xspace}  
\newcommand\gnees{géné\-ra\-lisées\xspace}  

\newcommand\gnl{géné\-ral\xspace}  
\newcommand\gnle{géné\-rale\xspace}  
\newcommand\gnls{géné\-raux\xspace}  
\newcommand\gnles{géné\-rales\xspace}  

\newcommand\gnlt{géné\-ra\-lement\xspace}  

\newcommand\gnn{géné\-ra\-li\-sa\-tion\xspace}  
\newcommand\gnns{géné\-ra\-li\-sa\-tions\xspace}  

\newcommand\gnq {géné\-rique\xspace}  
\newcommand\gnqs {géné\-riques\xspace}  

\newcommand\gnr{géné\-ra\-liser\xspace}  

\newcommand \gns{géné\-ra\-lise\xspace}

\newcommand \gnt{géné\-ra\-lité\xspace}
\newcommand \gnts{géné\-ra\-lités\xspace}

\newcommand \grl{groupe \rtl}
\newcommand \grls{groupes \rtls}

\newcommand \gRl {\hyperref[theorieGrl]{\grl}}
\newcommand \gRls {\hyperref[theorieGrl]{\grls}}

\newcommand\gtr{géné\-ra\-teur\xspace}  
\newcommand\gtrs{géné\-ra\-teurs\xspace}  


\newcommand \homo {ho\-mo\-mor\-phisme\xspace}
\newcommand \homos {ho\-mo\-mor\-phismes\xspace}

\newcommand \hmg {homo\-gène\xspace}
\newcommand \hmgs {homo\-gènes\xspace}

\newcommand \icftr {intervalle compact réticulé \ftm réel\xspace}
\newcommand \icftrs {intervalles compacts réticulés \ftm réels\xspace}

\newcommand \icl {inté\-gra\-lement clos\xspace}
\newcommand \icle {inté\-gra\-lement close\xspace}

\newcommand \icsr {intervalle compact \stm réticulé\xspace}
\newcommand \icsrs {intervalles compacts \stm réticulés\xspace}

\newcommand \icrc {intervalle compact réel clos\xspace}
\newcommand \icrcs {intervalles compact réels clos\xspace}

\newcommand \id {idéal\xspace}
\newcommand \ids {idéaux\xspace}

\newcommand \ida {\idt \agq}
\newcommand \idas {\idts \agqs}

\newcommand \idc  {\idt de Cramer\xspace}
\newcommand \idcs {\idts de Cramer\xspace}

\newcommand \idd {idéal déter\-minan\-tiel\xspace}
\newcommand \idds {idéaux déter\-minan\-tiels\xspace}

\newcommand \idema {idéal maxi\-mal\xspace}
\newcommand \idemas {idéaux maxi\-maux\xspace}

\newcommand \idep {idéal pre\-mier\xspace}
\newcommand \ideps {idéaux pre\-miers\xspace}

\newcommand \idemi {\idep minimal\xspace}
\newcommand \idemis {\ideps minimaux\xspace}

\newcommand \idf {idéal de Fitting\xspace}
\newcommand \idfs {idéaux de Fitting\xspace}

\newcommand \idif {idéal \dvlg fini\xspace}
\newcommand \idifs {idéaux \dvlgs finis\xspace}

\newcommand \idli {idéal \dvli\xspace} 
\newcommand \idlis {idéaux \dvlis\xspace} 

\newcommand \idm {idem\-po\-tent\xspace}
\newcommand \idms {idem\-po\-tents\xspace}
\newcommand \idme {idem\-po\-tente\xspace}
\newcommand \idmes {idem\-po\-tentes\xspace}

\newcommand \idp {idéal prin\-ci\-pal\xspace}
\newcommand \idps {idé\-aux prin\-ci\-paux\xspace}

\newcommand \idt {iden\-ti\-té\xspace}
\newcommand \idts {iden\-ti\-tés\xspace}

\newcommand \idtr {in\-dé\-ter\-mi\-née\xspace}
\newcommand \idtrs {in\-dé\-ter\-mi\-nées\xspace}

\newcommand \ifr {idéal frac\-tion\-nai\-re\xspace}
\newcommand \ifrs {idéaux frac\-tion\-nai\-res\xspace}

\newcommand \imd {immé\-diat\xspace}
\newcommand \imde {immé\-diate\xspace}
\newcommand \imds {immé\-diats\xspace}
\newcommand \imdes {immé\-diates\xspace}

\newcommand \imdt {immé\-dia\-te\-ment\xspace}

\newcommand \indtr {inf-demi-treillis\xspace} 

\newcommand \inteq {intui\-ti\-vement \eqve}
\newcommand \inteqs {intui\-ti\-vement \eqves}

\newcommand \Inteq {\hyperref[defextintequiv]{\inteq}}
\newcommand \Inteqs {\hyperref[defextintequiv]{\inteqs}}

\newcommand \ing {in\-ver\-se \gne}
\newcommand \ings {in\-ver\-ses \gnes}

\newcommand \iMP {in\-ver\-se de Moo\-re-Pen\-ro\-se\xspace}
\newcommand \iMPs {in\-ver\-ses de Moo\-re-Pen\-ro\-se\xspace}

\newcommand \ipp {\idep poten\-tiel\xspace}
\newcommand \ipps {\ideps poten\-tiels\xspace}

\newcommand \ird {irré\-duc\-tible\xspace}
\newcommand \irds {irré\-duc\-tibles\xspace}

\newcommand \iso {iso\-mor\-phisme\xspace}
\newcommand \isos {iso\-mor\-phismes\xspace}

\newcommand \itf {idéal \tf}
\newcommand \itfs {idé\-aux \tf}

\newcommand \itid {intui\-ti\-vement iden\-tique\xspace}
\newcommand \itids {intui\-ti\-vement iden\-tiques\xspace}

\newcommand \iv {inversible\xspace}
\newcommand \ivs {inversibles\xspace}

\newcommand \ivdg {inverse divisoriel\xspace} 
\newcommand \ivdgs {inverses divisoriels\xspace} 

\newcommand \ivde {inverse divisorielle\xspace} 
\newcommand \ivdes {inverses divisorielles\xspace} 

\newcommand \ivda {inverse divisoriel\xspace} 
\newcommand \ivdas {inverses divisoriels\xspace} 


\newcommand \lgb {local-global\xspace}
\newcommand \lgbe {locale-globale\xspace}
\newcommand \lgbs {local-globals\xspace}

\newcommand \LHe {Lemme de Hensel\xspace}
\newcommand \lHe {lemme de Hensel\xspace}
\newcommand \lHes {lemmes de Hensel\xspace}
\newcommand \LHm {\LHe multivarié\xspace}
\newcommand \lHm {\lHe multivarié\xspace}
\newcommand \lHms {\lHes multivariés\xspace}

\newcommand \lin {liné\-aire\xspace}
\newcommand \lins {liné\-aires\xspace}

\newcommand \lint {liné\-ai\-rement\xspace}

\newcommand \lnl {\lot \nl}
\newcommand \lnls {\lot \nls}

\newcommand \lot {loca\-lement\xspace}

\newcommand \lon {loca\-li\-sation\xspace}
\newcommand \lons {loca\-li\-sations\xspace}

\newcommand \lop {\lot prin\-cipal\xspace}
\newcommand \lops {\lot prin\-cipaux\xspace}

\newcommand \lsdz {\lot \sdz}

\newcommand \mdi {mo\-dule des \diles}

\newcommand \mlm {mo\-dule \lmo}
\newcommand \mlms {mo\-dules \lmos}

\newcommand \mlmo {ma\-tri\-ce de loca\-li\-sation
mono\-gène\xspace}
\newcommand \mlmos {ma\-tri\-ces de loca\-li\-sation
mono\-gène\xspace}

\newcommand \mlp {ma\-tri\-ce de loca\-li\-sation
prin\-ci\-pa\-le\xspace}
\newcommand \mlps {ma\-tri\-ces de loca\-li\-sation
prin\-ci\-pa\-le\xspace}

\newcommand \mo {mo\-no\"{\i}de\xspace}
\newcommand \mos {mo\-no\"{\i}des\xspace}

\newcommand \moco {\mos \com}

\newcommand \molo {morphisme de \lon\xspace}
\newcommand \molos {morphismes de \lon\xspace}

\newcommand \mom {mo\-nô\-me\xspace}
\newcommand \moms {mo\-nô\-mes\xspace}

\newcommand \moquo {morphisme de passage au quotient\xspace}
\newcommand \moquos {morphismes de passage au quotient\xspace}

\newcommand \mpf {mo\-dule \pf}
\newcommand \mpfs {mo\-dules \pf}

\newcommand \mpl {mo\-dule plat\xspace}
\newcommand \mpls {mo\-dules plats\xspace}

\newcommand \mpn {ma\-trice de \pn}
\newcommand \mpns {ma\-trices de \pn}

\newcommand \mprn {ma\-trice de \prn}
\newcommand \mprns {ma\-trices de \prn}

\newcommand \mptf {mo\-dule \ptf}
\newcommand \mptfs {mo\-dules \ptfs}

\newcommand \mrc {mo\-dule \prc}
\newcommand \mrcs {mo\-dules \prcs}

\newcommand \mtf {\aqo \tf}
\newcommand \mtfs {\aqos \tf}


\newcommand \ncr{néces\-saire\xspace}       
\newcommand \ncrs{néces\-saires\xspace}       

\newcommand \ncrt{néces\-sai\-rement\xspace}       

\newcommand \ndz {régu\-lier\xspace}
\newcommand \ndzs {régu\-liers\xspace}

\newcommand \nl {simple\xspace}
\newcommand \nls {simples\xspace}

\newcommand \noco {\noe \coh}
\newcommand \nocos {\noes \cohs}

\newcommand \Noe {Noether\xspace}

\newcommand \noe {noethé\-rien\xspace}
\newcommand \noes {noethé\-riens\xspace}
\newcommand \noee {noethé\-rienne\xspace}
\newcommand \noees {noethé\-riennes\xspace}

\newcommand \noet {noethé\-ria\-nité\xspace}

\newcommand \nst {Null\-stellen\-satz\xspace}
\newcommand \nsts {Null\-stellen\-s\"atze\xspace}

\newcommand \op{opé\-ra\-tion\xspace}  
\newcommand \ops{opé\-ra\-tions\xspace}
\newcommand \opari{\op\ari}  
\newcommand \oparis{\ops\aris}  
\newcommand \oparisv{\ops\arisv}  

\newcommand \oqc {ouvert \qc}
\newcommand \oqcs {ouverts \qcs}

\newcommand \ort{or\-tho\-go\-nal\xspace}  
\newcommand \orte{or\-tho\-go\-na\-le\xspace}  
\newcommand \orts{or\-tho\-go\-naux\xspace}  
\newcommand \ortes{or\-tho\-go\-na\-les\xspace}  


\newcommand \pa {couple saturé\xspace}
\newcommand \pas {couples saturés\xspace}
 
\newcommand \paral{paral\-lèle\xspace}  
\newcommand \parals{paal\-lèles\xspace}  

\newcommand \paralm{paral\-lè\-lement\xspace}

\newcommand \pb{pro\-blè\-me\xspace}  
\newcommand \pbs{pro\-blè\-mes\xspace}  

\newcommand \peq {purement équa\-tion\-nelle\xspace}
\newcommand \peqs {purement équa\-tion\-nelles\xspace}

\newcommand \pf {de \pn finie\xspace}

\newcommand \pgn {polygone de Newton\xspace}
\newcommand \pgns {polygones de Newton\xspace}

\newcommand \plc {rési\-duel\-lement \zed}
\newcommand \plcs {rési\-duel\-lement \zeds}

\newcommand \Plg {Prin\-cipe \lgb}
\newcommand \plg {prin\-cipe \lgb}
\newcommand \plgs {prin\-cipes \lgbs}

\newcommand \plga {\plg abs\-trait\xspace}
\newcommand \plgas {\plgs abs\-traits\xspace}

\newcommand \Plgc {\Plg con\-cret\xspace}
\newcommand \plgc {\plg con\-cret\xspace}
\newcommand \plgcs {\plgs con\-crets\xspace}

\newcommand \pn {présen\-ta\-tion\xspace}
\newcommand \pns {présen\-ta\-tions\xspace}

\newcommand \pog {\pol \hmg\xspace}
\newcommand \pogs {\pols \hmgs\xspace}

\newcommand \Pol {Poly\-nôme\xspace}
\newcommand \Pols {Poly\-nômes\xspace}

\newcommand \pol {poly\-nôme\xspace}
\newcommand \pols {poly\-nômes\xspace}

\newcommand \polH {\pol de Hensel\xspace}
\newcommand \polHs {\pols de Hensel\xspace}

\newcommand \poll{poly\-nomial\xspace}  
\newcommand \polls{poly\-nomiaux\xspace}  
\newcommand \polle{poly\-no\-miale\xspace}  
\newcommand \polles{poly\-no\-miales\xspace}  

\newcommand \pollt{poly\-no\-mia\-lement\xspace}  

\newcommand \polfon {\pol fon\-da\-men\-tal\xspace}
\newcommand \polmu {\pol rang\xspace}
\newcommand \polmus {\pols rang\xspace}
\newcommand \polcar {\pol carac\-té\-ris\-tique\xspace}
\newcommand \polcars {\pols carac\-té\-ris\-tiques\xspace}
\newcommand \polmin {\pol mini\-mal\xspace}
\newcommand \polmins {\pols mini\-maux\xspace}

\newcommand \polu {\pol unitaire\xspace}
\newcommand \polus {\pols unitaires\xspace}

\newcommand \prc {\pro de rang constant\xspace}
\newcommand \prcs {\pros de rang constant\xspace}

\newcommand \prcc {prin\-ci\-pe de \rcc}
\newcommand \prca {prin\-ci\-pe de \rca}
\newcommand \prce {prin\-ci\-pe de \rce}

\newcommand \prmt {préci\-sé\-ment\xspace}
\newcommand \Prmt {Préci\-sé\-ment\xspace}

\newcommand \prn {pro\-jec\-tion\xspace}
\newcommand \prns {pro\-jec\-tions\xspace}

\newcommand \pro {pro\-jec\-tif\xspace}
\newcommand \pros {pro\-jec\-tifs\xspace}

\newcommand \prr {pro\-jec\-teur\xspace}
\newcommand \prrs {pro\-jec\-teurs\xspace}

\newcommand \Prt {Pro\-pri\-été\xspace}
\newcommand \Prts {Pro\-pri\-étés\xspace}
\newcommand \prt {pro\-pri\-été\xspace}
\newcommand \prts {pro\-pri\-étés\xspace}

\newcommand \ptf {\pro \tf}
\newcommand \ptfs {\pros \tf}

\newcommand \qc {quasi-compact\xspace}
\newcommand \qcs {quasi-compacts\xspace}

\newcommand \qi {qua\-si in\-tè\-gre\xspace}
\newcommand \qis {qua\-si in\-tè\-gres\xspace}

\newcommand \qnl {quasi-\nl}
\newcommand \qnls {quasi-\nls}

\newcommand \ralg {règle \agq}
\newcommand \ralgs {règles \agqs}

\newcommand \rav {racine virtuelle\xspace}
\newcommand \ravs {racines virtuelles\xspace}

\newcommand \rcc {\rcm con\-cret\xspace}
\newcommand \rca {\rcm abs\-trait\xspace}
\newcommand \rce {\rcc des é\-ga\-li\-tés\xspace}

\newcommand \rcm {recol\-lement\xspace}
\newcommand \rcms {recol\-lements\xspace}

\newcommand \rcv {recou\-vrement\xspace} 
\newcommand \rcvs {recou\-vrements\xspace}

\newcommand \rde {rela\-tion de dépen\-dance\xspace}
\newcommand \rdes {rela\-tions de dépen\-dance\xspace}

\newcommand \rdi {\rde inté\-grale\xspace}
\newcommand \rdis {\rdes inté\-grales\xspace}

\newcommand \rdl {\rde \lin}
\newcommand \rdls {\rdes \lins}

\newcommand \rdt {rési\-duel\-lement\xspace}

\newcommand \rdy {règle dyna\-mique\xspace}
\newcommand \rdys {règles dyna\-miques\xspace}

\def \red {règle directe\xspace}
\newcommand \reds {règles directes\xspace}

\newcommand \rex {\hyperref[defexistsimple]{règle exis\-ten\-tielle simple\xspace}}
\newcommand \rexs {\hyperref[defexistsimple]{règles exis\-ten\-tielles simples\xspace}}

\newcommand \rexri {\hyperref[defitdyexrig]{règle exis\-ten\-tielle rigide\xspace}}
\newcommand \rexris {\hyperref[defitdyexrig]{règles exis\-ten\-tielles rigides\xspace}}

\newcommand \rsim {règle de simplification\xspace}
\newcommand \rsims {règles de simplification\xspace}

\newcommand \rtl {réti\-culé\xspace}
\newcommand \rtls {réti\-culés\xspace}

\newcommand \rmq {\rcm de quotients\xspace} 
\newcommand \rvq {\rcv par quotients\xspace} 
\newcommand \rmqs {\rcms de quotients\xspace} %
\newcommand \rvqs {\rcvs par quotients\xspace} %

\newcommand \rpf {réduite-de-présen\-tation-finie\xspace}
\newcommand \rpfs {réduites-de-présen\-tation-finie\xspace}


\newcommand \sad {\salg dynamique\xspace}
\newcommand \sads {\salgs dynamiques\xspace}

\newcommand \sagq {semi\agq}
\newcommand \sagqs {semi\agqs}

\newcommand \sagc {\sagq continue\xspace}
\newcommand \sagcs {\sagqs continues\xspace}

\newcommand \salg {structure \agq}
\newcommand \salgs {structures \agqs}

\newcommand \scentrel {relation semi-implicative\xspace}
\newcommand \scentrels {relations semi-implicatives\xspace}

\newcommand \scf {schéma fini\-taire\xspace}
\newcommand \scfs {schémas fini\-taires\xspace}

\newcommand \scl {schéma \elr}
\newcommand \scls {schémas \elrs}

\newcommand \sdo {\sdr \orte}
\newcommand \sdos {\sdrs \ortes}

\newcommand \sdr {somme directe\xspace}
\newcommand \sdrs {sommes directes\xspace}

\newcommand \sdz {sans \dvz}

\newcommand \sfio {sys\-tème fondamental d'\idms ortho\-gonaux\xspace}
\newcommand \sfios {sys\-tèmes fondamentaux d'\idms ortho\-gonaux\xspace}

\newcommand \sgr {\sys \gtr}
\newcommand \sgrs {\syss \gtrs}

\newcommand \slgb {stricte\-ment \lgb}
\newcommand \slgbs {stricte\-ment \lgbs}

\newcommand \sli {\sys \lin}
\newcommand \slis {\syss \lins}

\newcommand \smq {symé\-trique\xspace}
\newcommand \smqs {symé\-triques\xspace}

\newcommand \spb {sépa\-rable\xspace}  
\newcommand \spbs {sépa\-rables\xspace}

\newcommand \spe {spéci\-fi\-cation\xspace}
\newcommand \spes {spéci\-fi\-cations\xspace}

\newcommand \spi {\spe incomplète\xspace}
\newcommand \spis {\spes incomplètes\xspace}

\newcommand \spl {sépa\-rable\xspace}  
\newcommand \spls {sépa\-rables\xspace}

\newcommand \splc {\splt clos\xspace}
\newcommand \splce {\splt close\xspace}
\newcommand \splces {\splt closes\xspace}

\newcommand \splt {séparablement\xspace}  

\newcommand \spo {semipolynôme\xspace}
\newcommand \spos {semipolynômes\xspace}

\newcommand \spt{sépa\-ra\-bi\-lité\xspace}

\newcommand \srg {suite régu\-lière\xspace}
\newcommand \srgs {suites régu\-lières\xspace}

\newcommand \ste {strictement étale\xspace}
\newcommand \stes {strictement étales\xspace}

\newcommand \stf {strictement fini\xspace}
\newcommand \stfs {strictement finis\xspace}
\newcommand \stfe {strictement finie\xspace}
\newcommand \stfes {strictement finies\xspace}

\newcommand \stl {stablement libre\xspace}
\newcommand \stls {stablement libres\xspace}

\newcommand \stm {strictement\xspace}

\newcommand \str {\stm réticulé\xspace}
\newcommand \stre {\stm réticulée\xspace}
\newcommand \strs {\stm réticulés\xspace}
\newcommand \stres {\stm réticulées\xspace}

\newcommand \sul {supplé\-men\-taire\xspace}
\newcommand \suls {supplé\-men\-taires\xspace}

\newcommand \Sut {Support\xspace}
\newcommand \Suts {Supports\xspace}
\newcommand \sut {support\xspace}

\newcommand \syc {\sys de coordon\-nées\xspace}
\newcommand \sycs {\syss de coordon\-nées\xspace}

\newcommand \syp {\sys \poll}
\newcommand \Syp {\Sys \poll}
\newcommand \syps {\syss \polls}

\newcommand \sype {\syp étale\xspace}
\newcommand \sypes {\syps étales\xspace}

\newcommand \sys {système\xspace}
\newcommand \Sys {Système\xspace}
\newcommand \syss {systèmes\xspace}

\newcommand \sysN {\sys de Newton\xspace}
\newcommand \SysN {\Sys de Newton\xspace}
\newcommand \sysNs {\syss de Newton\xspace}

\newcommand \sysNe {\sysN étale\xspace}
\newcommand \sysNes {\sysNs étales\xspace}

\newcommand \talg {théorie \agq}
\newcommand \talgs {théories \agqs}

\newcommand \tco {théorie cohé\-rente\xspace}
\newcommand \tcos {théories cohé\-rentes\xspace}

\newcommand \tdy {théorie dyna\-mique\xspace}
\newcommand \tdys {théories dyna\-miques\xspace}

\newcommand \tel {\hyperref[defexistsimple]{théorie exis\-ten\-tielle\xspace}}
\newcommand \tels {\hyperref[defexistsimple]{théories exis\-ten\-tielles\xspace}}

\newcommand \telri {\hyperref[defitdyexrig]{théorie exis\-ten\-tielle rigide\xspace}}
\newcommand \telris {\hyperref[defitdyexrig]{théories exis\-ten\-tielles rigides\xspace}}

\newcommand \tf {de type fini\xspace}

\newcommand \tfo {théorie formelle\xspace}
\newcommand \tfos {théorie formelles\xspace}

\newcommand \tgm {théorie \gmq}
\newcommand \tgms {théories \gmqs}

\newcommand \Tho {Théo\-rème\xspace}
\newcommand \Thos {Théo\-rèmes\xspace}
\newcommand \tho {théo\-rème\xspace}
\newcommand \thos {théo\-rèmes\xspace}

\newcommand \thoc {théo\-rème$\mathbf{^*}$~}

\newcommand \tpe {théorie \peq}
\newcommand \tpes {théories \peqs}

\newcommand \trdi {treil\-lis dis\-tri\-bu\-tif\xspace}
\newcommand \trdis {treil\-lis dis\-tri\-bu\-tifs\xspace}

\newcommand \trel {trans\-for\-mation \elr}
\newcommand \trels {trans\-for\-mations \elrs}


\newcommand \ultm {ultramétrique\xspace}
\newcommand \ultms {ultramétriques\xspace}

\newcommand \umd {unimo\-du\-laire\xspace}
\newcommand \umds {unimo\-du\-laires\xspace}

\newcommand \unt {uni\-taire\xspace}
\newcommand \unts {uni\-taires\xspace}

\newcommand \uvl {uni\-versel\xspace}
\newcommand \uvle {uni\-ver\-selle\xspace}
\newcommand \uvls {uni\-versels\xspace}
\newcommand \uvles {uni\-ver\-selles\xspace}


\newcommand \vala {valeur absolue\xspace}
\newcommand \valas {valeurs absolues\xspace}

\newcommand \valn {valuation\xspace}
\newcommand \valns {valuations\xspace}

\newcommand \valu {\vala \ultm}
\newcommand \valus {\valas \ultms}

\newcommand \vfn {véri\-fi\-cation\xspace}
\newcommand \vfns {véri\-fi\-cations\xspace}

\newcommand \vmd {vec\-teur \umd}
\newcommand \vmds {vec\-teurs \umds}

\newcommand \vst {Valuativ\-stel\-lensatz\xspace}
\newcommand \vsts {Valuativ\-stel\-lensätze\xspace}

\newcommand \vstf {\vst formel\xspace}


\newcommand \zeH {zéro de Hensel\xspace}
\newcommand \zeHs {zéros de Hensel\xspace}

\newcommand \zed {z\'{e}ro-di\-men\-sionnel\xspace}
\newcommand \zede {z\'{e}ro-di\-men\-sion\-nelle\xspace}
\newcommand \zeds {z\'{e}ro-di\-men\-sion\-nels\xspace}
\newcommand \zedes {z\'{e}ro-di\-men\-sion\-nelles\xspace}

\newcommand \zedr {\zed réduit\xspace}
\newcommand \zedre {\zede réduite\xspace}
\newcommand \zedrs {\zeds réduits\xspace}

\newcommand \zmt {\tho de Zariski-Grothen\-dieck\xspace}


\newcommand \cof {cons\-truc\-tif\xspace}
\newcommand \cofs {cons\-truc\-tifs\xspace}

\newcommand \cov {cons\-truc\-tive\xspace}
\newcommand \covs {cons\-truc\-tives\xspace}

\newcommand \coma {\maths\covs}
\newcommand \clama {\maths clas\-siques\xspace}

\renewcommand \cot {cons\-truc\-ti\-vement\xspace}

\newcommand \matn {mathé\-ma\-ticien\xspace}
\newcommand \matne {mathé\-ma\-ti\-cienne\xspace}
\newcommand \matns {mathé\-ma\-ticiens\xspace}
\newcommand \matnes {mathé\-ma\-ti\-ciennes\xspace}

\newcommand \maths {mathé\-ma\-tiques\xspace}
\newcommand \mathe {mathé\-ma\-tique\xspace}

\newcommand \prco {démons\-tration \cov}
\newcommand \prcos {démons\-trations \covs}

\newcommand {\junk}[1]{}

\newcommand\rouge[1]{\textcolor{red}{#1}}
\newcommand\bleu[1]{\textcolor{blue}{#1}}
\newcommand\violet[1]{\textcolor{magenta}{#1}}

\newcommand{\Cadre}[2]{%

\medskip%
\newskip\oldleftskip
\newskip\oldrightskip
\oldleftskip=\leftskip%
\oldrightskip=\rightskip%
\leftskip=-\tabcolsep%
\rightskip=-\tabcolsep%
\begin{center}\fbox{%
\begin{tabular}%
{p{#1\textwidth}}
\setlength{\parindent}{5mm}%
\vspace{-1.5mm}#2\vspace{1mm}%
\end{tabular}}\end{center}\par\medskip%
\leftskip=\oldleftskip%
\rightskip=\oldrightskip%
\setlength{\parindent}{6mm}}

\newcommand\boite[2]{\begin{minipage}[c]{#1cm}
     \centering {#2} \end{minipage}}
\newcommand\Boite[3]{\parbox[t][#1cm][c]{#2cm}{\boite{#2}{#3}}}

\newcommand{\Encadre}[1]{\Cadre{.8}{#1}}

\newcommand{\Cencadre}[1]{\Encadre{\vspace{-3mm}
\begin{center}
#1 \end{center}\vspace{-8mm}}}

\newcommand{\cen}{\centerline}
\newcommand \Grandcadre[1]{%
\begin{center}
\begin{tabular}{|c|}
\hline
~\\[-3mm]
#1\\[-3mm]
~\\
\hline
\end{tabular}
\end{center}}

\newcommand\dsp{\displaystyle}
\newcommand\ndsp{\textstyle}

\newcommand{\eop}{\hfill \mbox{$\Box$}}

\newcommand \noi {\noindent}
\renewcommand \ss {\smallskip}
\newcommand \sni {\ss\noi}
\newcommand \snii {\noi}
\newcommand \ms {\medskip}
\newcommand \mni {\ms\noi}
\newcommand \bs {\bigskip}
\newcommand \bni {\bs\noi}
\newcommand \hs {\qquad}
\newcommand \alb {\allowbreak}
\newcommand \ce {\centerline}


\renewcommand \le{\leqslant}
\renewcommand \leq{\leqslant}
\renewcommand \preceq{\preccurlyeq}
\renewcommand \ge{\geqslant}
\renewcommand \geq{\geqslant}
\renewcommand \succeq{\succcurlyeq}
\newcommand   \nneq {\mathrel{\#}}
\newcommand   \ineq {$_{\,\mathrel{\#}}$}

\newcommand\eti{^\times}
\newcommand \epr{^\perp}
\newcommand \etl{^*}
\newcommand \sta{^\star}
\newcommand \bu {{$\bullet$}}
\newcommand \bl {^\bullet}
\newcommand{\bul}{^{\bullet}}
\newcommand \eci {^\circ}
\newcommand \uci{\mathring}
\newcommand \ep[1]{^{(#1)}}
\newcommand \esh{^\sharp}
\newcommand \efl{^\flat}
\newcommand \eto{$^*$ }
\newcommand \etoz{$^*$}
\newcommand \ist{_\star}

\newcommand \Ast {\gA^{\!\star}}
\newcommand \Bst {\gB^{\star}}
\newcommand \iBA {_{\gB/\!\gA}}
\newcommand \iWV {_{\gW\!/\gV}}
\newcommand \Bo{\BB\mathrm{o}}
\newcommand \Ati {\gA^{\!\times}}
\newcommand \Bti {\gB^{\times}}
\newcommand \Vti {\gV^{\times}}
\newcommand \Atl {\gA^{\!*}}
\newcommand \Btl {\gB^{*}}
\newcommand \Ktl {\gK^{*}}
\newcommand \Vtl {\gV^{*}}
\newcommand{\KAt}{\gK\etl\!\sur{\Ati}}
\newcommand{\AAt}{\Atl\!\sur{\Ati}}

\newcommand \divi {\mid}

\newcommand\equidef{\buildrel{{\rm def}}\over{\;\Longleftrightarrow\;}}
\newcommand\eqdef{\buildrel{\rm def}\over {\;=\;}}
\newcommand\eqdefi{\buildrel{\rm def}\over {\;=\;}}


\newcommand \fraC[2] {{{#1}\over {#2}}}
\newcommand \formule[1]{{\left\{ {\arraycolsep2pt\begin{array}{lll} #1 \end{array}}\right.}}
\newcommand \formul[1]{{\left\{ {\arraycolsep2pt\begin{array}{rcl} #1 \end{array}}\right.}}
\newcommand \formu[1]{\arraycolsep2pt\begin{array}{rcl} #1 \end{array}}

\newcommand\mapright[1]{\smash{\mathop{\longrightarrow}\limits^{#1}}}
\newcommand\maprightto[1]{\smash{\mathop{\longmapsto}\limits^{#1}}}
\newcommand\mapdown[1]{\downarrow\rlap{$\vcenter{\hbox{$\scriptstyle
#1$}}$}}
\newcommand{\pref}[1]{\textup{\hbox{\normalfont(\ref{#1})}}}

\newcommand \abs[1] {\left|{#1}\right|}
\newcommand \abS[1] {\big|{#1}\big|}
\newcommand \aqo[2] {#1\sur{\gen{#2}}\!}
\newcommand \aQo[2] {#1/{\gEn{#2}}\!}
\newcommand \Aqo[2] {#1\sur{\big\langle{#2}\big\rangle}\!}
\newcommand \Al[1] {\Vi^{\!#1}}
\newcommand \ci[1] {{{#1}^\circ}}
\newcommand \crac[2] {\cro {\frac{#1}{#2}}}
\newcommand \cro[1] {\left[#1\right]}
\newcommand \eqdf[1] {\buildrel{#1}\over =}
\newcommand \equivdf[1] {\buildrel{#1}\over \longleftrightarrow}
\newcommand \frt[1] {\!\left|_{#1}\right.\!}
\newcommand \impdef[1] {\buildrel{#1}\over \Longrightarrow}
\newcommand \norme[1]{\left\lVert #1 \right\rVert}
\newcommand \Norme[1]{\big\lVert #1 \big\rVert}
\newcommand \tra[1] {{\,^{\rm t}\!#1}}
\newcommand \gen[1] {\left\langle{#1}\right\rangle}
\newcommand \gEn[1] {\langle{#1}\rangle}
\newcommand \geN[1] {\big\langle{#1}\big\rangle}
\newcommand \sing[1] {\left\{{#1}\right\}}
\newcommand \so[1] {\left\{\,{#1}\, \right\}}
\newcommand \soo[1] {\{\,{#1}\,\}}
\newcommand \sO[1]{\big\{{#1}\big\}}
\newcommand \sotq[2]{\so{#1\mathrel{;}#2}}
\newcommand \sootq[2]{\soo{#1\mathrel{;}#2}}
\newcommand \sotQ[2]{\sO{#1\mathrel{;}#2}}
\newcommand \sur[1] {\!\left/#1\right.}
\newcommand \und[1] {\underline{#1}}

\newcommand \Sqr {\mathrm{Sqr}}

\newcommand \idg[1] {|\,#1\,|}
\newcommand \idG[1] {\big|\,#1\,\big|}

\newcommand \norm[1] {\Vert\,#1\,\Vert}

\newcommand \dex[1] {[\,#1\,]}
\newcommand \deX[1] {\big[\,#1\,\big]}

\newcommand \lst[1] {[\,#1\,]}
\newcommand \lsT[1] {\big[\,#1\,\big]}

\newcommand{\mt}{\mapsto}

\newcommand{\llongrightarrow}{\relbar\joinrel\mkern-1mu\longrightarrow}
\newcommand{\lllongrightarrow}{\relbar\joinrel\mkern-1mu\llongrightarrow}
\newcommand{\llllongrightarrow}{\relbar\joinrel\mkern-1mu\lllongrightarrow}
\newcommand\simarrow{\vers{_\sim}}
\newcommand\isosim{\buildrel{_\sim}\over \longleftrightarrow }
\newcommand\vers[1]{\buildrel{#1}\over \longrightarrow }
\newcommand\vvers[1]{\buildrel{#1}\over \llongrightarrow }
\newcommand\vvvers[1]{\buildrel{#1}\over \lllongrightarrow }
\newcommand \lora {\longrightarrow}
\newcommand \llra {\llongrightarrow}
\newcommand \lllra {\lllongrightarrow}

\renewcommand \leq{\leqslant}
\renewcommand \geq{\geqslant}

\newcommand \som {\sum\nolimits}
\newcommand \Ex {{\exists}}
\newcommand \Tt {{\forall}}
\newcommand \te {\otimes}
\newcommand \vep{{\varepsilon}}


\newcommand\lra[1]{\langle{#1}\rangle}
\newcommand\lrb[1] {\llbracket #1 \rrbracket}
\newcommand\lrbd {\lrb{1..d}}
\newcommand\lrbn {\lrb{1..n}}
\newcommand\lrbl {\lrb{1..\ell}}
\newcommand\lrbm {\lrb{1..m}}
\newcommand\lrbk {\lrb{1..k}}
\newcommand\lrbp {\lrb{1..p}}
\newcommand\lrbq {\lrb{1..q}}
\newcommand\lrbr {\lrb{1..r}}
\newcommand\lrbs {\lrb{1..s}}


\newcommand \vda {\,\vdash\,}

\newcommand \vdw {\,\vdash_w\,}

\newcommand\Aq[2]{#1_{[#2]}}
\newcommand\Aqj[2]{#1_{\{#2\}}}

\newcommand \dar[1] {\MA{\downarrow \!#1}}
\newcommand \uar[1] {\MA{\uparrow \!#1}}
\newcommand \clps[1] {{\downarrow #1\,\downarrow}}

\newcommand\tsbf[1]{\textsf{\textbf{\textup{#1}}}}
\newcommand\lab[1]{\item[\tsbf{#1}]}
\newcommand\Lab[1]{\rdb\item[\tsbf{#1}]\label{Ax#1}}
\newcommand\Tsbf[1]{\hyperref[Ax#1]{\tsbf{#1}}}

\newcommand\SA[1]{\rdb\sa{#1}\label{theorie#1}}
\newcommand\Sa[1]{\hyperref[theorie#1]{\sa{#1}}}
\newcommand\sa[1]{\hbox{\usefont{T1}{pzc}{m}{it}#1}\,}
\newcommand\sab[1]{\hbox{\usefont{T1}{pzc}{m}{it}#1}\,\,}
\newcommand\sA[1]{\hbox{\small\usefont{T1}{pzc}{m}{it}#1}\,}
\newcommand\sAb[1]{\hbox{\small\usefont{T1}{pzc}{m}{it}#1}\,\,}

\newcommand \snic[1] {\sni\centerline{$#1$}

\ss}

\newcommand \snac[1]{\sni
{\small\centering$#1$\par}

\ss}

\newcommand \eoe {\hbox{}\nobreak\hfill
\vrule width .5em height .5em depth 0mm \par \smallskip}

\newcommand \bal[1] {^\rK_{#1}}
\newcommand \ul[1] {_\rK^{#1}}

\newcommand \ov[1] {\overline{#1}}

\newcommand \wh[1] {\widehat{#1} }
\newcommand \wi[1] {\widetilde{#1} }

\newcommand\dessus[2]{{\textstyle {#1} \atop \textstyle {#2}}}

\newcommand\carray[2]{{\left[\begin{array}{#1} #2 \end{array}\right]}}
\newcommand\cmatrix[1]{\left[\matrix{#1}\right]}
\newcommand\clmatrix[1]{{\left[\begin{array}{lllllll} #1 \end{array}\right]}}
\newcommand\dmatrix[1]{\abs{\matrix{#1}}}
\newcommand\Cmatrix[2]{\setlength{\arraycolsep}{#1}\left[\matrix{#2}\right]}
\newcommand\Dmatrix[2]{\setlength{\arraycolsep}{#1}\left|\matrix{#2}\right|}

\newcommand \bloc[4]{\left[
\begin{array}{cc}
#1 & #2   \\
#3 & #4
\end{array}
\right]}

\newcommand \cm{em}

\newcommand{\blocs}[8]{%
{\setlength{\unitlength}{.0833\textwidth}
\tabcolsep0pt\renewcommand{\arraystretch}{0}%
\begin{tabular}{|c|c|}
\hline
\parbox[t][#3\cm][c]{#1\cm}{\begin{minipage}[c]{#1\cm}
\centering#5
\end{minipage}}&
\parbox[t][#3\cm][c]{#2\cm}{\begin{minipage}[c]{#2\cm}
\centering#6
\end{minipage}}\\
\hline
\parbox[t][#4\cm][c]{#1\cm}{\begin{minipage}[c]{#1\cm}
\centering#7
\end{minipage}}&
\parbox[t][#4\cm][c]{#2\cm}{\begin{minipage}[c]{#2\cm}
\centering#8
\end{minipage}}\\
\hline
\end{tabular}
}}

\newcommand \UneCol[1]{%
\sni\mbox{\hspace{.02\textwidth}%
\parbox[t]{.98\textwidth}{#1}%
}}

\newcommand \Unecol[1]{%
\sni\mbox{\hspace{.1\textwidth}%
\parbox[t]{.9\textwidth}{#1}%
}}

\newcommand \Deuxcol[4]{%
\sni\mbox{\parbox[t]{#1\textwidth}{#3}%
\hspace{.05\textwidth}%
\parbox[t]{#2\textwidth}{#4}}}

\newcommand \DeuxCol[2]{%
\Deuxcol{.475}{.475}{#1}{#2}}

\newcommand \DeuxCols[2]{%
\sni\mbox{\hspace{.02\textwidth}%
\parbox[t]{.475\textwidth}{#1}%
\hspace{.03\textwidth}%
\parbox[t]{.475\textwidth}{#2}}}

\newcommand \DeuxRegles[2]{%
\vspace{-1em}\DeuxCols
{\begin{enumerate}  #1
\end{enumerate}
}
{\begin{enumerate}  #2
\end{enumerate}
}
\vspace{-.3em}
}

\newcommand \UneRegle[2]{%
\vspace{-1em}\UneCol{
\begin{enumerate}
\lab{#1}{#2}
\end{enumerate}
}
\vspace{-.3em}
}

\newcommand \Regles[1]{%
\vspace{-1em}\UneCol{
\begin{enumerate}
{#1}
\end{enumerate}
}
\vspace{-.3em}
}

\newcommand \regles[1]{%
\vspace{-1em}\Unecol{
\begin{enumerate}
{#1}
\end{enumerate}
}
\vspace{-.3em}
}

\newcommand \itbu {\item[$\bullet$]}
\newcommand \labu {\lab{$\bullet$}}

\newcommand \Deuxbu[2]{%
\vspace{-.6em}\DeuxCols
{\begin{itemize}  #1
\end{itemize}}
{\begin{itemize}  #2
\end{itemize}}}

\makeatletter
\def\revddots{\mathinner{\mkern1mu\raise\p@
\vbox{\kern7\p@\hbox{.}}\mkern2mu
\raise4\p@\hbox{.}\mkern2mu\raise7\p@\hbox{.}\mkern1mu}}
\makeatother

\newcommand \BB{\mathbb {B}}
\newcommand \CC{\mathbb {C}}
\newcommand \FF{\mathbb {F}}
\newcommand \II{\mathbb {I}}
\newcommand \KK{\mathbb {K}}
\newcommand \MM{\mathbb {M}}
\newcommand \NN{\mathbb {N}}
\newcommand \ZZ{\mathbb {Z}}
\newcommand \OO{\mathbb {O}}
\newcommand \PP{\mathbb {P}}
\newcommand \QQ{\mathbb {Q}}
\newcommand \RR{\mathbb {R}}

\newcommand \gk {\mathbf{k}}
\newcommand \gkb {\ov\gk}
\newcommand \gl {\mathbf{l}}
\newcommand \gA {\mathbf{A}}
\newcommand \gB {\mathbf{B}}
\newcommand \gC {\mathbf{C}}
\newcommand \gD {\mathbf{D}}
\newcommand \gd {\mathbf{d}}
\newcommand \gDd {\mathbf{Dd}}
\newcommand \gE {\mathbf{E}}
\newcommand \gF {\mathbf{F}}
\newcommand \gG {\mathbf{G}}
\newcommand \gI {\mathbf{I}}
\newcommand \gIo {\mathbf{Io}}
\newcommand \gK {\mathbf{K}}
\newcommand \gKb {\ov\gK}
\newcommand \gAh {\widehat\gA}
\newcommand \gKh {\widehat\gK}
\newcommand \gKt {\widetilde\gK}
\newcommand \rhe {^{\mathrm{h}}}
\newcommand \rHe {^{\mathrm{H}}}
\newcommand \rNe {^{\mathrm{N}}}
\newcommand \rhs {^{\mathrm{hs}}}
\newcommand \Khe {\gK\rhe}
\newcommand \KHe {\gK\rHe}
\newcommand \KNe {\gK\rNe}
\newcommand \Kxi {\gK[\xi]}
\newcommand \gKp {{\gK'}}
\newcommand \gKw {\widetilde\gK}
\newcommand \gL {\mathbf{L}}
\newcommand \gLb {\ov\gL}
\newcommand \gLh {\widehat\gL}
\newcommand \gLw {\widetilde\gL}
\newcommand \Lhe {\gL\rhe}
\newcommand \gM {\mathbf{M}}
\newcommand \gP {\mathbf{P}}
\newcommand \gR {\mathbf{R}}
\newcommand \gRa {\mathbf{R}_\mathrm{a}}
\newcommand \gS {\mathbf{S}}
\newcommand \gT {\mathbf{T}}
\newcommand \gV {\mathbf{V}}
\newcommand \Vhe {\gV\rhe}
\newcommand \VHe {\gV\rHe}
\newcommand \VNe {\gV\rNe}
\newcommand \fmhe {\fm\rhe}
\newcommand \fmHe {\fm\rHe}
\newcommand \mhe {\fmhe}
\newcommand \mHe {\fmHe}
\newcommand \fmNe {\fm\rNe}
\newcommand \Ahe {\gA\!\rhe}
\newcommand \fmhs {\fm\rhs}
\newcommand \Ahs {\gA\rhs}
\newcommand \fmh {\widehat\fm}

\newcommand \KV {(\gK,\gV)}
\newcommand \LW {(\gL,\gW)}

\newcommand \fmA {{\fm_\gA}}
\newcommand \fmB {{\fm_\gB}}
\newcommand \fmC {{\fm_\gC}}
\newcommand \fmV {{\fm_\gV}}

\newcommand \gVh {\widehat\gV}
\newcommand \fmVh {\widehat{\fm_\gV}}
\newcommand \fmAh {\widehat{\fm_\gA}}
\newcommand \gAt {\widetilde\gA}
\newcommand \gVt {\widetilde\gV}
\newcommand \fmti {\widetilde\fm}
\newcommand \gVp {{\gV'}}
\newcommand \gW {\mathbf{W}}
\newcommand \gX {\mathbf{X}}
\newcommand \gZ {\mathbf{Z}}

\newcommand \Gat {\wh{\Gamma}}
\newcommand \Ksep {{\gK^\mathrm{sep}}}
\newcommand \ksep {{\gk^\mathrm{sep}}}
\newcommand \Kac {{\gK^\mathrm{ac}}}
\newcommand \Vsep {{\gV^\mathrm{sep}}}
\newcommand \fmsep {{\fm^\mathrm{sep}}}
\newcommand \vsep {{v_\mathrm{sep}}}
\newcommand \Vac {{\gV^\mathrm{ac}}}
\newcommand \fmac {{\fm^\mathrm{ac}}}

\newdimen\xyrowsp
\xyrowsp=3pt
\newcommand{\SCO}[6]{
\xymatrix @R = \xyrowsp {
                                  &1 \ar@{-}[dl] \ar@{-}[dr] \\
#3 \ar@{-}[ddr]                   &   & #6 \ar@{-}[ddl] \\
                                  &\bullet\ar@{-}[d] \\
                                  &\bullet   \\
#2 \ar@{-}[ddr] \ar@{-}[uur]      &   & #5 \ar@{-}[ddl] \ar@{-}[uul] \\
                                  &\bullet \ar@{-}[d] \\
                                  &\bullet  \\
#1 \ar@{-}[uur]                   &   & #4 \ar@{-}[uul] \\
                                  & 0 \ar@{-}[ul] \ar@{-}[ur] \\
}
}

\newcommand \Adj {\MA{\mathrm{Adj}}}
\newcommand \adj {\MA{\mathrm{adj}}}
\newcommand \Adu {\MA{\mathrm{Adu}}}
\newcommand \Ann {\mathrm{Ann}}
\newcommand \Atom {\mathrm{Atom}}
\newcommand \Aut {\MA{\mathrm{Aut}}}
\newcommand \BZ {\MA{\mathrm{BZ}}}
\newcommand \car {\MA{\mathrm{car}}}
\newcommand \Cl {\MA{\mathrm{Cl}}}
\newcommand \ClW {\MA{\mathrm{Cl}_{\mathrm{W}}}}
\newcommand \Coker {\MA{\mathrm{Coker}}}
\newcommand \Cont{\mathrm{Co}}
\newcommand \Dc {\MA{\mathrm{Dc}}}
\newcommand \DDiv {\MA{\mathrm{Dv}}}
\renewcommand \det {\MA{\mathrm{det}}}
\renewcommand \deg {\MA{\mathrm{deg}}}
\newcommand \Diag {\MA{\mathrm{Diag}}}
\newcommand \disc {\MA{\mathrm{disc}}}
\newcommand \Disc {\MA{\mathrm{Disc}}}
\newcommand \Div {\MA{\mathrm{Div}}}
\newcommand \DivA {\Div\gA }
\newcommand \DivAp {(\Div\gA)^{+} }
\newcommand \DivB {\Div\gB }
\newcommand \DivBp {(\Div\gB)^{+} }
\newcommand \DkM {\MA{\mathrm{DkM}}}
\newcommand \dv {\MA{\mathrm{div}} }
\newcommand \dvA {\dv_\gA }
\newcommand \dvB {\dv_\gB }
\newcommand \ev {{\mathrm{ev}}}
\newcommand \End {\MA{\mathrm{End}}}
\newcommand \fsac {\MA{\mathrm{fsa}}}
\newcommand \Fix {\MA{\mathrm{Fix}}}
\newcommand \Frac {\MA{\mathrm{Frac}}}
\newcommand \Gal {\MA{\mathrm{Gal}}}
\newcommand \Gfr {\MA{\mathrm{Gfr}}}
\newcommand \Gr {\MA{\mathrm{Gr}}}
\newcommand \gr {\MA{\mathrm{gr}}}
\newcommand \Gram {\MA{\mathrm{Gram}}}
\newcommand \gram {\MA{\mathrm{gram}}}
\newcommand \Grl {\MA{\mathrm{Grl}}}
\newcommand \hauteur {\mathrm{hauteur}}
\newcommand \Hom {\MA{\mathrm{Hom}}}
\newcommand \Id {\MA{\mathrm{Id}}}
\newcommand \Iv {\MA{\mathrm{Iv}}}
\newcommand \I {\mathrm{I}}
\newcommand \Idif {\MA{\mathrm{Idif}}}
\newcommand \Idv {\MA{\mathrm{Idv}}}
\newcommand \Ifr {\MA{\mathrm{Ifr}}}
\newcommand \Icl {\MA{\mathrm{Icl}}}
\renewcommand \Im {\MA{\mathrm{Im}}}
\newcommand \Inf {\MA{\mathrm{Inf}}}
\newcommand \Itf {\MA{\mathrm{Itf}}}
\newcommand \Ker {\MA{\mathrm{Ker}}}
\newcommand \Lst {\MA{\mathrm{Lst}}}
\newcommand \LIN {\mathrm{Lin}}
\newcommand \Lsf {\MA{\mathrm{Lsf}}}
\newcommand \Mat {\MA{\mathrm{Mat}}}
\newcommand \Mip {\mathrm{Min}}
\newcommand \md {\mathrm{md}}
\newcommand \Mgcd {\MA{\mathrm{Mgcd}}}
\renewcommand \mod {\;\mathrm{mod}\;}
\newcommand \Mor {\MA{\mathrm{Mor}}}
\newcommand \NDc {\MA{\mathrm{NDc}}}
\newcommand \poids {\mathrm{poids}}
\newcommand \poles {\hbox {\rm p\^oles}}
\newcommand \pgcd {\MA{\mathrm{pgcd}}}
\newcommand \ppcm {\MA{\mathrm{ppcm}}}
\newcommand \Rad {\MA{\mathrm{Rad}}}
\newcommand \Reg {\MA{\mathrm{Reg}}}
\newcommand \rg{\MA{\mathrm{rg}}}
\newcommand \rgst {\mathrm{rgst}}
\newcommand \Res {\mathrm{Res}}
\newcommand \Rs {\MA{\mathrm{Rs}}}
\newcommand \rPr{\MA{\mathrm{Pr}}}
\newcommand \Rv {\mathrm{Rv}}
\newcommand \Sli {\MA{\mathrm{Sli}}}
\newcommand \Som {\MA{\mathrm{Som}}}
\newcommand \Sup {\MA{\mathrm{Sup}}}
\newcommand \Sace {\MA{\mathrm{Sace}}}
\newcommand \Smtf {\MA{\mathrm{Smtf}}}
\newcommand \Stp {\MA{\mathrm{Stp}}}
\newcommand \St {\mathrm{St}}
\newcommand \Tri {\MA{\mathrm{Tri}}}
\newcommand \Tor {\MA{\mathrm{Tor}}}
\newcommand \tr {\MA{\mathrm{tr}}}
\newcommand \Tr {\MA{\mathrm{Tr}}}
\newcommand \Tsc {\MA{\mathrm{Tsch}}}
\newcommand \Um {\MA{\mathrm{Um}}}
\newcommand \val {\MA{\mathrm{val}}}

\newcommand \SIPD {\MA{\mathrm{SIPD}}}
\newcommand \ARC {\MA{\mathrm{ARC}}}
\newcommand \AFR {\MA{\mathrm{AFR}}}
\newcommand \AFRRV {\MA{\mathrm{AFRRV}}}
\newcommand \AFRNZ {\MA{\mathrm{AFRNZ}}}
\newcommand \PPM {\MA{\mathrm{PPM}}}
\newcommand \PB {\MA{\mathrm{PB}}}

\newcommand \Suslin{{\rm Suslin}}
\newcommand{\DBxk}{{\Der \gk\gB\xi}}%
\newcommand{\DAxk}{{\Der \gk\gA\xi}}%
\newcommand{\DkXxk}{{\Der \gk\kuX\xi}}%
\newcommand \DAbul {\rD_{\!\Abul}}

\newcommand\MA[1]{\mathop{#1}\nolimits}

\newcommand \sfa {\mathsf{a}}
\newcommand \sfb {\mathsf{b}}
\newcommand \sfc {\mathsf{c}}
\newcommand \sfd {\mathsf{d}}
\newcommand \sfe {\mathsf{e}}
\newcommand \sff {\mathsf{f}}
\newcommand \sfx {\mathsf{x}}
\newcommand \sfy {\mathsf{y}}
\newcommand \sfz {\mathsf{z}}
\newcommand \sbv {\tsbf{v}}
\newcommand \sbw {\tsbf{w}}
\newcommand \sfv {\mathsf{v}}
\newcommand \sfu {\mathsf{u}}
\newcommand \sft {\mathsf{t}}
\newcommand \sfw {\mathsf{w}}

\newcommand \Cdim {\MA{\mathsf{Cdim}}}
\newcommand \Divdim {\MA{\mathsf{Divdim}}}
\newcommand \Glo {\MA{\mathsf{Glo}}}
\newcommand \GK {\MA{\mathsf{GK}}}
\newcommand \GKO {\MA{\mathsf{GK}_0}}
\newcommand \HO {\MA{\mathsf{H}_0}}
\newcommand \HOp {\MA{\mathsf{H}_0^+}}
\newcommand \Hdim {\MA{\mathsf{Hdim}}}
\newcommand \HeA {{\Heit\gA}}
\newcommand \Heit {\MA{\mathsf{Heit}}}
\newcommand \Hspec {\MA{\mathsf{Hspec}}}
\newcommand \Jdim {\MA{\mathsf{Jdim}}}
\newcommand \jdim {\MA{\mathsf{jdim}}}
\newcommand \Jspec {\MA{\mathsf{Jspec}}}
\newcommand \jspec {\MA{\mathsf{jspec}}}
\newcommand \KO {\MA{\mathsf{K}_0}}
\newcommand \KOp {\MA{\mathsf{K}_0^+}}
\newcommand \KTO {\wi{\mathsf{K}}_0}
\newcommand \Kdim {\MA{\mathsf{Kdim}}}
\newcommand \Lin {\mathsf{L}}
\newcommand \Max {\MA{\mathsf{Max}}}
\newcommand \Min {\MA{\mathsf{Min}}}
\newcommand \OQC {\MA{\mathsf{Oqc}}}
\newcommand \Pic {\MA{\mathsf{Pic}}}
\newcommand \Reel {\MA{\mathsf{Reel}}}
\newcommand \Spec {\MA{\mathsf{Spec}}}
\newcommand \Speclin {\MA{\mathsf{Speclin}}}
\newcommand \Spv {\MA{\mathsf{Spv}}}
\newcommand \Spev {\MA{\mathsf{Spev}}}
\newcommand \Sper {\MA{\mathsf{Sper}}}
\newcommand \Val {\MA{\mathsf{Val}}}
\newcommand \Valp {\MA{\mathsf{Val'}}}
\newcommand \SpecA {\Spec\gA}
\newcommand \SpevA {\Spev\gA}
\newcommand \SpvA {\Spv\gA}
\newcommand \SperA {\Sper\gA}
\newcommand \SpecT {\Spec\gT}
\newcommand \Zar {\MA{\mathsf{Zar}}}
\newcommand \ZF {\MA{\mathsf{ZF}}}
\newcommand \ValA {{\Val\gA}}
\newcommand \ValpA {{\Valp\gA}}
\newcommand \ZarA {{\Zar\gA}}

\newcommand \cA {{\cal A}}
\newcommand \cB {{\cal B}}
\newcommand \cC {{\cal C}}
\newcommand \cD {{\cal D}}
\newcommand \cI {{\cal I}}
\newcommand \cJ {{\cal J}}
\newcommand \cF {{\cal F}}
\newcommand \cH {{\cal H}}
\newcommand \cK {{\cal K}}
\newcommand \cL {{\cal L}}
\newcommand \cM {{\cal M}}
\newcommand \cN {{\cal N}}
\newcommand \cP {{\cal P}}
\newcommand \cQ {{\cal Q}}
\newcommand \cR {{\cal R}}
\newcommand \cS {{\cal S}}
\newcommand \cT {{\cal T}}
\newcommand \cV {{\cal V}}

\newcommand \ccd{\mathcal{CD}}
\newcommand \cco{\mathcal{CO}}

\newcommand \Cin{C^{\infty}}

\newcommand \SK {\cS^\rK}
\newcommand \IK {\cI^\rK}
\newcommand \JK {\cJ^\rK}
\newcommand \IH {\cI\rHe}
\newcommand \JH {\cJ\rHe}

\newcommand \JAC {J}
\newcommand \Jac {\mathrm{Jac}}
\newcommand \rc {\mathrm{c}}
\newcommand \Df {\MA{\mathrm{Df}}}
\newcommand \Dfr {\MA{\mathrm{Df}}^\mathrm{R}}
\newcommand \Dfmc {\MA{\mathrm{Dfmc}}}
\newcommand \rd {\mathrm{d}}
\newcommand \rv {\mathrm{v}}
\newcommand \rI {\mathrm{I}}
\newcommand \Ic {\mathrm{Ic}}
\newcommand \rC {\mathrm{C}}
\newcommand \rD {\mathrm{D}}
\newcommand \rF {\mathrm{F}}
\newcommand \rG {\mathrm{G}}
\newcommand \rH {\mathrm{H}}
\newcommand \rJ {\mathrm{J}}
\newcommand \Li {\MA{\mathrm{Li}}}
\newcommand \rK {\mathrm{K}}
\newcommand \rL {\mathrm{L}}
\newcommand \Mc {\mathrm{Mc}}
\newcommand \rN {\mathrm{N}}
\newcommand \rP {\mathrm{P}}
\newcommand \rR {\mathrm{R}}
\newcommand \rmSa {\MA{\mathrm{Sa}}}
\newcommand \rmSamc {\MA{\mathrm{Samc}}}
\newcommand \rS {\mathrm{S}}
\newcommand \rU {\mathrm{U}}
\newcommand \rV {\mathrm{V}}
\newcommand \DA {\rD_{\!\gA}}
\newcommand \JA {\rJ_\gA}
\newcommand \JT {\rJ_\gT}

\newcommand\fa{\mathfrak{a}}
\newcommand\fb{\mathfrak{b}}
\newcommand\fc{\mathfrak{c}}
\newcommand\fA{\mathfrak{A}}
\newcommand\fB{\mathfrak{B}}
\newcommand\fD{\mathfrak{D}}
\newcommand\fI{\mathfrak{i}}
\newcommand\fII{\mathfrak{I}}
\newcommand\fj{\mathfrak{j}}
\newcommand\fJ{\mathfrak{J}}
\newcommand\fF{\mathfrak{F}}
\newcommand\ff{\mathfrak{f}}
\newcommand\ffg{\mathfrak{g}}
\newcommand\fG{\mathfrak{G}}
\newcommand\fh{\mathfrak{h}}
\newcommand\fl{\mathfrak{l}}
\newcommand\fm{\mathfrak{m}}
\newcommand\mV{{\fm_\gV}}
\newcommand\mW{{\fm_\gW}}
\newcommand\fM{\mathfrak{M}}
\newcommand\fp{\mathfrak{p}}
\newcommand\fP{\mathfrak{P}}
\newcommand\fq{\mathfrak{q}}
\newcommand\fU{\mathfrak{U}}
\newcommand\fV{\mathfrak{V}}
\newcommand\fx{\mathfrak{x}}
\newcommand\fy{\mathfrak{y}}

\newcommand \scC{\mathscr{C}}
\newcommand \scR{\mathscr{R}}

\newcommand{\bma}{\bm{a}}
\newcommand{\bmb}{\bm{b}}
\newcommand{\bmc}{\bm{c}}
\newcommand{\bmd}{\bm{d}}
\newcommand{\bme}{\bm{e}}
\newcommand{\bmf}{\bm{f}}
\newcommand{\bmu}{\bm{u}}
\newcommand{\bmv}{\bm{v}}
\newcommand{\bmw}{\bm{w}}
\newcommand{\bmy}{\bm{y}}
\newcommand{\bmx}{\bm{x}}
\newcommand{\bmz}{\bm{z}}

\newcommand \LLPO {\tsbf{LLPO}}
\newcommand \LPO  {\tsbf{LPO}}

\newcommand \Zg {{\Z[G]}}

\newcommand \vu {\vee} 
\newcommand \vi {\wedge} 
\newcommand \Vu {\bigvee}
\newcommand \Vi {\bigwedge}
\newcommand \im {\rightarrow} 
\newcommand \da {\,\downarrow\!}

\newcommand \vdu[1] {\vdash^{#1}}
\newcommand \vdb[1] {\vdash_{#1}}

\newcommand \Vrai {\mathsf{Vrai}}
\newcommand \Faux {\mathsf{Faux}}
\newcommand \Un {\mathbf{1}}
\newcommand \Deux {\mathbf{2}}
\newcommand \Trois {\mathbf{3}}
\newcommand \Quatre {\mathbf{4}}
\newcommand \Cinq {\mathbf{5}}

\newcommand \una {{\underline{a}}}
\newcommand \ual {{\underline{\alpha}}}
\newcommand \ua  {{\underline{a}}}
\newcommand \ub  {{\underline{b}}}
\newcommand \ube {{\underline{\beta}}}
\newcommand \uc  {{\underline{c}}}
\newcommand \ud  {{\underline{d}}}
\newcommand \udel{{\underline{\delta}}}
\newcommand \ue  {{\underline{e}}}
\newcommand \uf  {{\underline{f}}}
\newcommand \uF  {{\underline{F}}}
\newcommand \ug  {{\underline{g}}}
\newcommand \uh  {{\underline{h}}}
\newcommand \uga {{\underline{\gamma}}}
\newcommand \uP  {{\underline{P}}}
\newcommand \ur{{\underline{r}}}
\newcommand \ut{{\underline{t}}}
\newcommand \uu{{\underline{u}}}
\newcommand \ux {{\underline{x}}}
\newcommand \uxi {{\underline{\xi}}}
\newcommand \uX {\underline{X}}
\newcommand \uy{{\underline{y}}}
\newcommand \uY  {{\underline{Y}}}
\newcommand \uz{{\underline{z}}}
\newcommand \uze {{\underline{0}}}

\newcommand \ak {a_1,\ldots,a_k}
\newcommand \am {a_1,\ldots,a_m}
\newcommand \an {a_1,\ldots,a_n}
\newcommand \aq {a_1,\ldots,a_q}
\newcommand \aln {\alpha_1,\ldots,\alpha_n}
\newcommand \bn {b_1,\ldots,b_n}
\newcommand \bzn {b_0,\ldots,b_n}
\newcommand \bbm {b_1,\ldots,b_m}
\newcommand \ck {c_1,\ldots,c_k}
\newcommand \rcr {c_1,\ldots,c_r}
\newcommand \cq {c_1,\ldots,c_q}
\newcommand \gan {\gamma_1,\ldots,\gamma_n}
\newcommand \un {u_1,\ldots,u_n}
\newcommand \xk {x_1,\ldots,x_k}
\newcommand \Xk {X_1,\ldots,X_k}
\newcommand \xm {x_1,\ldots,x_m}
\newcommand \Xm {X_1,\ldots,X_m}
\newcommand \fn {f_1,\ldots,f_n}
\newcommand \lfm {f_1,\ldots,f_m}
\newcommand \Fn {F_1,\ldots,F_n}
\newcommand \lFm {F_1,\ldots,F_m}
\newcommand \Fp {\FF_p}
\newcommand \Zp {\ZZ_p}
\newcommand \Qp {\QQ_p}
\newcommand \Cp {\CC_p}
\newcommand \sn {s_1,\ldots,s_n}
\newcommand \xn {x_1,\ldots,x_n}
\newcommand \xzn {x_0,\ldots,x_n}
\newcommand \xhn {x_0:\ldots:x_n}
\newcommand \Xn {X_1,\ldots,X_n}
\newcommand \xr {x_1,\ldots,x_r}
\newcommand \Xr {X_1,\ldots,X_r}
\newcommand \xin {\xi_1,\ldots,\xi_n}
\newcommand \xizn {\xi_0,\ldots,\xi_n}
\newcommand \xihn {\xi_0:\ldots:\xi_n}
\newcommand \ym {y_1,\ldots,y_m}
\newcommand \yr {y_1,\ldots,y_r}
\newcommand \Yr {Y_1,\ldots,Y_r}
\newcommand \Yn {Y_1,\ldots,Y_n}
\newcommand \Ym {Y_1,\ldots,Y_m}
\newcommand \yn {y_1,\ldots,y_n}

\newcommand \AT {\gA[T]}
\newcommand \AX {\gA[X]}
\newcommand \Ax {\gA[x]}
\newcommand \Ared {\gA\red}
\newcommand \AuX {\gA[\uX]}
\newcommand \Aux {\gA[\ux]}
\newcommand \ArX {\gA\lra X}
\newcommand \Axn {\gA[\xn]}
\newcommand \AXn {\gA[\Xn]}

\newcommand \AXm {\gA[\Xm]}
\newcommand \KKXm {{\KK[\Xm]}}
\newcommand \KKuX {{\KK[\uX]}}

\newcommand \AY {\gA[Y]}
\newcommand \Ayn {\gA[\yn]}

\newcommand \BuX {\gB[\uX]}
\newcommand \BuY {\gB[\uY]}
\newcommand \BX {{\gB[X]}}
\newcommand \BT {{\gB[T]}}
\newcommand \BY {{\gB[Y]}}
\newcommand \Bxn {\gB[\xn]}
\newcommand \BXn {{\gB[\Xn]}}
\newcommand \BYm {\gB[\Ym]}

\newcommand \CT {{\gC[T]}}
\newcommand \CX {{\gC[X]}}
\newcommand \CXn {{\gC[\Xn]}}

\newcommand \kX {{\gk[X]}}
\newcommand \KX {\gK[X]}
\newcommand \Kx {\gK[x]}
\newcommand \VX {\gV[X]}
\newcommand \Vx {\gV[x]}
\newcommand \VT {\gV[T]}
\newcommand \KT {\gK[T]}
\newcommand \KuX {\gK[\uX]}
\newcommand \kuX {\gk[\uX]}
\newcommand \VuX {\gV[\uX]}
\newcommand \Vuxi {\gV[\uxi]}
\newcommand \Vux {\gV[\ux]}
\newcommand \Kux {\gK[\ux]}
\newcommand \kux {\gk[\ux]}
\newcommand \KXk {\gK[\Xk]}
\newcommand \KXm {\gK[\Xm]}
\newcommand \KXn {\gK[\Xn]}
\newcommand \kxm {\gk[\xm]}
\newcommand \kxn {\gk[\xn]}
\newcommand \Ky {\gK[y]}
\newcommand \KY {\gK[Y]}
\newcommand \Kz {\gK[z]}
\newcommand \KZ {\gK[Z]}
\newcommand \kZ {\gk[Z]}
\newcommand \KYn {\gK[\Yn]}
\newcommand \KYm {\gK[\Ym]}
\newcommand \kXr {\gk[\Xr]}
\newcommand \KXr {\gK[\Xr]}
\newcommand \Kxr {\gK[\xr]}

\newcommand \Vxn {\gV[\xn]}
\newcommand \VXn {\gV[\Xn]}

\newcommand \KuY {\gK[\uY]}
\newcommand \Kuy {\gK[\uy]}
\newcommand \Kyn {\gK[\yn]}
\newcommand \Kyr {\gK[\yr]}
\newcommand \kYr {\gk[\Yr]}
\newcommand \KYr {\gK[\Yr]}
\newcommand \Kxn {\gK[\xn]}

\newcommand \LuX {\gL[\uX]}
\newcommand \lXn {\gl[\Xn]}
\newcommand \lxn {\gl[\xn]}
\newcommand \LXn {\gL[\Xn]}
\newcommand \lXr {\gl[\Xr]}
\newcommand \LXr {\gL[\Xr]}
\newcommand \lYr {\gl[\Yr]}
\newcommand \LYr {\gL[\Yr]}

\newcommand \QQXn {\QQ[\Xn]}

\newcommand \Rx {\gR[x]}
\newcommand \Rux {\gR[\ux]}
\newcommand \RuX {\gR[\uX]}
\newcommand \RXk {{\gR[\Xk]}}
\newcommand \RXm {{\gR[\Xm]}}
\newcommand \Rxm {{\gR[\xm]}}
\newcommand \RXn {{\gR[\Xn]}}
\newcommand \Rxn {{\gR[\xn]}}
\newcommand \RXzn {{\gR[\Xzn]}}
\newcommand \Rxzn {{\gR[\xzn]}}
\newcommand \RXr {\gR[\Xr]}

\newcommand \Ruy {\gR[\uy]}
\newcommand \Ryn {{\gR[\yn]}}
\newcommand \RRX {\RR[X]}
\newcommand \RRXn {\RR[\Xn]}
\newcommand \RRxn {\RR[\xn]}
\newcommand \RYr {\gR[\Yr]}
\newcommand \RRuX {\RR[\uX]}
\newcommand \RRux {\RR[\ux]}

\newcommand \CCX {\CC[X]}

\newcommand \ZZXn {\ZZ[\Xn]}
\newcommand \ZG {\ZZ[G]}

\newcommand \RRXk {{\RR[\Xk]}}
\newcommand \RRxk {{\RR[\xk]}}
\newcommand \RRXm {{\RR[\Xm]}}
\newcommand \RRxm {{\RR[\xm]}}

\newcommand \lfs {f_1,\ldots,f_s}
\newcommand \lfn {f_1,\ldots,f_n}

\newcommand \Gn  {\gG_n}
\newcommand \Gnk {\gG^n_{k}}
\newcommand \Gnr {\gG^n_{r}}
\newcommand \cGn {\cG_n}
\newcommand \cGnk{\cG_{n,k}}
\newcommand \GGn {\GG^n}
\newcommand \GGnk{\GGn_{k}} 
\newcommand \GGnr{\GGn_{r}}
\newcommand \GA  {\mathbb{GA}}
\newcommand \GAn {\GA^n}  
\newcommand \GAq {\GA^q}
\newcommand \GAnk{\GAn_{k}}
\newcommand \GAnr{\GAn_{r}}
\newcommand \GL {\mathbb{GL}}
\newcommand \GLn {{\GL_n}}
\newcommand \SL {\mathbb{SL}}
\newcommand \SLn {{\SL_n}}
\newcommand \EE {\mathbb{E}}
\newcommand \En {\EE_n}
\newcommand \Pn {\PP^n}
\newcommand \An {\AA^n}
\newcommand \Sl {\mathbf{SL}}
\newcommand \Sln {{\Sl_n}}

\newcommand \Mm {\MM_{m}}
\newcommand \Mn {\MM_{n}}
\newcommand \Mk {\MM_{k}}
\newcommand \Mq {\MM_{q}}
\newcommand \Mr {\MM_{r}}
\newcommand \MMn {\MM_{n}}

\newcommand \Pf {{{\cal P}_{\mathrm{f}}}}
\newcommand \Pfe {{\rm P}_{{\rm fe}}}

\newcommand\hsz{\\ }
\newcommand\hsu{\\ \hspace*{4mm}}
\newcommand\hsd{\\ \hspace*{8mm}}
\newcommand\hst{\\ \hspace*{1,2cm}}
\newcommand\hsq{\\ \hspace*{1,6cm}}
\newcommand\hsc{\\ \hspace*{2cm}}
\newcommand\hsix{\\ \hspace*{2,4cm}}
\newcommand\hsept{\\ \hspace*{2,8cm}}

%


\renewcommand \sibrouillon[1]{}

\renewcommand \hum[1] {\sibrouillon{\noindent {\sf hum: #1}}}

\stMF

\title{Lemme de Hensel multivarié pour les corps ultramétriques}

\author{M.-E. Alonso\thanks{Universidad Complutense, Madrid, Espa\~na. {\tt mariemi@mat.ucm.es}}
 \and  Henri Lombardi
\thanks{Université de Franche-Comté, Laboratoire de mathématiques de Besançon, UMR CNRS 6623, 16 route de Gray, 25000
Besançon, France. {\tt henri.lombardi@univ-fcomte.fr}}
\and
Stefan Neuwirth
\thanks{Université de Franche-Comté, Laboratoire de mathématiques de Besançon, UMR CNRS 6623, 16 route de Gray, 25000
Besançon, France. {\tt stefan.neuwirth@univ-fcomte.fr}}}

\clearpage\setcounter{section}{0}\selectlanguage{french}\def\frenchproofname{\textsl{Démonstration}}
\emptythanks\setcounter{footnote}{0}
\setcounter{equation}{0}
\clearpage\setcounter{section}{0}\selectlanguage{french}\def\frenchproofname{\textsl{Démonstration}}

\date{\today}
\maketitle

\vspace{-1em}

\begin{abstract} 
Le \lHm affirme que sur un anneau local, tout système polynomial de Newton admet un zéro à \coos dans le hensélisé de l'anneau local.
Le \lHm pour les anneaux locaux est usuellement démontré comme une conséquence de la version Grothendieck du Zariski Main Theorem (ZMT). Ce ZMT traite une situation plus générale à priori beaucoup plus difficile.
Nous proposons ici une \prco de ce lemme indépendante du ZMT pour le cas des \cvdus. En \clama cela implique le lemme pour les corps valués de rang $1$.
\end{abstract}

\rdb
\label{beginfrench}

\medskip \noindent {\bf Mots clés:} \LHm, corps valué, corps valorisé discret, corps ultramétrique, hensélisé d’un anneau local, hensélisé d’un corps valué,  \coma. 

\smallskip  \noindent {\bf MSC2020:} 13B40, 13J15, 03F65

\setcounter{tocdepth}{4}
\markboth{Table des matières}{Table des matières}

\small
\printcontents[french]{}{1}{}
\newpage
\normalsize

\markboth{Introduction}{Introduction}
\Section{Introduction}
Cet article est écrit dans le cadre des \coma à la Bishop (\citealt*{fBi67,fBB85,fBR1987,fMRR,fCACM,fYen2015,fACMC}).

Nous nous situons dans la continuation naturelle des articles \citealt*{fCLR01,fKL00,fKLP03,fALP08} et, dans une moindre mesure, \citealt*{fALN2021,fCL2016,fCL2016b,fLM2022}.

\smallskip
Le lemme de Hensel apparait en \maths dans le cadre des travaux de Hensel sur les corps $p$-adiques $\Qp$ et leurs extensions finies. Il s’agit alors de corps valorisés munis d’une valeur absolue ultramétrique $x\mapsto \abs x, \,\gK\to\RR$ (voir la section \ref{fseccvdu}).

Dans les années 30  \cite{fKru1930} et \cite{fDeu1931} introduisent une \gnn des \cvus aux \cvals \gnls. Cette \gnn est \ncr pour obtenir le \tho fondamental de Krull qui affirme que la \cli d’un anneau intègre dans son \cdf est l’intersection des suranneaux qui sont des \advs. Les \cvus correspondent aux \cvals de rang $1$.

Une première notion de corps valorisé hensélien (un \cvu où tout \polH admet un \zeH) et du hensélisé d'un \cvu est introduite par Ostrowski  (avec une terminologie différente) dans l'article fondateur \citealt*{fOst1934}.
Ostrowski définit le hensélisé comme la \cli du corps dans son complété.
Pour plus de détails voir \url{https://www.mathi.uni-heidelberg.de/~roquette/manu.html#Valuations} et \citealt*{fRoq02}.
D'un point de vue moderne, l'étude d'Ostrowski se limite aux \cvals avec une valuation de rang~1. 

Bien qu’apparemment presque évident pour les \cvus dans la mesure où le zéro d'un \sysN peut être calculé par la méthode de Newton dans un complété du \cvar, il reste que les \coos de ce séro appartiennent au hensélisé s'avère difficile à démontrer même pour des \advds les plus simples. 
Des ouvrages de référence en \clama pour la théorie des \advs sont  \citealt*{fNagata62}, \citealt*{fBou-AC-5-6} et \citealt*{fEP2005}. Ils ne  prêtent cependant pas attention au \lHm, même dans les exercices.
Dans les années 50, la notion \gnle d’\alo hensélien est introduite par \cite{fAzu1951} et \cite{fNag1953}. Elle est devenue essentielle en géométrie algébrique.
  
Le \lHm pour les \alos est usuellement démontré comme une conséquence de la version Grothendieck du Zariski Main Theorem (ZMT). Cette version du~ZMT traite une situation plus générale à priori beaucoup plus difficile 
(pour un traitement \cof, voir \citealt*{fACL2014}).

\smallskip  Nous proposons dans cet article, pour le cas des \cvdus, une \prco de ce lemme indépendante du~ZMT, avec un contenu \algq précis. Cela fournit aussi une \demo en \clama pour les \advs de rang 1. 

\smallskip Dans la section \ref{fsecEtales} nous introduisons le \pb \gnl du \lHm pour les anneaux locaux.
Nous rappelons les notions de  \sysN, de \sysN étale et d'\alg étale. 
Le \thref{fthstrucste} rappelle des versions \covs des \thos de  structure pour les \algs \stfes étales sur un \cdi  (\citet*[\thos \hbox{VI-1.7} et \hbox{VI-1.9}]{fCACM}).  
Un outil décisif est alors le \tho de structure pour les \algs \pf nettes (non ramifiées) sur un \cdi (\thref{fthNewtonCorpsDiscret}) qui
affirme qu'une telle \alg est toujours étale et \stfe.
Une \prco se trouve dans \citet*[Corollary~\hbox{VI-6.15}]{fCACM}. Dans \citealt*{fACMC}, une \prco plus \elr est donnée à la fin de la section~VI-6.

\smallskip Dans la section \ref{fseccvdu} nous étudions le cas des \cvdus. 

\setcounter{section}{2}
\setcounter{subsection}{2}
\setcounter{ftheorem}{4}
Notre résultat crucial est le \tho suivant, dans lequel nous comparons le hensélisé  $(\KHe,\VHe)$ construit dans \citealt*{fKL00} avec~$(\gKt,\gVt)$, où $\gKt$ la clôture séparable de $\gK$ dans $\KHe$ et $\gVt=\gKt\cap\VHe$.
\begin{ftheorem}[deux versions \eqves du hensélisé d'un \cud] 
\label{fthMRRKLintro}~\\
On considère un  \cud $(\gK,\abs\cdot)$.
Le hensélisé  $(\KHe,\VHe)$ de $\KV$ est isomorphe à  $(\gKt,\gVt)$. Plus \prmt, il existe un unique $\gK$-\homo  $\KHe\to\gKt$ qui envoie~$\VHe$ dans $\gVt$, et cet \homo est un  \iso.
\end{ftheorem}

On a obtenu \prmt que tout \elt de $\gKt$ est l'image d'un  \elt $\gamma$
dans un corps $\gK[\xi]\subseteq \KHe$ où 
$\xi$ est le zéro spécial d'un \pol spécial.

Nous terminons en démontrant notre \lHm.

\setcounter{ftheorem}{5}
\begin{ftheorem}[\lHm  pour un \cud] \label{fLHM2intro} ~ \\
Soit $(\gK,\abs\cdot)$ un \cud et 
$(f_1,\dots,f_n)$ un \sysN sur $\KV$. Ce système admet un unique zéro à coordonnées dans $\fmti$.  Il admet aussi un unique zéro à coordonnées dans $\fm\VHe$.
\end{ftheorem}

\setcounter{section}{0}
\setcounter{subsection}{0}
\setcounter{ftheorem}{0}

\section{Systèmes de Newton, \algs étales}\label{fsecEtales}

\Subsubsection{Terminologie \cov}
En \coma un \textsl{\alo} est un anneau pour lequel, $\forall x$, $x$ \hbox{ou $1-x$} est \iv (avec un \gui{ou} explicite). 

Le \textsl{radical de Jacobson d’un anneau $\gA$} est l’\id $\Rad(\gA)=\sotq{x\in\gA}{1+x\gA\subseteq \Ati}$. 

Pour un \alo le radical de Jacobson est l’unique \idema\footnote{En \clama.}, noté \gnlt $\fm_\gA$ ou $\fm$. On dire simplement \textsl{l'anneau local $(\gA,\fm)$}.

Un \textsl{corps de Heyting} est un \alo non trivial dont le radical de Jacobson est nul. 

Un \textsl{\cdi} est un anneau non trivial dans lequel tout \elt est nul ou \iv. C’est la même chose qu’un corps de Heyting avec un test à $0$. 

Le \textsl{corps résiduel} d’un \alo non trivial $(\gA,\fm)$
est le corps de Heyting $\gA/\fm$ souvent noté~\(\ov\gA\).

L’\alo est dit \textsl{\dcd} si son corps résiduel est un \cdi. Cela revient à dire que nous avons explicitement la disjonction $x\in\Ati$ ou $x\in\fm_\gA$ pour tout $x\in \gA$. 

Pour deux \alos $(\gA,\fmA)$ et $(\gB,\fmB)$ un morphisme d’anneaux $\varphi\colon\gA\to\gB$ est dit \textsl{local} lorsque $\varphi^{-1}(\gB^\times)\subseteq \Ati$. Dans ce cas on dit que $(\gB,\fmB)$ est une $(\gA,\fmA)$-\alg.

\Subsubsection{Codes de Hensel sur un \alo}

Un \textsl{\cdH} sur un \alo $(\gA,\fm)$ est un couple $(f,a)\in\AX\times \gA$ où $f$ est un \polu, $f(a)\in \fm$ et $f'(a)\in \Ati$, i.e. $\ov a$ est un zéro simple de $\ov f$. On dit aussi dans ce cas que $f$ est un \textsl{\polH} ou un \textsl{\pol de Nagata}. 

Un \textsl{\pol spécial} est un \pol 
$h(X)=X^n-X^{n-1}+\sum_{k=n-2}^0 a_kX^k$ avec les $a_k\in\fm$. Il est clair dans ce cas que $(h,1)$ est un \cdH.

Un \textsl{\zeH du \cdH $(f,a)$} dans une $(\gA,\fmA )$-\alg $(\gB,\fmB)$ est un \elt $\xi\in\gB$ tel que $\varphi_\star(f)(\xi)=0$ et $\xi\in\varphi(a)+\fmB$.

Un \alo  $(\gA, \fm)$ est dit \textsl{hensélien} si tout code de Hensel $(f,a)$ a un \zeH $\alpha$ dans~$(\gA, \fm)$. On dit dans ce cas que $\alpha$ relève le zéro simple $a$ de $\ov f$ dans $\ov\gA$. Le \zeH pour le \cdH $(h,1)$ d'un \pol spécial $h$ est appelé le \textsl{zéro spécial} du \pol.
 
Un \zeH pour un code donné est \ncrt unique, sans besoin de supposer que $f$ est unitaire: on écrit 
\[0=f(\alpha')-f(\alpha)=f'(\alpha)\mu+b\mu^2=\mu\cdot(f'(\alpha)+b\mu), \quad \hbox{avec }b\in\gA\eqno(+)\]
 où $\mu=\alpha'-\alpha\in\fm$, et comme $f'(\alpha)+b\mu\in \Ati$ on obtient $\mu=0$. 

Le \textsl{complété} de l’\alo $(\gA,\fm)$ pour la topologie $\fm$-adique, \cad la limite projective des $\gA/\fm^k$ est noté $\gAh$.
On a un morphisme naturel $\varphi\colon\gA\to\gAh$ et on note $\fmAh=\varphi(\fm)\gAh$. On obtient un morphisme naturel  $(\gA,\fm)\to(\gAh,\fmAh)$. 
 Le morphisme~$\varphi$ est injectif si la topologie est séparée, \cad si $\bigcap_{n\in\NN}\fm^n=0$.

Si $(\gA,\fm)$ est un \alo \dcd,  $(\gAh,\fmh)$ est un \alo hensélien: 
la méthode de Newton (voir le \thref{fthNewtonQuad}) permet de calculer un \zeH pour n'importe quel \cdH.

\Subsubsection{Hensélisé d'un \alo}

La notion de hensélisé d'un \alo correspond à la solution du \pb \uvl associé à la sous-catégorie pleine des \alos henséliens de la catégorie des \alos (les flèches dans la catégorie sont les morphismes locaux).

L'article \citealt*{fALP08} construit ce hensélisé dans le cas d'un \alrd en ajoutant itérativement un \zeH formel pour tout \cdH.
L’ajout d’un \zeH formel se fait par une \textsl{extension de Nagata} de $\gA$ \citep*[Definition~6.1]{fALP08}. 
On note $\Ahe$ le \textsl{hensélisé} de $(\gA,\fmA)$.


\smallskip Voici maintenant un premier résultat qui étend la possibilité de relever un zéro simple au cas d'un \pol non \ncrt unitaire. 

Ce résultat parle d'un zéro d'un \pol $f\in\AX$ dans un \alo hensélien $(\gB,\fmB)$.
Le résultat est établi pour un \alrd et son hensélisé dans \cite[Lemma~5.3 et Proposition~5.4]{fALP08}.  
\begin{flemma}[l’astuce de Hervé] \label{flemTrick1} Soit $(\gA,\fm)$ un \alo, $\varphi\colon\gA\to\gB$ a morphisme local avec $(\gB,\fmB)$ hensélien, 
et $f(X)=\sum_{k=0}^na_kX^k\in\AX$ avec $a_0\in\fmA$ et $a_1\in\Ati$.
Alors, le \pol $f$ admet un zéro  $\gamma$ dans $a_0\gB\subseteq \fmB\subseteq \Rad(\gB)$, avec $f'(\gamma)$ inversible dans~$\gB$. C'est l'unique zéro de $f$
dans $\fmB$.
En particulier, si $(\gA,\fmA)$ est hensélien,~$f$ admet un zéro dans $a_0\gA\subseteq \fmA$, et c'est l'unique zéro de $f$ dans $\fmA$.
\end{flemma}
%
\begin{proof} On définit le \pol spécial 
\[
\formu{  
    g(X)&=&X^n-X^{n-1}+a_{0}\cdot
\big(\som_{j=2}^{n}(-1)^j   a_{j}a_{0}^{j-2} a_{1}^{-j}  X^{n-j}\big) 
\\[.4em]
&=& X^n-X^{n-1}+a_{0}\ell(X)  \quad \hbox{avec } \ell(X)\in\AX
  } 
\]
L'\egt suivante  a lieu dans $\gA[X,1/X]$
\[  
  a_0g(X) = X^nf\left(\frac{-a_0a_1^{-1}}X\right)  \eqno (*)
  \]   
Soit $\delta=1+\alpha$ où  $\alpha\in\fmB\subseteq \Rad(\gB)$ le zéro spécial du \pol spécial $g$. Alors $\delta\in
{\gB}^\times$. 
Posons $\gamma={-a_0  a_1^{-1}\over \delta}=-a_0  (a_1\delta)^{-1}\in\fmB$. En appliquant~$(*)$ nous avons
$-a_0  g(\delta)=\delta^n  f(\gamma)$, donc $f(\gamma)=0$. En outre $f'(\gamma)\in{\gB}^\times$ parce que $f'(0)\in\Ati$ et $\gamma\in\fmB=\Rad(\gB)$.

\noindent L'unicité est déjà établie plus haut, cf. $(+)$.
\end{proof}
%

\begin{fremark}  
Notez que
pour un \alrd $(\gA,\fmA)$ 
l'\elt $\delta$ apparait dans la construction de~$\Ahe$ en tant qu'\elt de l'anneau noté $\gA_g:=S^{-1}\Ax$  où~$\Ax=\aqo \AX g$ et $S=\sotq{s(x)\in\Ax}{s(1)\in\Ati}$. L'anneau $\gA_g$ est un \alrd \fpt sur $\gA$
et $\Rad(\gA_g)=\fm\gA_g$. L'\elt $\delta$ est juste l'image de $x\in\Ax$ dans $\gA_g$ via le morphisme de localisation.
Comme  le morphisme canonique $\gA\to\gA_g$ est injectif, on peut identifier~$\gA$ à un sous anneau de $\gA_g$. On constate alors que $\gA[\delta]$ est (isomorphe~à) un quotient de $\Ax$. Ce genre d'anneau $\gA_g$ est une brique \elr de la construction du hensélisé d'un \alrd: voir  \citet*[Definition~6.1]{fALP08}. \eoe
\end{fremark}

\subsection{\Syp de Newton dans un \alo}

Sur un \alo $(\gA,\fm)$, un \textsl{\sysN} (ou \textsl{de Hensel)  au point $(\ua)=(\an)$} de $\gA^n$  est un \syp $(\uf)=(f_1,\dots,f_n)\in{\AXn}^n$ qui admet $(\ua)$  comme
\textsl{zéro simple approché} au sens suivant: 
\begin{itemize}
\item les $f_j(\ua)$ sont dans $\fm$
\item  la matrice jacobienne du \sys est inversible au point $(\ua)$ modulo $\fm$.
%
\end{itemize}
La deuxième condition revient à dire que le \deter jacobien $\Jac(\ua)$ est inversible dans $\gA$.

Le \textsl{\lHm} classique dit qu'un \sysN  au point $(\ua)$ sur un \alo $(\gA,\fmA )$ admet un zéro $(\uxi)$ à \coos dans le hensélisé~$\Ahe$ de l'\alo avec les $\xi_j-a_j\in\fm\Ahe$.

Un \cdH est la même chose qu'un \sysN à une seule variable.

\smallskip Notons aussi que  \cite{fLafon63} et \citetalias[Section 15.11]{fstacks-project}  
donnent sous forme à peine cachée le \lHm pour les couples henséliens (par exemple l'implication $(5)\Rightarrow(2)$ dans le Lemma 15.11.6 de Stacks) mais qu'ils utilisent le ZMT pour leur \demo.

\Subsubsection{Unicité du \zeH}

\begin{flemma} \label{flem-Unicite-zeH}
Soit $(\gA,\fm)$ un \alo et $(\uf)=(f_1,\dots,f_n)$ un \sysN au point $(\ua)\in\gA^n$.  Si $(\ual)$ et $(\uga)$ sont des \zeHs au point $(\ua)$ pour ce \sys, alors $(\ual)=(\uga)$.      
\end{flemma}
\begin{proof} Notons $J(\uX)$ la matrice jacobienne du \sys.  On écrit $(\uga)=(\ual)+(\udel)$ avec les $\delta_i\in\fm$. 
On regarde $(\uf(\ual))=(f_i(\ual))_{i\in\lrbn}$, $(\uf(\uga))$ et $(\udel)$ comme des vecteurs colonnes.
Les formules de Taylor en plusieurs variables pour les \pols donnent une \egt
\[(\uf(\uga))=(\uf(\ual))+J(\ual)(\udel)+ M(\udel)
\]
avec la matrice carrée $M$ à \coes dans l'\id $\gen{{\udel}}$ et 
$(\uf(\uga))=(\uf(\ual))=(\uze)$.
\\
 Donc $(\rI_n+J(\ual)^{-1}M)(\udel)=(\uze)$, d'où
$(\udel)=(\uze)$.
\end{proof}

L'image d'un \sysN par un morphisme local est un \sysN.
L'unicité est prouvée de la même manière pour un zéro éventuel dans le but du morphisme local.

\subsection{La méthode de Newton}

Voici un énoncé précis pour ce qui est couramment appelé la méthode de Newton quadratique, dans un cadre purement \agq, cf. \citet*[Theorem III-10.3]{fCACM}.
\begin{ftheorem}[méthode de Newton quadratique] \label{fthNewtonQuad}~ \\
Considérons un \alo~$(\gA,\fm)$, et soit
 $(\uf)=(f_1,\dots,f_n)$ un \sysN dans  $\AXn$ au point $(\ua)=(a_1,\ldots ,a_n)\in\gA^n$. Notons $\JAC(\uxi)$ la matrice jacobienne de $(\uf)$ au point $(\uxi)$. Soit $U$ un inverse de $\JAC(\ua)$ modulo $\fm$.  
Définissons les suites $(\ua^{(m)})_{m\geq
0}$ dans $\gA^n$  et $(U^{(m)})_{m\geq 0}$ dans $\Mn(\gA)$
par l'itération de
Newton quadratique suivante:
$$\begin{array}{lcl}
\ua^{(0)}=\ua,&  \quad \quad  & \ua^{(m+1)}=\ua^{(m)}- U^{(m)} \cdot
\uf(\ua^{(m)}),    \\[1mm]
U^{(0)}=U,&   & U^{(m+1)}=U^{(m)}\,\left(2\I_n-\JAC(\ua^{(m+1)})U^{(m)}\right).
\end{array}$$
On obtient alors pour tout $m$ les congruences suivantes:
$$
\begin{array}{lcll}
 \ua^{(m+1)}\equiv\ua^{(m)} & \,\,\mathrm{et}\,\, &
  U^{(m+1)}\equiv U^{(m)} &\,\mod \,\fm^{2^m}    \\[1mm]
  \uf(\ua^{(m)})\equiv 0 &  \,\,\mathrm{et}\,\,  &
   U^{(m)}\,\JAC(\ua^{(m)})\equiv \rI_n  &\,\mod \,\fm^{2^m}.
\end{array}
$$
\end{ftheorem}

Ce \tho nous dit que si l'on a un \syp   $(\uf)=(f_1,\dots,f_n)$ dans $\AXn$ avec $(\ua)$ comme zéro simple $1$-approché pour  $(\gA,\fm)$, alors~$(\ua^{(m)})$ est  un zéro simple $2^m$-approché de $(\uf)$ avec $\ua^{(m)}\equiv\ua \mod \,\fm$ ($m> 1$).

En fait, le \tho précédent s’applique pour n’importe quel \id $\fm$ de n’importe quel anneau $\gA$.
Et dans la plupart des cas, l'\id $\fm$ à considérer est plutôt celui engendré par les $f_i(\ua)$, car il donne une description plus précise des approximations successives du zéro recherché (il se peut par exemple que $\fm^2=\fm$). 

\smallskip Nous utiliserons la terminologie suivante  en \coma: un \alo $(\gA,\fm)$  
est dit \textsl{quasi-\noe} lorsque $\gA$ est  \coh \fdi, $\fm$ est un \itf, et 
$\bigcap_{n\in\NN}\fm^n=0$\footnote{Dans certaines \prcos, on demande en outre que soit connu, pour n'importe quel $x\neq 0$, l'entier $k$ tel que $x\in\fm^k\setminus\fm^{k+1}$.}. Dans ce cas chaque $\fm^n$ est un \Amo \coh \fdi  et le morphisme naturel $\gA\to\gAh$ est injectif.
Cela est vérifié dans le cas suivant:  $\gA=\gB_{1+\fmB }$ où $\gB$ est une \alg \pf sur~$\ZZ$  ou sur un \cdi, $\fmB $ est un \itf, et~$\fmA =\fmB \gA$\footnote{Un \alo quasi-\noe est dit \noe si l'anneau $\gA$ est \noe. Nous n'utiliserons pas ici cette notion.}.

\smallskip Voici deux \crls du \thref{fthNewtonQuad}
\begin{fcorollary}[\lHm, version faible \num1] \label{fLHM1} ~\\ 
Soit $(\gK,\abs \cdot )$ un \cvdu\footnote{Voir la section \ref{fseccvdu}.} et soit $(f_1,\dots,f_n)$ un \sys \poll comme dans le \tho \ref{fthNewtonQuad}, avec $\gA=\gV=\sotq{x\in\gK}{\abs x \leq 1}$ et $\fm=\sotq{x\in\gK}{\abs x < 1}$. Alors le système admet un unique zéro $(\uxi)=(\xi_1,\dots,\xi_n)$ à coordonnées dans $\gKh$ satisfaisant $\xi_i-a_i\in \fmh$ pour $i=1,\dots,n$.  
\end{fcorollary}

\begin{fcorollary}[\lHm, version faible \num2] \label{fLHMquasinoet} ~\\
Soit $(\gA,\fm )$ un \alo quasi-noethérien et soit $(f_1,\dots,f_n)$ un \syp comme dans le \tho \ref{fthNewtonQuad}. Alors le système admet un unique zéro $(\uxi)=(\xi_1,\dots,\xi_n)$ à coordonnées dans $\gAh$ satisfaisant $\xi_i-a_i\in \fm\gAh$ pour $i=1,\dots,n$.  
\end{fcorollary}

\subsection{Algèbres étales}

Le contexte de la méthode de Newton est généralisé et formalisé sous le nom d'algèbre étale-de-base.

\begin{fdefinition} \label{fdefiBasicEtale} Soit~$\gA$ un anneau commutatif.
\begin{enumerate}
\item Soient
 $(\uf)= (f_1,\dots,f_n)$ un \sys de $n$ \pols dans  $\AXn$ et soit
 $\gB=\aQo\AuX \uf$. 
 L'\Alg $\gB=\Axn$ est dite \textsl{étale-de-base} si la matrice jacobienne $\JAC({\ux})$ du \sys $(\uf)$  est
\iv dans~$\gB$. On dit dans ce cas que le \sys $(f_1,\dots,f_n)$ est un \textsl{\sype}.
\item Une \Alg \pf $\gC=\aqo \AXm \lfs$ est dite \textsl{étale} si on~a des \eco $(u_i)_{i\in I}$ de $\gC$ tels que chacune des \Algs $\gC[1/u_i]$ est étale-de-base pour une \pn convenable.
%
%
\end{enumerate}
\end{fdefinition}

La notion d'\alg étale est un concept fondamental de l'\alg commutative, introduit par Grothendieck. 

On démontre \cot que toute \alg étale est 
étale-de-base pour une \pn convenable. Nous n'utiliserons pas ce résultat dans cet article.


Un produit fini d'\algs étales est étale.
Une localisée $\gC[1/s]$ d'une \alg étale est étale. En particulier l'\alg triviale est étale.

\Subsubsection{\SysN versus \sype}

\begin{flemma} \label{fnote0}
On peut toujours remplacer un \sysN $(f_1,\dots,f_n)$ sur un \alo $(\gA,\fm)$ 
par un \emph{\sysNe}, \cad qui est à la fois un \sysN sur $(\gA,\fm)$ et un \sype sur $\gA$. En outre
la notion de \zeH  dans une $(\gA,\fm)$-\alg est inchangée.
\end{flemma}

%
\begin{proof}
En effet, en notant $\Jac(\uX)$ le \deter jacobien, on ajoute une \idtr $X_{n+1}$ et le \pol $f_{n+1}:=(1+X_{n+1})\Jac(\uze)^{-1}\Jac(\uX)-1$.
On a $f_{n+1}(\uze,0)=0$. Alors le nouveau \syp a son \deter jacobien~$\Jac_1(\ux,x_{n+1})$ \iv dans la nouvelle \Alg quotient, et tout \zeH~$(\uxi)$ du premier \sys à \coos dans une \Alg donne un \zeH $(\uxi,\eta)$  du nouveau \sysN avec $(1+\eta)\Jac(\uze)^{-1}\Jac(\uxi)=1$, i.e. 
$\eta=\Jac(\uze)/\Jac(\uxi)-1$.

\end{proof}

Cette manipulation \elr permet de démontrer certaines \prts des \coos d’un \zeH. Par exemple, comme conséquence du \thref{fthNewtonCorpsDiscret}, si $(\gA,\fm)$ est un anneau local intègre et en notant~\hbox{$\gK=\Frac\gA$}, les \coos d’un \zeH d'un \sysN sur $(\gA,\fm)$ dans une \Klg sont toujours des \elts \spls sur $\gK$.

Dans le cas des \crls \ref{fLHM1} et \ref{fLHMquasinoet}, si l'on considère $\gB=\aQo\AuX\uf=\Aux$, on est certain que la \Alg 
$\gB[1/\Jac(\ux)]$,
qui est étale-de-base, n'est pas triviale parce qu'on~a un morphisme de $\gB[1/\Jac(\ux)]$ vers le complété $\gKh$, qui est une \Alg non triviale.

\Subsubsection{Structure des \algs étales sur un \cdi} 

On considère un \cdi $\gK$ et une extension $\gKp$ qui est un \Kev \tf. 
On dit alors que $\gKp$ est une \textsl{extension \agq finie} de $\gK$. 

Plus \gnlt une \Klg $\gC$ qui est un \Kev \tf est dite \textsl{finie sur $\gK$}. Les \elts de $\gC$ sont tous \agqs sur $\gK$, mais on ne connait pas forcément la dimension de $\gC$ comme \Kev.
Si $\gC$ est engendré par~$n$ \elts  comme \Kev, on écrit $[\gC:\gK]\leq n$.
Si on connait $m$ \elts $\gK$-\lint indépendants, on écrit $[\gC:\gK]\geq m$.
Enfin  $\gC$ est dite 
 \textsl{\stfe} sur $\gK$ lorsque que l'on connait une base de $\gC$ sur $\gK$,
\cad si la dimension de $\gC$ comme \Kev est connue, et l'on écrit alors $[\gC:\gK]= n$. Dans ce cas:
\begin{itemize}
\item on sait calculer le \polmin sur $\gK$ de tout \elt de $\gC$;
\item  si $\gD$
est une \Klg intermédiaire \tf, elle est \stfe sur $\gK$;
\item  si en outre $\gD$ est un corps, $\gC$ est \stfe sur $\gD$ et l'on a la relation classique $[\gC:\gK]= [\gC:\gD][\gD:\gK]$.
\end{itemize}

\begin{fdefinition}[\algs \stes sur un \cdi] \label{fdefiaste}~\\
Soit $\gK$  un \cdi. Une \Klg \stfe $\gB$ est dite \textsl{\ste} si  la forme trace $\phi(x,y)=\Tr_{\gB/\gK}(xy)\colon\gB\times \gB\to\gK$ est \textsl{non dégénérée}, ce qui signifie qu'en posant $\varphi (x):=\phi(x,\bullet)$, l'\Kli $\varphi$ définit un \iso du \Kev $\gB$ sur son dual\footnote{Ces \dfns se généralisent aux \algs sur un anneau commutatif arbitraire \cite[Theorem VI-5.5]{fCACM}.}.  
\end{fdefinition}

Un premier \tho de structure pour les \Klgs \stes est le suivant.
Une \demo \cov se trouve dans \citet*[théorèmes VI-1.7 et \hbox{VI-1.9}]{fCACM}.
En particulier on voit que les \algs \stes sur un \cdi sont des \algs étales \stfes. En outre le \thref{thNewtonCorpsDiscret} établit que tout \ale sur un \cdi est une \alg \ste.

Notez que le \tho \cof \ref{thstrucste} est plus précis que ses versions classiques, et que la \prco est plutôt subtile. En fait, les hypothèses sont données sous une forme qui permet d'obtenir les conclusions sous forme d'\algos.

\begin{ftheorem}[théorème de l'\elt primitif] \label{fthstrucste}~\\  
Soit $\gK$  un \cdi et soit $\gB$ une \Klg \stfe. 
\begin{enumerate}\itemsep=-.2em
\item 
\Propeq
\begin{enumerate}\itemsep=.01em
\item  $\gB$ est \ste.
\item $\gB$ est engendrée par des \elts \spls sur $\gK$.
\item  Tous les \elts de $\gB$ sont \spls sur $\gK$.
\item $\gB$ est isomorphe à un produit fini de \Klgs $\aqo{\KX}{h_i}$ pour des \pols  unitaires \spls $h_i$.
\end{enumerate}
En particulier toute \Klg \ste est étale.
\item Lorsque  $\gB$ est un \cdi ou si $\gK$ est infini, les \prts du point 1 sont aussi \eqves à: 
\begin{enumerate}\setcounter{enumii}{4}\itemsep=.01em
\item $\gB$ est isomorphe à une \Klg $\gK[\zeta]=\aqo \KZ g$ où 
$g$ est un \pol unitaire \spl de~$\KZ$. 
\end{enumerate}
\end{enumerate}
\end{ftheorem}

Dans le dernier cas, si $g$ se factorise sous la forme $g=g_1\cdots g_r$,  on a un \iso canonique $\gB\simeq \prod_{j=1}^r\aqo\KZ{g_j}$.

\smallskip
Une \Alg $\gB$  est dite \textsl{nette} si elle est \pf et si son module des \diles est nul. Le \Bmo des \diles est isomorphe au conoyau de la transposée de la matrice jacobienne (vue dans $\gB$). Autrement dit, le module des \diles est nul \ssi la transposée de la matrice jacobienne est surjective. 
Il est clair qu'une \alg étale sur un anneau commutatif arbitraire est nette.

On a le \tho important suivant qui donne une réciproque très forte du  
\thref{fthstrucste}. Pour la notion de zéro isolé simple dans le \tho on peut voir le traitement \cof dans \citet*[section IX-4]{fCACM}.

\begin{ftheorem}[\alg nette sur un \cdi] \label{fthNewtonCorpsDiscret}~\\
Sur un \cdi $\gK$ toute \Klg \pf nette est \stfe, étale, \ste.  
En particulier, pour un \syp étale, tous les zéros de la variété correspondante  sur un surcorps \ac sont isolés, simples et leurs \coos sont des \elts \spls sur $\gK$.
\end{ftheorem}
%

Une démonstration constructive se trouve dans \citet*[Corollary~\hbox{VI-6.15}]{fCACM}. Dans la deuxième édition française \citealt*{fACMC}, une démonstration \cov plus élémentaire est proposée à la fin de la section VI-6.

\smallskip Donnons maintenant une description plus précise de la situation. 

\newcommand{\Kual}{\gK[\ual]}
\begin{fdescri}[\syp étale sur un \cdi, précisions] \label{fdescriEtaleCdi} ~\\
Soit 
 $(\uf)= (f_1,\dots,f_n)$ un \syp étale de $n$ \pols dans  $\KXn$ sur un \cdi infini $\gK$\footnote{Lorsque $\gK$ est fini, ou plus \gnlt lorsqu'on ne sait pas s'il est infini, de légères modifications doivent être introduites dans cette description en tenant compte du point \textsl{1d)} du \thref{fthstrucste}.}  
 et notons 
 \[\gD=\Aqo\KuX \uf=\Kxn\] 
 l'\alg quotient sur $\gK$. 
\begin{itemize}
\item On sait construire un \elt primitif $\zeta$ de $\gD$, son \polmin $g$ sur $\gK$ est \spl, de sorte que $\gD=\gK[z]\simeq \aqo{\KZ}{g(Z)}$. 
On a donc des \pols $q_i\in\KZ$ tels que $x_i=q_i(z)$ dans~$\gD$.
\item  Un zéro $(\ual)=(\alpha_1,\dots,\alpha_n)$ du \syp dans une \Klg arbitraire $\gC$ donne un $\gK$-morphisme $\varphi\colon\gD\to\Kual \subseteq \gC$ vérifiant $\varphi(\ux)=(\ual)$.
L'\alg $\Kual$ est alors isomorphe à
un quotient de $\gD$.
Si $\Kual$ est connexe non triviale, c'est un \cdi, car il est \zedr \citep[Fact IV-8.8]{fCACM}.
\item  Nous supposons dans la suite que $\Kual$ est connexe non triviale.
\begin{itemize}
\item [$\bullet$] Si l'on connait une \fcn de $g$ en un produit de $r$ \pols~$g_j$ \irds sur~$\gK$, alors $\gD\simeq \gL_1\times \dots\times \gL_r$ avec $\gL_j\simeq \aqo \gD{g_j(z)}$, et l'on a un \sfio $(e_1(z),\dots,e_r(z))$ correspondant dans $\gD$\footnote{L'\id $\gen{g_j(z)}$ est engendré par l'\idm $1-e_j(z)$. Notons que le cas avec $r=0$ n'est pas exclu. Un \sype peut être impossible.}. Et $\gK[\ual]$ est isomorphe à l'un des \cdis $\gL_j$ via $\varphi$, avec $\varphi(e_j(z))=1$.
\item [$\bullet$] Par ailleurs si $\gL$ est une extension \spb de $\gK$
sur laquelle $g$ se factorise complètement,
l'\alg quotient vue sur~$\gL$, i.e. $\gL\otimes_\gK\gD\simeq\aqo{\gL[Z]}{g(Z)} $, est isomorphe  à~$\gL^{d}$, \hbox{où $d=\deg(g)$}. Ainsi le \syp considéré possède exactement $d$ zéros à \coos dans $\gL$.
Si l'on immerge $\gL\subseteq\Ksep\subseteq\Kac$ on a ainsi tous les zéros à \coos dans~$\Kac$\footnote{Ici \(\Ksep\) et \(\Kac\) sont respectivement une \clsep et une \cla de~\(\gK\).}. \eoe
\end{itemize}
\end{itemize}
\end{fdescri}

Que se passe-t-il dans la situation plus \gnle où l'on ne sait pas factoriser $g$ sur~$\gK$? Ce cas n'est pas vraiment différent du précédent. En effet, si à un moment donné, on a mis en évidence un \idm $e\neq 0,1$ de $\gD$\footnote{Cela se produit chaque fois qu'un \elt $\neq 0$ de $\gD$ n'est pas inversible, \cad lorsque son \polmin est de degré $>1$ et a son \coe constant nul.}, il est égal à $1$ ou $0$ dans $\Kual$ et l'on peut remplacer~$\gD$ par $\aqo\gD{1-e}$ ou par $\aqo\gD{e}$ (ce qui revient à remplacer $g$ par un diviseur strict). La nouvelle \alg est encore étale (c'est une localisée de $\gD$). Tous les calculs précédents restent valables, et la nouvelle version de $g$ est meilleure. Le nombre d'améliorations possibles est clairement borné par $\deg(g)$. En conséquence la plupart des résultats concrets obtenus en supposant que l'on sait factoriser $g$ resteront valables sans cette hypothèse. C'est ici l'essence même de la méthode dynamique qui entre en jeu.    

\section{Le cas des \cvdus}\label{fseccvdu}

\subsection{Définitions}

\Subsubsection{Corps valués discrets}

Rappelons d'abord la \dfn \cov d'un \textsl{\cvd} $\KV$ comme dans \cite{fKL00} ou \citet*{fCLR01}:
\begin{itemize}
\item $\gK$ est un \cdi;
\item $\gV$ est un sous-anneau de $\gK$; 
\item pour tout $x\in\gK\eti$, on a  $x$ ou $1/x\in\gV$; 
\item la \dve dans $\gV$ est explicite%
\footnote{Cela signifie que pour $x,y\in\gV$ on a un test pour l'existence d'un $z\in\gV$ tel que $yz=x$. Cela revient à dire que  $\gV$ est un sous-anneau détachable de $\gK$. C'est toujours vrai en \clama par la loi du tiers exclu.}. 
\end{itemize}

Dans ce cas~$\gV$ est un \alrd   \icl. 

Cette \dfn équivaut en \clama à la \dfn usuelle. 
On a ajouté des hypothèses de décidabilité pour faciliter les calculs.

\begin{fdefinition}[hensélisé d'un \cvd]
Soit $\KV$ un \cvd.  
\begin{itemize}

\item  Une \textsl{extension} de $\KV$ est un \cvd $(\gL,\gW)$ avec un \homo $\phi\colon\gK \rightarrow \gL$ qui vérifie $\gV = \gK \cap  \phi^{-1}(\gW)$.

\item  Un \textsl{hensélisé} de $\KV $ est une extension $(\KHe,\VHe)$
qui est un \cvd hensélien, 
tel que le morphisme
$\phi^H\colon\gK \rightarrow \KHe$ se factorise de manière unique à travers toute extension de $(\gK, \gV)$ qui est un \cvd hensélien.
\end{itemize}

\end{fdefinition}

En tant que solution d'un \pb \uvl, le hensélisé d'un \cvd est unique~à \iso unique près.
 
\cite{fKL00} construisent un hensélisé $(\KHe,\VHe)$ d'un \cvd $\KV$. Il est obtenu en ajoutant succcessivement des \zeHs de \polHs.
 
Pour un \pol de Hensel $ f = X^n + a_{n-1}X^{n-1} + \cdots + a_0 $ dans $\gV[X]$, les auteurs décrivent explicitement une extension $ (\gK[\alpha], \gV_{\alpha})$
de $\KV$ pour laquelle $f$ vu dans $\gK [\alpha][X]$
a un \zeH $\alpha$, et telle que le morphisme naturel $\gK \rightarrow  \gK[\alpha]$ factorise de manière unique toute extension  $(\gL,\gW)$ telle que  $f$
vu dans $\gL[X]$ a un \zeH. En outre le corps résiduel et le groupe de valeurs de $(\gK[\alpha], \gV_{\alpha}) $  sont canoniquement isomorphes au corps résiduel et au groupe de valeurs de  $\KV$. Cette construction explicite est basée sur l’étude du \pgn de~$f$.
Notons que l'on ne suppose pas savoir si le corps de départ contient un zéro spécial de~$f$. 
Par suite, l'extension finie $\gK[\alpha]$ que l'on construit est un \cdi mais elle n'est pas \ncrt une extension \stfe de $\gK$.  
La construction du hensélisé est donc tout à fait similaire à la construction de la clôture réelle d'un \codi,
qui fonctionne même lorsque l'on ne sait pas décider si un \pol arbitraire admet un zéro réel ou pas dans le corps de départ. 

\begin{fremark}
La construction du hensélisé dans \cite{fKL00} se situe entièrement dans le cadre des \cvds, et elle est à priori moins \gnle que la construction dans le cadre des \alrds donnée dans \cite{fALP08}. 
Les outils utilisés dans ces deux constructions sont très différents.
Cela justifie que nous prenions deux notations distinctes: $\Vhe$ pour l’henselisation en tant qu’\alo, et $\VHe$ en tant qu’\adv du hensélisé du \cvd. Ce n’est pas évident que ces deux hensélisés coïncident, mais nous n’avons pas trouvé de \demo de cette coïncidence dans la littérature. 
\eoe
\end{fremark}

\Subsubsection{Corps ultramétriques}

Dans l'ouvrage \citealt*{fMRR}, la théorie des \valas est développée \cot en utilisant la \dfn suivante, qui est la \dfn usuelle des corps munis d'une valeur absolue en \clama\footnote{En fait, ils partent avec un corps de Heyting $\gK$, mais les deux \dfns sont clairement \eqves. Et en \clama on utilise aussi habituellement un corps $\gK$ au départ.}. 

\begin{fdefinition}[\cvar, \cvu] \label{fdefivala}~
\begin{enumerate}
\item 
Une \textsl{\vala} sur un anneau $\gK$ est une fonction
$\gK\to \RR^{\geq 0}, \,x \mt |x|$ qui satisfait les \prts suivantes:
\begin{itemize}
\item $|x| = 0$ \ssi $x=0$;
\item $|x| > 0$ \ssi $x$ est inversible;
\item $|xy| = |x||y|$;
\item $|x + y| \leq |x| + |y|$.
\end{itemize}
 On dit alors que $(\gK,\abs\cdot)$ est un \textsl{corps muni d'une valeur absolue}, ou un \textsl{\cvar}\footnote{L'ouvrage \citealt*{fMRR} utilise le terme \gui{\cval}, comme souvent dans la littérature anglaise en théorie des nombres, mais cela entre en conflit avec la terminologie que nous adoptons ici pour les corps valués généraux à la suite de Krull et Bourbaki.}.
%
%
\item 
La \vala est dite \textsl{ultramétrique} si $|x + y| \leq \sup(|x|,|y|)$\footnote{Notez que d'un point de vue \cof, $z=\sup(x,y)$ est bien défini pour des nombres réels $x,y$, mais on ne peut par prouver que $z=x$ ou $z=y$ avec un \gui{ou} explicite.}.
Dans ce cas on parle de \textsl{\cvu} ou plus simplement de \textsl{corps \ultm}. On note que si $\abs x<\abs y$ alors $\abs{x+y}=\abs y$.
\item 
Toute \valu définit un \gui{\adv} $\gV$: 
\[
\formule{
\gV&:=&\sotq{x\in\gK}{|x|\leq 1} ,\hbox{ avec} 
\\ \gV\eti&:=&\sotq{x\in\gK}{|x|= 1}\hbox{ et }
\\ \fm_\gV&:=&\Rad(\gV)=\sotq{x\in\gK}{|x|< 1}.
}
\]  
\item 
Deux \valus non triviales sur $\gK$ qui ont le même \adv sont dites \textsl{\eqves}, et chacune est une puissance strictement positive de l'autre \citep[voir][Theorem XII-1.2]{fMRR}. 
\item 
La distance $d(x,y)=|x-y|$ fait de $\gK$ un espace métrique, dont le complété est noté $\gKh$. La \vala s'étend de manière unique à $\gKh$, et $(\gKh,\abs{\cdot})$ est aussi un corps \ultm\footnote{Pour le cas d'une \vala archimédienne, $(\gKh,\abs{\cdot})$ est aussi un corps valorisé.}.  Les complétés de $\gV$ et $\fm_\gV$ sont notés $\gVh$ et $\fmVh$.
L'image de $\gK\eti$ dans $(\RR^{>0},\times )$ est le \textsl{groupe des valeurs} de $(\gK,\abs\cdot)$.

\end{enumerate}
\end{fdefinition}

Un \cvar $\gK$ est un \alo non trivial dont le radical de Jacobson est réduit~à~$0$, \cad un \textsl{corps de Heyting}%
\footnote{Comme $\gK$ n'est pas \ncrt un \cdi, nous avons préféré introduire la \dfn pour un anneau $\gK$, sans avoir à rappeler avant de commencer ce qu'est un corps de Heyting.}.
 
La \dfn ne nécessite pas que $\gK$ soit un \cdi. En général le complété~$\gKh$ n'est pas discret, même lorsque $\gK$ est discret.

\begin{fremark} \label{fremdefivala} 
Dans le point 4. nous avons mis \gui{\adv} entre guillemets car l'anneau $\gV$ n'est pas \ncrt, du point de vue \cof, ni local, ni discret, ni \dcd.
\end{fremark}

\smallskip Considérons un \cvu $(\gK,\abs\cdot)$

Le corps de Heyting $\gK$ est discret \ssi pour tout $x\in\gK$ on a la disjonction \hbox{\gui{$\abs x = 0 \vuu \abs x >0$}}. Cela revient à dire que $\gV$ est intègre. 

Le couple $\KV$ est un \cvd au sens \cof si l'on a en outre la disjonction \gui{$\abs x = 1 \vuu \abs x <1 \vuu \abs x >1$} pour tout $x\in\gK$.

Cela revient à dire que $\gV$ est un \alo intègre avec un \crdl discret.
On dit dans ce cas que $(\gK,\abs\cdot)$ est un \textsl{\cvdu}, ou un \textsl{\cud}.
  
\Subsubsection{Traduction en termes de valuations}

Considérons l’application \(\ell:(\RR^{\geq 0},\times)\to(\RR\cup\sing{+\infty},+)\) définie par \(\ell(r)=-\log(r)\) pour $r\neq 0$ et par $\ell(0)=+\infty$; munissons $\RR\cup\sing{+\infty}$ de la toplogie qui fait de  $\ell$ un homéomorphisme.
On définit alors la valuation $v: \gK\to (\RR\cup\sing{+\infty},+)$ par $v(x)=\ell(\abs x)$ 

 On est donc ici simplement en train de traduire les \prts de $(x\mapsto \abs x,\,\gK\to \RR^{\geq 0})$  en les \prts de $v$, en renversant l'ordre, en remplaçant la multiplication par l'addition et en remplaçant $\sup$ par $\inf$.
On a donc les \prts suivantes:
\begin{itemize}
\item $ v(x) = \infty$ \ssi $x=0$;
\item $v(x) \neq  \infty$ \ssi $x$ est inversible;
\item $v(xy) = v(x)+v(y)$;
\item $v(x + y) \geq \inf(v(x),v(y))$, avec \egt si $v(x)\neq v(y)$;
\item $ v(x) \geq 0$ \ssi $x\in \gV$;
\item $ v(x) > 0$ \ssi $x\in \fm_\gV$;
%
%
\end{itemize}

Si l'on remplace la \vala par une \vala \eqve, la valuation $v$ est simplement multipliée par une constante $r>0$.

\smallskip Lorsque l'on a un \cud, $\gK$ et l'anneau résiduel $\gV/\fm$ sont des \cdis, le sous-groupe $\sotq{\abs x}{x\in\gK\eti}$ est un sous-groupe multiplicatif discret $\Delta$ de  $\RR^{\geq 0}$, \hbox{$\Gamma:=\sotq{v (x)}{x\in\gK\eti}$} est un sous-groupe additif discret de $(\RR,+)$ et la réunion
$\Gamma\cup\sing{+\infty}$ est une réunion disjointe: la topologie est discrète. 

Dans ce cas, on a le résultat suivant très utile valable pour les \cvds en \gnl:
\begin{itemize}
\item si $\sum_{i=1}^nx_i=0$, avec les $x_i$ non tous nuls, l'inf des $v(x_i)$ est atteint au moins pour deux indices $i$ distincts.
\end{itemize}

\Subsubsection{Trois exemples de base}

En \coma on définit un  \textsl{\advd}  comme un anneau intègre  $\gV$ (de \cdf $\gK$) dans lequel on donne un \elt  $\pi\neq 0$ (appelé une \textsl{uniformisante}) tel que tout \elt de $\Vtl$ s'écrive sous forme  $a=u\pi^k$ avec $u\in\gV\eti$ et $k\in\NN$. 
Cela fait de $\KV$ un \cvd avec la \valn $v(u\pi^k)=k$. En posant $\abs {u\pi^k}= e^{-k}$ pour un réel $e>0$, $(\gK,\abs\cdot)$ est un \cud.   

Trois exemples de base sont donnés comme suit. Dans les exemples 2 et 3, la \vala n'est pas à valeurs dans $\RR^{\geq 0}$ mais dans un sous-\mo de 
$(\gK,\times)$ isomorphe~à l'adhérence de $\sotq{x=1/2^n}{n\in\NN}$ 
dans $\RR^{\geq 0}$. 

\begin{enumerate}
\item Ici $\gK=\QQ$, la \vala $p$-adique standard est $|r|_p=p^{-k}$ pour~$r=\frac m n\,p^k$ avec $m$ et $n\in\ZZ$ étrangers à $p$. L'\advd correspondant est $\gV=\ZZ_{1+p\ZZ}$ (le localisé de $\ZZ$ en $\gen{p}$) avec $\Rad\gV=p\gV$, l'uniformisante $p$ et le \crdl~$\Fp$.
\item Ici $\gK=\QQ(t)$, la \gui{\vala} $t$-adique standard est $|r|_t=t^{-k}$ pour~$r=\frac m n\,t^k$ \hbox{avec $m,n\in\QQ[t]$}, $m(0)$ et $n(0)\neq 0$ dans $\QQ$.
L'\advd correspondant est \hbox{$\gV=(\QQ[t])_{1+t\QQ[t]}$}  (le localisé de $\QQ[t]$ en $\gen{t}$) avec $\Rad\gV=t\gV$, l'uniformisante $t$ et le \crdl~$\QQ$. 
\item Ici $\gK=\Fp(t)$, la \gui{\vala} $t$-adique standard est $|r|_t=t^{-k}$ pour~$r=\frac m n\,t^k$ avec avec $m,n\in \Fp[t]$, 
$m(0)$ et $n(0)\neq 0$ dans $\Fp$.
L'\advd  est \hbox{$\gV=(\Fp[t])_{1+t\Fp[t]}$}
  (le localisé de $\Fp[t]$ en $\gen{t}$) avec $\Rad\gV=t\gV$, l'uniformisante~$t$ et le \crdl~$\Fp$. 
\end{enumerate}

Dans le deuxième exemple, le corps $\QQ$ peut être remplacé par un \cdi arbitraire (comme dans le troisième exemple).

\subsection{Le \lHm pour les \cuds}

\Subsubsection{Un résultat essentiel dans \citet*{fMRR}}

\begin{fnotation} \label{fnotaKtilde}
Si $(\gK,\abs\cdot)$ est un \cvu avec une \vala non triviale\footnote{\Cad qu'il existe un $x$ tel que $\abs x\neq 0,1$.} on note 
\begin{itemize}
\item $\gKt$ la clôture séparable de~$\gK$ dans son complété $\gKh$,
\item $\gVt=\sootq{x\in\gKt}{\abs x \leq 1}$ son \gui{anneau de valuation},   et
\item $\fmti=\sootq{x\in\gKt}{\abs x < 1}$.
\end{itemize}
\end{fnotation}

\citet*{fMRR} démontrent que pour un \cud $(\gK,\abs \cdot)$,  $(\gKt,\gVt)$ est un \cvdh au sens usuel (tout \polH admet un \zeH).

La \demo de ce résultat est assez compliquée parce que \citet{fMRR} ont en vue des résultats \gnls concernant le cas non discret.
Nous proposons en conséquence une \demo plus simple dans le paragraphe qui suit.
 
\Subsubsection{Une \demo plus simple}

Nous pourrions donner cette \demo dans le cadre des \cvds dont le groupe de valeurs $\Gamma$ est archimédien. Dans le cadre d'un \cvd, le groupe de valeurs est totalement ordonné discret, et s'il est archimédien il est isomorphe à un sous-groupe de $(\RR,+)$, la situation est donc la même que pour un \cud. 

Nous démontrons deux lemmes.

\begin{flemma} \label{flemComplete}
Soit $(\gK,\abs\cdot)$ un \cud. 
\begin{enumerate}
\item Alors  $(\gKh,\abs\cdot)$ est un \cvu, $\gVh$ est un \alo de radical $\fmVh$, et l'anneau résiduel $\gVh/\fmVh$  est isomorphe à $\gV/\fm$ (c'est un \cdi).
\item L'\alo  $(\gVh,\fmVh)$ est hensélien. Plus \gnlt tout \sysN $(\lfn)\in\gVh[\uX]^n$ au point $(\an)\in {\gVh}^n$ possède un zéro $(\xin)$
avec les $\xi_k\in a_k+\fmVh$.
%
%
\end{enumerate}
\end{flemma}
%
\begin{proof} 
Le premier point est facile. En particulier $\gKh$ est un corps de Heyting, donc un \alo. Le deuxième point est le \crl \ref{fLHM1}, une conséquence de la méthode de Newton expliquée en \ref{fthNewtonQuad}. 
\end{proof}
%

\begin{flemma} \label{flemKtilde} 

\noindent 
Les \elts de $\gKt$ forment un sous-anneau discret de $\gKh$ et $(\gKt,\abs\cdot)$ est un \cud hensélien. 
\end{flemma}
%
\begin{proof} Le fait que $\gKt$ est un sous-anneau est classique.
Soit maintenant  un \elt~\hbox{$\xi\in\gKh$} qui annule un \pol \spl $f\in\KX$. 
Posons $\Kx=\aqo\KX f$. On~a un morphisme
de \Klgs, 
$\varphi\colon\Kx \to \gK[\xi]$ tel que $\varphi(x)=\xi$. 
Comme~$\gK[\xi]$ est connexe (c'est un sous-anneau de l'\alo $\gKh$) 
c'est un \cdi en tant que quotient d'une \Klg \spb \stfe: un anneau \zedr connexe est un \cdi \citep[Fact IV-8.8]{fCACM}.  Ainsi $\gKt$ est un \cdi. 
\end{proof}
%

\Subsubsection{L'\iso entre deux variations sur le thème du hensélisé}

Nous commençons par rappeler une  variante de style Hensel-Newton du lemme de Hensel pour les \pols univariés \cite[Proposition XII-7.6]{fLang}.
Cela s'applique pour tous les \cvdhs.

\begin{flemma}[lemme de Hensel-Newton pour un \cvd] \label{flemNewtonHensel}~\\
Soit $\KV$ un \cvd. On note $v:\gK\to\Gamma\cup\so{+\infty}$ la valuation associée. \\
Soit $F(x)=\sum_{k=0}^na_kx^k\in\Vx$ avec $a_1\neq 0$ et  $v(a_0)>2v(a_1)$. 
\begin{enumerate}
\item Le \pol $F$ admet un zéro dans $\frac{a_0}{a_1}\cdot\VHe\subseteq a_{1}\cdot\fm \VHe$, et c'est l'unique zéro de $F$ dans $a_{1}\cdot\fm \VHe$.
\item En particulier si $\KV$ est hensélien, $F$ admet un zéro $\xi\in\frac{a_0}{a_1}\cdot\gV\subseteq a_{1}\cdot\fm \gV$, et c'est l'unique zéro de $F$ dans $a_{1}\cdot\fm \gV$
\end{enumerate}
\end{flemma}
%
\begin{proof}
On considère le \pol 
\[
f(x)=\frac{1}{a_1^2}\,F(a_1x)=\frac{a_0}{a_1^2}+X+\som_{k=2}^na_ka_1^{k-2}X^k.
\]
Les hypothèses du lemme~\ref{flemTrick1} sont satisfaites, donc le \pol $f$ a un zéro~$\zeta\in\frac{a_0}{a_1^2}\cdot\VHe\subseteq \fm\VHe$, et c'est
l'unique zéro de $f$ dans $\fm\VHe$. Cela donne le zéro $\xi=a_1\zeta$ de $F$ dans 
$\frac{a_0}{a_1}\cdot\VHe\subseteq a_1\cdot\fm\VHe$, 
et c'est l'unique zéro de $F$ dans $a_1\cdot\fm \VHe$.    
\end{proof}

En relisant la \demo du lemme \ref{flemTrick1}, on voit que le calcul qui a été fait montre que~$\zeta$ est  dans l'image d'un étage initial de la construction du hensélisé $\Vhe$, obtenu en ajoutant le zéro spécial d'un \pol spécial, sans que $\zeta$ lui-même soit le \zeH d'un \polH. En outre $v(\zeta)>0$.     

\begin{ftheorem} [deux versions \eqves du hensélisé d'un \cud] 
\label{fthMRRKL}~\\
On considère un  \cud $(\gK,\abs\cdot)$.
Le hensélisé  $(\KHe,\VHe)$ de $\KV$ construit comme dans \citealt*{fKL00} est isomorphe à  $(\gKt,\gVt)$. Plus \prmt: 
\begin{enumerate}
\item il existe un unique $\gK$-\homo  $\KHe\to\gKt$ qui envoie $\VHe$ dans $\gVt$, et cet \homo est un  \iso;
\item le morphisme naturel $\Vhe\to\gVt$ est surjectif et le morphisme $\VHe\to\gVt$ est un \iso.
\end{enumerate}
\end{ftheorem}
\begin{proof}
Comme $(\gKt,\gVt)$ est un \cvdh extension de~$\KV$, on a un unique  $\KV$-morphisme $\varphi\colon(\KHe,\VHe)\to(\gKt,\gVt)$. Il est injectif parce que~$\KHe$ est un \cdi et~$\gKt$ n'est pas trivial.
Nous devons montrer qu'il est surjectif. 

Il suffit de montrer le point \textsl{2}.

On considère un \elt $\xi\in\gVt$: $\xi\in\gKh$, $v(\xi)\geq 0$ et $f(\xi)=0$ pour un \pol $f\in\VX$ \spl dans $\KX$. En particulier $f'(\xi)\neq 0$, i.e. $v(f'(\xi))<+\infty$.

Puisque $\xi\in\gKh$ on connait des approximations arbitrairement précises de $\xi$ dans $\gK$, i.e. de la forme  \hbox{$a=\xi+\zeta\in\gK$}, avec $v(\zeta)$ arbitrairement grand. Puisque $v(\xi)\geq 0$, on a $v(a)\geq 0$.

On considère alors pour un tel $a$ le \pol 
$F_a(X)=f(X+a)\in\VX$.
Le \coe $c_0=F_a(0)=f(a)$ est une approximation arbitrairement précise de $f(\xi)=0$, i.e. $v(c_0)$ est arbitrairement grand. 
Le \coe $c_1=F’_a(0)=f'(a)$ est une approximation arbitrairement précise de~$f'(\xi)$, donc pour $v(\zeta)$ suffisamment grand $v(c_1)=v(f'(\xi))<+\infty$. Pour $v(\zeta)$ suffisamment grand on a donc $v(c_0)>2v(c_1)$. Le lemme \ref{flemNewtonHensel} nous dit que $F_a$ admet un zéro  $\alpha\in\frac{c_0}{c_1}\cdot\VHe\subseteq\fm \VHe$. 

L'image de $\alpha$ dans $\gKt$ est égale à $\zeta$ car ce sont deux \zeHs dans $\gKt$ du même \cdH $(F_a,0)$.
Donc l'image  de $a-\alpha$ est $a-\zeta=\xi$.
En outre comme il est expliqué juste après le lemme \ref{flemNewtonHensel}, $\xi$ est non seulement dans l'image du morphisme naturel $\VHe\to\gVt$
mais aussi dans celle  du morphisme naturel $\Vhe\to\gVt$.    
\end{proof}

On a obtenu \prmt que tout \elt de $\gKt$ est l'image d'un  \elt $\gamma$
dans un corps $\gK[\xi]\subseteq \KHe$ où 
$\xi$ est le zéro spécial d'un \pol spécial.
Ce n'était pas clair à priori dans \cite{fKL00} (où le \cvd était arbitraire).

\Subsubsection{Le \lHm}

Enfin on a le résultat souhaité.

\begin{ftheorem}[\lHm  pour un \cud] \label{fLHM2} ~ \\
Soit $(\gK,\abs\cdot)$ un \cud et 
$(f_1,\dots,f_n)$ un \sysN sur $\KV$. Ce système admet un unique zéro à coordonnées dans $\fmti$.  Il admet aussi un unique zéro à coordonnées dans $\fm \VHe$. 
\end{ftheorem}
%
\begin{proof}
On peut supposer que le \sysN est étale (note \ref{fnote0}). 
La méthode de Newton construit un zéro $(\ual)$ à \coos dans $\fmVh\subseteq \gKh$ (corolaire \ref{fLHM1}).
On note $\gD=\Kxn=\aqo\KXn{\lfn}$ l'\alg quotient sur $\gK$ du \syp.
Elle est \stfe, \ste.
On a le morphisme naturel de \Klgs $\varphi\colon\gD\to\gK[\ual]\subseteq \gKh$, où $(\ual)$ est le \zeH du \syp, ($\varphi(x_i)=\alpha_i$
pour chaque $i$).  
Donc les~$\alpha_i$ sont \spls sur $\gK$.  Ainsi pour chaque $i$,  $\alpha_i\in\fmti$: les \coos du \zeH $(\ual)$ sont  dans $\fmti$. 
\\
Par ailleurs $\gVt$ est canoniquement isomorphe à~$\VHe$ (\thref{fthMRRKL}). Les \elts  dans $\VHe$ qui correspondent aux $\alpha_i$  donnent le \zeH du \syp à \coos dans~$\mHe$. 
\end{proof}

Notons que d'après la description \ref{fdescriEtaleCdi},
comme la \Klg $\gK[\ual]$ est connexe non triviale\footnote{En tant que sous-\alg de $\gKh$ qui est un \cvar, donc un \alo.}, $\gK[\ual]$ est un \cdi isomorphe à un quotient de $\gD$, mais il ne semble pas qu'on puisse
déterminer à coup sûr la dimension de  $\gK[\ual]$ comme \Kev. 

\medskip 
\centerline{--------------------------}

\medskip \noindent \textsl{Remarque finale.}\\ 
Notez que de nombreux textes de classiques qui utilisent ou proposent des variantes du \lHm,
comme par exemple  \cite{fFis1997} et \cite{fSma1998} 
donnent la solution sous forme de zéros dont les \coos sont des \elts d'une certaine \gui{completion} du corps valué considéré et non dans son hensélisé.
\\
Les articles \cite{fKuh2011} et \cite{fPCR2000} donnent la solution, de manière non \algq,
en utilisant la notion de corps sphériquement complet. Leurs \demos sont pour le moment difficiles à interpréter \cot.


\newpage
\addcontentsline{toc}{section}{Références}
\rdb
\small
\bibliographystyle{plainnat-fr}

\end{document}